\documentstyle[11pt]{article}

\newtheorem{eg}{{\sc Example}}
\newtheorem{lemma}{{\sc Lemma}}[section]
\newtheorem{theo}{{\sc Theorem}}
\newtheorem{prop}{{\sc Proposition}}
\newtheorem{cor}[theo]{{\sc Corollary}}

\font\bbb=msbm10 scaled\magstep1

\newcommand{\RR}{\mbox{\bbb R}}
\newcommand{\ZZ}{\mbox{\bbb Z}}

\newcommand{\si}{\sigma}

\begin{document}

\title{Degree-regular triangulations of torus and Klein bottle}

\author{
{Basudeb Datta, Ashish Kumar Upadhyay}\\[2mm]
{\normalsize Department of Mathematics, Indian Institute of Science, }  \\
{\normalsize Bangalore 560\,012,  India.}\\[1mm]
{\small dattab@math.iisc.ernet.in, upadhyay@math.iisc.ernet.in} }
\date{To appear in `{\bf Proc. Indian Acad. Sci. $($Math. Sci.$)$}'}
\maketitle

\vspace{-5mm}

\hrule

\begin{abstract} A triangulation of a connected closed surface is
called weakly regular if the action of its automorphism group on
its vertices is transitive. A triangulation of a connected closed
surface is called degree-regular if each of its vertices have the
same degree. Clearly, a weakly regular triangulation is
degree-regular. In \cite{l}, Lutz has classified all the weakly
regular triangulations on at most 15 vertices. In \cite{dn},
Datta and Nilakantan have classified all the degree-regular
triangulations of closed surfaces on at most 11 vertices.

In this article, we have proved that any degree-regular
triangulation of the torus is weakly regular. We have shown that
there exists an $n$-vertex degree-regular triangulation of the
Klein bottle if and only if $n$ is a composite number $\geq 9$.
We have constructed two distinct $n$-vertex weakly regular
triangulations of the torus for each $n \geq 12$ and a $(4m +
2)$-vertex weakly regular triangulation of the Klein bottle for
each $m \geq 2$. For $12 \leq n \leq 15$, we have classified all
the $n$-vertex degree-regular triangulations of the torus and the
Klein bottle. There are exactly 19 such triangulations, 12 of
which are triangulations of the torus and remaining 7 are
triangulations of the Klein bottle. Among the last $7$, only one
is weakly regular.
\end{abstract}

{\small

{\bf AMS classification\,:} 57Q15, 57M20, 57N05.

{\bf Keywords\,:} Triangulations of $2$-manifolds, regular
simplicial maps, combinatorially  \newline \mbox{} \hspace{25mm}
regular triangulations, degree-regular triangulations.

}

\bigskip

\hrule

\section{Introduction and results.}

Recall that a {\em simplicial complex} is a collection of
non-empty finite sets (set of {\em vertices}) such that every
non-empty subset of an element is also an element.  For $i \geq
0$, the elements of size $i+1$ are called the {\em $i$-simplices}
of the simplicial complex. $1$-simplices are also called the {\em
edges} of the simplicial complex. For a simplicial complex $X$,
the maximum of $k$ such that $X$ has a $k$-simplex is called the
{\em dimension} of $X$. The set $V(X)$ of vertices of $X$ is
called the {\em vertex-set} of $X$. A simplicial complex $X$ is
called finite if $V(X)$ is finite.

If $X$ and $Y$ are two simplicial complexes, then a (simplicial)
isomorphism from $X$ to $Y$ is a bijection $\varphi :  V(X) \to
V(Y)$ such that for $\si \subseteq V(X)$, $\si$ is a simplex of
$X$ if and only if $\varphi(\si)$ is a simplex of $Y$. Two
simplicial  complexes $X, Y$ are called (simplicially) {\em
isomorphic} (and is denoted by $X \cong Y$) when such an
isomorphism exists. We identify two complexes if they are
isomorphic. An isomorphism from a simplicial complex $X$ to
itself is called an {\em automorphism} of $X$. All the
automorphisms of $X$ form a group, which is denoted by ${\rm
Aut}(X)$.

A simplicial complex $X$ is usually thought of as a prescription
for constructing a topological space (called the {\em geometric
carrier} of $X$ and is denoted by $|X|$) by pasting together
geometric simplices. Formally, $|X|$ is the subspace of
$[0,1]^{V(X)}$ consisting of the functions $f: V(X) \to [0,1]$
such that the support $\{ v \in V (X) : f (v) \neq 0 \}$ is a
simplex of $X$ and $\sum_{v \in V (X)} f(v) = 1$.  If $\sigma $
is a simplex then $|\sigma  |:= \{ f\in |X| : \sum_{v \in \sigma }
f(v) = 1 \}$ is called the {\em geometric carrier} of $\sigma $.
We say that a simplicial complex $X$ {\em triangulates} a
topological space $P$ (or $X$ is a {\em triangulation} of $P$) if
$P$ is homeomorphic to $|X|$. A simplicial complex $X$ is called
{\em connected} if $|X|$ is connected. A 2-dimensional simplicial
complex is called a {\em combinatorial $2$-manifold} if it
triangulates a closed surface. A combinatorial $2$-manifold $X$
is called {\em orientable} if $|X|$ is an orientable 2-manifold.

If $v$ is a vertex of a simplicial complex $X$, then the number
of edges containing $v$ is called the {\em degree} of $v$ and is
denoted by $\deg_X(v)$ (or $\deg(v)$). If the number of
$i$-simplices of an $m$-dimensional finite simplicial complex $X$
is $f_i(X)$ ($0\leq i\leq m$), then the number $\chi(X) :=
\sum_{i= 0}^{m}(- 1)^i f_i(X)$ is called the {\em Euler
characteristic} of $X$. A simplicial complex is called {\em
neighbourly} if each pair of vertices form an edge.

A {\em combinatorially regular} combinatorial 2-manifold is a
connected combinatorial 2-manifold with a flag-transitive
automorphism group (a {\em flag} is a triple $(u, e, F)$, where
$e$ is an edge of the face $F$ and $u$ is a vertex of $e$). A
connected combinatorial $2$-manifold $X$ is said to be {\em
weakly regular} (or a {\em weakly regular triangulation} of
$|X|$) if the automorphism group of $X$ acts transitively on
$V(X)$. Clearly, a combinatorially regular combinatorial
$2$-manifold is weakly regular. Well known examples of
combinatorially regular combinatorial $2$-manifolds are the
boundaries of the tetrahedron, the octahedron, the icosahedron and
the 6-vertex real projective plane (\cite{bd, dn}). The
combinatorial manifolds $T_{3, 3, 0}$ and $T_{6, 2, 2}$ (in
Examples 2 \& 3) are combinatorially regular. Schulte and Wills
(\cite{sw1, sw2}) have constructed two combinatorially regular
triangulations of the orientable surface of genus 3. In \cite{l},
Lutz has shown that there are exactly $14$ combinatorially regular
combinatorial $2$-manifolds on at most $22$ vertices. By using
computer, Lutz has shown the following\,:

\begin{prop}$\!\!${\bf .} \label{lutz1}
There are exactly $77$ weakly regular combinatorial $2$-manifolds
on at most $15$ vertices; $42$ of these are orientable and $35$
are non-orientable. Among these $77$ combinatorial $2$-manifolds,
$20$ are of Euler characteristic $0$. These $20$ are $T_{7, 1, 2},
\dots, T_{15, 1, 2}$, $T_{12, 1, 3}, \dots, T_{15, 1, 3}$, $T_{12,
1, 4}, T_{15, 1, 4}$, $T_{15, 1, 5}$, $T_{6, 2, 2}$, $T_{3, 3,
0}$, $Q_{5, 2}$ and $Q_{7, 2}$ of Examples $1, 2, 3, 6$.
\end{prop}

A connected combinatorial $2$-manifold $X$ is said to be {\em
degree-regular of type $d$} if each vertex of $X$ has degree $d$.
A combinatorial $2$-manifold $X$ is said to be {\em
degree-regular} (or a {\em degree-regular triangulation} of
$|X|$) if it is degree-regular of type $d$ for some $d$. So,
trivial examples of degree-regular combinatorial $2$-manifolds
are weakly regular and neighbourly combinatorial $2$-manifolds.

If $K$ is an $n$-vertex degree-regular of type $d$ combinatorial
2-manifold then $nd = 2f_1(K)$ $= 3f_2(K)$ and $\chi(K) = f_0(K)
- f_1(K) + f_2(K) = n - \frac{nd}{2} + \frac{nd}{3} = \frac{n(6 -
d)}{6}$. So, if $\chi(K) \neq 0$ then only finitely many $(n, d)$
satisfies the above equation and hence only finitely many
degree-regular combinatorial 2-manifolds of a given non-zero Euler
characteristic. If $K$ is degree-regular and $\chi(K) > 0$ then
$(n, d) = (4, 3), (6, 4), (6, 5)$ or $(12, 5)$. For each $(n, d)
\in \{(4, 3), (6, 4), (6, 5), (12, 5)\}$, there exists unique
combinatorial 2-manifold, namely, the 4-vertex 2-sphere, the
boundary of the octahedron, the 6-vertex real projective plane
and the boundary of the icosahedron (cf. \cite{bd, dn}). These 4
combinatorial 2-manifold are combinatorially regular. For the
existence of degree-regular of type $d$ combinatorial 2-manifolds
of negative Euler characteristic, $d$ must be at least 7. Since
$\frac{n(6-d)}{6} \neq -1$ for $n > d \geq 7$, there does not
exist any degree-regular combinatorial 2-manifolds of Euler
characteristic $-1$. If $\chi(K) = -2$ then $(f_0(K), d) = (12,
7)$. In \cite{du}, we have seen that there are exactly 6
degree-regular triangulations of the orientable surface of genus
2, three of which are weakly regular and none of them are
combinatorially regular.

For the existence of an $n$-vertex neighbourly combinatorial
2-manifold, $n(n-1)$ must be divisible by 6, equivalently,
$n\equiv 0$ or 1 mod 3. Ringel and Jungerman (\cite{jr, r}) have
shown that {\em there exist neighbourly combinatorial
$2$-manifolds on $3k$ and $3k+1$ vertices, for each $k\geq 2$}.
By using computer, Altshuler {\it et al} (\cite{abs}) have shown
that {\em there are exactly $59$ orientable neighbourly
combinatorial $2$-manifolds on $12$ vertices.} In \cite{a2},
Altshuler describe two operations by which one gets many
neighbourly combinatorial 2-manifolds from one such on same
number of vertices. Using this he has constructed 40615 distinct
non-orientable neighbourly combinatorial $2$-manifolds on $12$
vertices.

Here we are interested in the cases when the Euler characteristic
is 0 (i.e., triangulations of the torus and the Klein bottle).
Clearly, If $K$ is an $n$-vertex degree-regular combinatorial
2-manifold and $\chi(K) = 0$ then $n > d = 6$. From \cite{dn}, we
know the following\,:

\begin{prop}$\!\!${\bf .} \label{datta1}
$(a)$ For each $n\geq 7$, there exists an $n$-vertex weakly
regular triangulation of the torus. \vspace{1mm} \newline $(b)$
For each $k, l \geq 3$, there exists a $kl$-vertex degree-regular
triangulation of the Klein bottle.
\end{prop}

\begin{prop}$\!\!${\bf .} \label{datta2}
There are exactly $27$ degree-regular combinatorial $2$-manifolds
on at most $11$ vertices; $8$ of which are of Euler characteristic
$0$. These $8$ are $T_{7, 1, 2}, \dots, T_{11, 1, 2}$, $T_{3, 3,
0}$, $B_{3, 3}$ and $Q_{5, 2}$ of Examples $1, 3, 4, 6$.
\end{prop}

Here we prove.

\begin{theo}$\!\!${\bf .} \label{t1} Any degree-regular
triangulation of the torus is weakly regular.
\end{theo}

\begin{theo}$\!\!${\bf .} \label{t2} There exists an $n$-vertex
degree-regular triangulation of the Klein bottle if and only if
$n$ is a composite number $\geq 9$.
\end{theo}

\begin{theo}$\!\!${\bf .} \label{t3} $(a)$ For each $n \geq 12$
there exist $2$ distinct $n$-vertex weakly regular triangulations
of the torus. \vspace{0.5mm}
\newline $(b)$ For each $n \geq 18$ there exist $3$ distinct
$n$-vertex weakly regular triangulations of the torus.
\newline $(c)$ For each $m \geq 2$ there exists a $(4m + 2)$-vertex
weakly regular triangulation of the Klein bottle.
\end{theo}

\begin{theo}$\!\!${\bf .} \label{t4}
Let $T_{n, 1, k}$ be as in Example $1$. For a prime $n \geq 7$, if
$M$ is an $n$-vertex weakly regular triangulation of the torus
then $M$ is isomorphic to $T_{n, 1, k}$ for some $k$.
\end{theo}

\begin{cor} $\!\!${\bf .} \label{t5} $(a)$ For $n = 13$ or $17$,
there are exactly $2$ distinct $n$-vertex degree-regular
combinatorial $2$-manifolds of Euler characteristic $0$. These
are $T_{n, 1, 2}$ and $T_{n, 1, 3}$. \vspace{1mm}
\newline $(b)$ There are exactly $3$ distinct $19$-vertex
degree-regular combinatorial $2$-manifolds of Euler
characteristic $0$. These are $T_{19, 1, 2}$, $T_{19, 1, 3}$ and
$T_{19, 1, 7}$.
\end{cor}

From Theorem \ref{t1} and Proposition \ref{lutz1} we know all the
degree-regular triangulations of the torus on at most 15
vertices. Here we present (without using computer) the
following\,:

\begin{theo}$\!\!${\bf .} \label{t6} Let $M$ be an $n$-vertex
degree-regular combinatorial $2$-manifold of Euler characteristic
$0$. If $n = 12, 14$ or $15$ then $M$ is isomorphic to $T_{12, 1,
2}, \dots, T_{12, 1, 4}$, $T_{6, 2, 2}$, $T_{14, 1, 2}$, $T_{14,
1, 3}$, $T_{15, 1, 2}, \dots, T_{15, 1, 5}$, $Q_{7, 2}$, $Q_{5,
3}$, $B_{3, 4}$, $B_{4, 3}$, $B_{3, 5}$, $B_{5, 3}$ or $K_{3,
4}$. These $17$ combinatorial $2$-manifolds are pairwise
non-isomorphic. The first $10$ triangulate the torus and the
remaining $7$ triangulate the Klein bottle. Among the last $7$,
only $Q_{7, 2}$ is weakly regular.
\end{theo}

\section{Examples.}

In this section we present some degree-regular combinatorial
$2$-manifolds of Euler characteristic 0. First we give some
definitions and notations which will be used throughout the paper.

A 2-simplex in a 2-dimensional simplicial complex is also said to
be a {\em face}. We denote a face $\{u, v, w\}$ by $uvw$. We also
denote an edge $\{u, v\}$ by $uv$.

A {\em graph} is a simplicial complex of dimension at most one.
The complete graph on $n$ vertices is denoted by $K_n$. Disjoint
union of $m$ copies of $K_n$ is denoted by $mK_n$. A graph
without any edge is called a {\em null} graph. An $n$-vertex null
graph is denoted by $\emptyset_{n}$.

If $G$ is a graph and $n \geq 0$ is an integer then we define the
graph $G_n(G)$ as follows. The vertices of $G_n(G)$ are the
vertices of $G$. Two vertices $u$ and $v$ form an edge in
$G_n(G)$ if the number of common neighbours of $u$ and $v$ is
$n$. Clearly, if $G$ and $H$ are isomorphic then $G_n(G)$ and
$G_n(H)$ are isomorphic for all $n \geq 0$.

A connected finite graph is called a {\em cycle} if the degree of
each vertex is 2. An $n$-cycle is a cycle on $n$ vertices and is
denoted by $C_n$ (or by $C_n(a_1, \dots, a_n)$ if the edges are
$a_1a_2, \dots, a_{n-1}a_n, a_na_1$).  Disjoint union of $m$
copies of $C_n$ is denoted by $mC_n$.

For a simplicial complex $K$, the graph consisting of the edges
and vertices of $K$ is called the {\em edge-graph} of $K$ and is
denoted by ${\rm EG}(K)$. The complement of ${\rm EG}(K)$ is
called the {\em non-edge graph} of $K$ and is denoted by ${\rm
NEG}(K)$. Let $K$ be a simplicial complex with vertex-set $V(K)$.
If $U \subseteq V(K)$ then the {\em induced subcomplex} of $K$ on
$U$, denoted by $K[U]$, is the subcomplex whose simplices are
those of $K$ which are subsets of $U$.

If $v$ is a vertex of a simplicial complex $X$, then the {\em
link} of $v$ in $X$, denoted by ${\rm lk}_X(v)$ (or ${\rm
lk}(v)$), is the simplicial complex $\{\tau \in X ~ \colon ~ v
\not\in \tau, ~ \{v\}\cup\tau\in X\}$. If $v$ is a vertex of a
simplicial complex $X$, then the {\em star} of $v$ in $X$,
denoted by ${\rm st}_X(v)$ (or ${\rm st}(v)$), is the simplicial
complex $\{\{v\}, \tau, \tau \cup \{v\} ~ \colon ~ \tau \in {\rm
lk}_X(v)\}$. Clearly, a finite simplicial complex $K$ is a
combinatorial $2$-manifold if and only if ${\rm lk}_K(v)$ is a
cycle for each vertex $v$ of $K$.

\begin{eg}$\!\!${\rm {\bf :}} \label{e1} {\rm A series of weakly regular
orientable combinatorial 2-manifolds of Euler characteristic $0$.
For each $n\geq 7$ and each $k\in \{2, \dots, \lfloor\frac{n
-3}{2}\rfloor\}\cup \{\lceil\frac{n+1}{2}\rceil, \dots, n-3\}$,}
$$
T_{n,1,k} =  \{\{i, i+k, i+k+1\}, \{i, i+1, i+k+1\} ~ : ~ 1\leq
i\leq n\},
$$
{\rm where $V(T_{n,1,k}) = \{1, \dots, n\}$. Since ${\rm lk}(i) =
C_6(i+k, n+i-1, n+i-k-1, n+i-k, i+1, i+k+1)$, $T_{n, 1, k}$ is a
combinatorial $2$-manifold. Clearly, $T_{n, 1, k}$ triangulates
the torus and hence it is orientable. Since $\ZZ_n$ acts
transitively (by addition) on vertices, $T_{n, 1, k}$ is weakly
regular. (Here addition is modulo $n$.) In \cite{a1}, Altshuler
has shown that $T_{n, 1, k}$ is a subcomplex of an $n$-vertex
cyclic polytopal 3-sphere.}
\end{eg}

\begin{lemma}$\!\!${\bf .} \label{l2.1}
Let $T_{n, 1, k}$ be as above. We have the following\,:
\begin{enumerate}
     \item[$(a)$] $T_{n, 1, k}= T_{n, 1, n-k-1}$ for all $n$ and $k$.
     \item[$(b)$] $T_{n, 1, 2}\not\cong T_{n, 1, 3}$ for all
     $n\geq 12$.
     \item[$(c)$] $T_{n, 1, 2} \not\cong T_{n, 1, 4} \not\cong
     T_{n, 1, 3}$ for all $n \geq 20$.
     \item[$(d)$] $T_{12, 1, 2} \not\cong T_{12, 1, 4} \not\cong
     T_{12, 1, 3}$.
     \item[$(e)$] $T_{13, 1, 4} \cong T_{13, 1, 2} \cong T_{13, 1,
     5}$.
     \item[$(f)$] $T_{15, 1, k}\not\cong T_{15, 1, j}$ for $j, k
     \in \{2, 3, 4, 5\}$  and $j\neq k$.
     \item[$(g)$] $T_{16, 1, 5} \cong T_{16, 1, 2} \not\cong T_{16, 1, 6}
     \not\cong T_{16, 1, 3} \cong T_{16, 1, 4}$.
     \item[$(h)$] $T_{17, 1, 5} \cong T_{17, 1, 7} \cong T_{17, 1,
     2}$ and $T_{17, 1, 4} \cong T_{17, 1, 6} \cong T_{17, 1, 3}$.
     \item[$(i)$] $T_{18, 1, 6} \cong T_{18, 1, 5}$ and $T_{18, 1, k}
     \not\cong T_{18, 1, j}$ for $j, k \in \{2, 3, 4, 5, 7\}$, $j\neq
     k$.
     \item[$(j)$] $T_{19, 1, 6} \cong T_{19, 1, 8} \cong T_{19, 1, 2}
     \not\cong T_{19, 1, 7} \not\cong T_{19, 1, 3} \cong T_{19, 1, 4}
     \cong T_{19, 1, 5}$.
     \item[$(k)$] $T_{20, 1, 6} \cong T_{20, 1, 2}$, $T_{20, 1, 7} \cong
     T_{20, 1, 3}$ and $T_{20, 1, k}\not\cong T_{20, 1, j}$ for $j, k
     \in \{2, 3, 4, 5, 8\}$, $j\neq k$.
     \end{enumerate}
\end{lemma}

\noindent {\bf Proof.} Observe that ${\rm lk}_{T_{n, 1, k}}(i) =
C_6(i+1, i+k+1, i+k, n+i-1, n+i-k-1, n+i-k) = {\rm lk}_{T_{n, 1,
n-k-1}}(i)$. So, the faces in both $T_{n, 1, k}$ and $T_{n, 1,
n-k-1}$ are same. This proves $(a)$.

Then $G_4({\rm EG}(T_{n, 1, 2})) = C_n(1, \dots, n)$ for $n\geq
11$ and $G_4({\rm EG}(T_{n, 1, 3}))$ is a null graph for $n=13$
and for $n\geq 15$. Also, $G_4({\rm EG}(T_{14, 1, 2}))$ is $7K_2$.
So, $T_{n, 1, 2}\not\cong T_{n, 1, 3}$ for all $n\geq 13$.

Observe that $G_4({\rm EG}(T_{12, 1, 3}))$ is a 12-cycle with
edges $\{i, i+5\}$, $1\leq i\leq 12$. So, the edges of $G_4({\rm
EG}(T_{12, 1, 3}))$ are non-edges of $T_{12, 1, 2}$. But
$G_4({\rm EG}(T_{12, 1, 2}))$ is a subgraph of ${\rm EG}(T_{12, 1,
2})$. So, $T_{12, 1, 2} \not\cong T_{12, 1, 3}$. This proves
$(b)$.

For all $n\geq 20$, $G_3({\rm EG}(T_{n, 1, 4}))$ is a null graph,
but $\{i, i+2\}$ is an edge in $G_3({\rm EG}(T_{n, 1, 3}))$. So,
$G_3({\rm EG}(T_{n, 1, 4})) \not\cong G_3({\rm EG}(T_{n, 1, 3}))$
and hence $T_{n, 1, 4} \not\cong T_{n, 1, 3}$ for $n \geq 20$.

Again, for all $n\geq 20$, $G_4({\rm EG}(T_{n, 1, 4})) =
\emptyset_n$. So, $G_4({\rm EG}(T_{n, 1, 4})) \not\cong G_4({\rm
EG}(T_{n, 1, 2}))$ and hence $T_{n, 1, 4} \not\cong T_{n, 1, 2}$
for $n\geq 20$. This proves $(c)$.

Observe that $G_4({\rm EG}(T_{12, 1, 4})) = 3K_4$ (with edges
$\{i, j\}$, $i-j \equiv 0$ (mod 3)). Since, $G_4({\rm EG}(T_{12,
1, 2}))$ and $G_4({\rm EG}(T_{12, 1, 3}))$ are 12-cycles, $T_{12,
1, 2} \not \cong T_{12, 1, 4} \not\cong T_{12, 1, 3}$. This
proves $(d)$.

The map $i \mapsto 4i$ (mod 13) defies an isomorphism from
$T_{13, 1, 2}$ to $T_{13, 1, 4}$ and the map $i \mapsto 7i$ (mod
13) defies an isomorphism from $T_{13, 1, 2}$ to $T_{13, 1, 5}$.
This proves $(e)$.

Observe that $G_0({\rm EG}(T_{15, 1, 2}))$ is a 15-cycle (with
edges $\{i, i+7\}$, $1\leq i\leq 15$) and $G_0({\rm EG}(T_{15, 1,
i}))$ is a null graph for $i=3, 4, 5$. So, $T_{15, 1, i}\not\cong
T_{15, 1, 2}$ for $i=3, 4$ or 5. Again, $G_4({\rm EG}(T_{15, 1,
3}))$ is a null graph, $G_4({\rm EG}(T_{15, 1, 4})) = 3C_5$ (with
edges $\{i, i+6\}$, $1\leq i\leq 15$) and $G_4({\rm EG}(T_{15, 1,
5})) = C_{15}$ (with edges $\{i, i+4\}$, $1\leq i\leq 15$). So,
$T_{15, 1, 3}\not\cong T_{15, 1, 4}\not\cong T_{15, 1,
5}\not\cong T_{15, 1, 3}$. These prove $(f)$.

The map $i\mapsto 3i$ (mod 16) defines an isomorphism from $T_{16,
1, 4}$ to $T_{16, 1, 3}$ and an isomorphism from $T_{16, 1, 5}$ to
$T_{16, 1, 2}$.

Observe that $G_4({\rm EG}(T_{16, 1, 6})) = 8K_2$ (with edges
$\{i, i+8\}$, $1\leq i\leq 8$). Since, $G_4({\rm EG}(T_{16, 1,
2})) = C_{16}$ and $G_4({\rm EG}(T_{16, 1, 3}))$ is a null graph,
$T_{16, 1, 2} \not \cong T_{16, 1, 6} \not\cong T_{16, 1, 3}$.
This proves $(g)$.

The map $i \mapsto 14i$ (mod 17) defines an isomorphism from
$T_{17, 1, 5}$ to $T_{17, 1, 2}$. The map $i \mapsto 2i$ (mod 17)
defines an isomorphism from $T_{17, 1, 7}$ to $T_{17, 1, 2}$. The
map $i \mapsto 13i$ (mod 17) defines an isomorphism from $T_{17,
1, 4}$ to $T_{17, 1, 3}$. The map $i \mapsto 3i$ (mod 17) defines
an isomorphism from $T_{17, 1, 6}$ to $T_{17, 1, 3}$. This proves
$(h)$.

The map $i \mapsto 5i$ (mod 18) defines an isomorphism from
$T_{18, 1, 6}$ to $T_{18, 1, 5}$.

Now, $G_3({\rm EG}(T_{18, 1, 3})) = 2C_9$ (with edges $\{i,
i+2\}$, $1\leq i\leq 18$) and $G_3({\rm EG}(T_{18, 1, 7})) = 9K_2$
(with edges $\{i, i+9\}$, $1\leq i\leq 9$). So, $T_{18, 1, 3} \not
\cong T_{18, 1, 7}$. Again, $G_4({\rm EG}(T_{18, 1, 2})) =
C_{18}(1, 2, \dots, 18) \subseteq {\rm EG}(T_{18, 1, 2})$,
$G_4({\rm EG}(T_{18, 1, 3}))$ and $G_4({\rm EG}(T_{18, 1, 7}))$
are null graphs, \linebreak $G_4({\rm EG}(T_{18, 1, 4})) = 9K_2$
(with edges $\{i, i+9\}$, $1\leq i\leq 9$) and $G_4({\rm
EG}(T_{18, 1, 5})) = C_{18}$ (with edges $\{i, i+7\}$, $1\leq
i\leq 18$). Thus, $G_4({\rm EG}(T_{18, 1, 5}))$ is not a subgraph
of ${\rm EG}(T_{18, 1, 5})$. These imply $(i)$.

The map $i\mapsto 3i$ (mod 19) defines an isomorphism from $T_{19,
1, 5}$ to $T_{19, 1, 3}$. The map $i\mapsto 15i$ (mod 19) defines
an isomorphism from $T_{19, 1, 4}$ to $T_{19, 1, 3}$. The map
$i\mapsto 6i$ (mod 19) defines an isomorphism from $T_{19, 1, 2}$
to $T_{19, 1, 6}$. The map $i\mapsto 9i$ (mod 19) defines an
isomorphism from $T_{19, 1, 2}$ to $T_{19, 1, 8}$.

The graph $G_4({\rm EG}(T_{19, 1, 7}))$ is null, where as
$G_4({\rm EG}(T_{19, 1, 2})) = C_{19}$. So, $T_{19, 1, 7}\not\cong
T_{19, 1, 2}$. Again, $G_0({\rm EG}(T_{19, 1, 7}))$ is null, where
as $G_0({\rm EG}(T_{19, 1, 3})) = C_{19}$ (with edges $\{i,
i+9\}$, $1\leq i\leq 19$). So, $T_{19, 1, 7}\not\cong T_{19, 1,
3}$. This proves $(j)$.

The map $i\mapsto 3i$ (mod 20) defines an isomorphism from $T_{20,
1, 6}$ to $T_{20, 1, 2}$ and an isomorphism from $T_{20, 1, 7}$ to
$T_{20, 1, 3}$.

Observe that $G_4({\rm EG}(T_{20, 1, 2})) = C_{20}$, $G_4({\rm
EG}(T_{20, 1, 3}))$, $G_4({\rm EG}(T_{20, 1, 4}))$, $G_4({\rm
EG}(T_{20, 1, 5}))$ are null graphs and $G_4({\rm EG}(T_{20, 1,
8})) = 10K_2$ (with edges $\{i, i+10\}$, $1\leq i\leq 10$). So,
$T_{20, 1, 2}\not \cong T_{20, 1, i}$ for $i=3, 4, 5, 8$ and
$T_{20, 1, 8}\not \cong T_{20, 1, i}$ for $i=3, 4, 5$.

Again, $G_3({\rm EG}(T_{20, 1, 3})) = 2C_{10}$ (with edges $\{i,
i+2\}$, $1\leq i\leq 20$) but $G_3({\rm EG}(T_{20, 1, 4}))$ and
$G_3({\rm EG}(T_{20, 1, 5}))$ are null graphs. So, $T_{20, 1,
3}\not \cong T_{20, 1, i}$ for $i= 4, 5$.

Finally, if possible let $\varphi$ be an isomorphism from $T_{20,
1, 4}$ to  $T_{20, 1, 5}$. Then $\varphi$ induces isomorphism
between $G_n({\rm EG}(T_{20, 1, 4}))$ and $G_n({\rm EG}(T_{20, 1,
5}))$ for each $n$. Since, ${\rm Aut}(T_{20, 1, 4})$ acts
transitively on $V(T_{20, 1, 4})$, we can assume that
$\varphi(20) = 20$. Since $G_0({\rm EG}(T_{20, 1, 4})) =
C_{20}(20, 7, 14, \dots, 13)$ and $G_0({\rm EG}(T_{20, 1, 5})) =
C_{20}(20, 3, 6, \dots, 17)$, $\varphi(7) = 17$ or 3. If
$\varphi(7) = 17$ then $\varphi(14) = 14$, $\varphi(1) = 11,
\dots$. In that case, $\varphi(\{20, 1\}) = \{20, 11\}$. This is
a contradiction since $\{20, 1\}$ is an edge in $T_{20, 1, 4}$
but $\{20, 11\}$ is not an edge in $T_{20, 1, 5}$. Similarly, we
get a contradiction if $\varphi(7)= 3$. This proves $(k)$. \hfill
$\Box$

\medskip

Let $D_n$ denote the dihedral group of order $2n$ and $\ZZ_{m^2 -
m + 1}\colon \ZZ_{6} := < \rho, \mu ~ \colon ~ \rho^{m^2 - m + 1}
= 1 = \mu^6, \mu^{- 1} \rho \mu = \rho^m>$ for $m \geq 3$. In
\cite{l}, Lutz has shown that ${\rm Aut}(T_{n, 1, k}) = D_n$ for
$(n, k) = (9, 2), \dots, (15, 2), (12, 3), (14, 3), (15, 5)$,
${\rm Aut}(T_{12, 1, 4}) = D_4\times D_3$, ${\rm Aut}(T_{15, 1,
k}) = D_5 \times D_3$ for $k = 3, 4$ and ${\rm Aut}(T_{m^2-m+1, 1,
m-1}) = \ZZ_{m^2-m+1}\colon \ZZ_6$ for $m = 3, 4$. Here we prove
the following.

\begin{lemma}$\!\!${\bf .} \label{l2.2} $(a)$ $D_n$ acts
face-transitively on $T_{n, 1, k}$ for all $n\geq 7$ and for all
$k$. \vspace{1mm}
\newline $(b)$ ${\rm Aut}(T_{n, 1,2}) = D_n$ for all $n \geq 9$.
\vspace{1mm}
\newline $(c)$ $D_{2m}\times D_{m+1} \leq {\rm Aut}(T_{2m^2+2m, 1,
2m})$ for $m\geq 2$. \vspace{1mm}
\newline $(d)$ $D_{m+1}\times D_{m-1} \leq {\rm Aut}(T_{m^2-1, 1,
m-1})$ for $m\geq 4$. \vspace{1mm}
\newline $(e)$ $D_{m+1}\times D_{m-1} \leq {\rm Aut}(T_{m^2-1, 1,
m})$ for $m\geq 4$. \vspace{1mm}
\newline $(f)$ $\ZZ_{m^2-m+1}\colon \ZZ_6 \leq {\rm Aut}(T_{m^2-m+1,
1, m-1})$ for $m\geq 3$. \vspace{1mm}
\newline Here $H\leq G$ means $G$ has a subgroup isomorphic to
$H$.
\end{lemma}

\noindent {\bf Proof.} Let $\alpha_n, \beta_n \colon V(T_{n, 1,
k}) \to V(T_{n, 1, k})$ be given by $\alpha_n(i) = i+1$ and
$\beta_n(i) = n-i$ (modulo $n$). Let $A_{n, k, i} := \{i, i+1,
i+k+1\}$ and $B_{n, k, i} := \{i, i+k, i+k+1\}$. Then
$\alpha_n(A_{n, k, i}) = A_{n, k, i+1}$, $\alpha_n(B_{n, k, i}) =
B_{n, k, i+1}$, $\beta_n(A_{n, k, i}) = B_{n, k, n-i-k-1}$ and
$\beta_n(B_{n, k, i}) = A_{n, k, n-i-k-1}$. So, $\alpha_n,
\beta_n\in {\rm Aut}(T_{n, 1, k})$. Clearly, the order of
$\alpha_n$ is $n$, the order of $\beta_n$ is 2 and $\beta_n
\alpha_n \beta_n = \alpha_n^{-1}$. Thus, $<\alpha_n, \beta_n>$ is
isomorphic to $D_n$. Clearly, the action of $<\alpha_n, \beta_n>$
on $T_{n, 1, k}$ is transitive on the faces. This proves $(a)$.

For $n\geq 11$, $G_4({\rm EG}(T_{n, 1, 2})) = C_n(1, 2, \dots,
n)$. Therefore, $<\alpha_n, \beta_n> \leq {\rm Aut}(T_{n, 1, 2})
\leq {\rm Aut}(G_4({\rm EG}(T_{n, 1, 2}))) = {\rm Aut}(C_n(1,
\dots, n)) = <\alpha_n, \beta_n>$. Thus, ${\rm Aut}(T_{n, 1, 2})
= <\alpha_n, \beta_n>$ for $n\geq 11$. Since $G_2({\rm EG}(T_{10,
1, 2})) = C_{10}(1, 4, 7, 10, 3, 6, 9, 2, 5, 8)$,
$<\alpha_{10}^3, \beta_{10}> \leq {\rm Aut}(T_{10, 1, 2}) \leq
{\rm Aut}(G_2({\rm EG}(T_{10, 1, 2}))) = {\rm Aut}(C_{10}(1, 4,
\dots, 5, 8)) = <\alpha_{10}^3, \beta_{10}>$. Thus, ${\rm
Aut}(T_{10, 1, 2}) = <\alpha_{10}^3, \beta_{10}> = <\alpha_{10},
\beta_{10}>$. Observe that ${\rm NEG}(T_{9, 1, 2}) = C_9(1, 5, 9,
4, 8, 3, 7, 2, 6)$. Therefore, $< \alpha_9^4, \beta_9>$ $\leq
{\rm Aut}(T_{9, 1, 2}) \leq {\rm Aut}({\rm NEG}(T_{9, 1, 2})) =
{\rm Aut}(C_9(1, 5, 9, 4, 8, 3, 7, 2, 6)) = <\alpha_9^4,
\beta_9>$. Thus, ${\rm Aut}(T_{9, 1, 2})  = <\alpha_9^4, \beta_9>
= <\alpha_9, \beta_9>$. This proves $(b)$.

Let $\alpha, \beta, \gamma, \delta \colon V(T_{2m^2+2m, 1, 2m})
\to V(T_{2m^2+2m, 1, 2m})$ be given by $\alpha(i) = i+m+1$,
$\beta(i) = i+2m$, $\gamma(i) = (2m + 1)i$ and $\delta(i) = (2m^2
- 1)i$ (i.e., $\alpha = \alpha_{2m^2 + 2m}^{m+1}$, $\beta =
\alpha_{2m^2 + 2m}^{2m}$). Then $\alpha, \beta, \gamma, \delta \in
{\rm Aut}(T_{2m^2+2m, 1, 2m})$, order of $\alpha$ is $2m$, order
of $\beta$ is $m+1$, order of $\gamma$ is 2, order of $\delta$ is
2, $\alpha\beta = \beta\alpha$, $\alpha\gamma = \gamma\alpha$,
$\beta\delta = \delta\beta$, $\gamma\delta = \delta\gamma=
\beta_{2m^2+2m}$, $\delta\alpha\delta = \alpha^{-1}$,
$\gamma\beta\gamma = \beta^{-1}$. Therefore, $<\alpha, \beta,
\gamma, \delta> = <\alpha, \delta> \times <\beta, \gamma> \cong
D_{2m}\times D_{m+1}$. This proves $(c)$.

\smallskip

\noindent {\bf Claim.} If $m^2 \equiv 1$ (mod $n$) then $\mu(i) =
mi$ and $\nu(i) = (n - m)i$ define two distinct involutions
(automorphisms of order 2) of $T_{n, 1, k}$ for each $k \in \{m-1,
m\}$.

Let $A_{n, k, i}$ and $B_{n, k, i}$ be as above. Then $\mu(A_{n,
m-1, i}) = A_{n, m-1, mi}$, $\mu(B_{n, m-1, i}) = B_{n, m-1,
mi-m+1}$, $\mu(A_{n, m, i}) = B_{n, m, mi}$ and $\mu(B_{n, m, i})
= A_{n, m, mi}$. Thus $\mu\in {\rm Aut}(T_{n, 1, k})$ and hence
$\nu = \beta_n\mu = \mu\beta_n\in {\rm Aut}(T_{n, 1, k})$ for $k =
m-1, m$. Since $\mu^{-1}=\mu$, $\nu^{-1}=\nu$ and $\mu\nu =
\beta_n\neq $ the identity, $\mu\neq \nu$. This proves the claim.

For $k = m-1, m$, let $\alpha, \lambda, \mu, \nu \colon V(T_{m^2 -
1, 1, k}) \to V(T_{m^2 - 1, 1, k})$ be given by $\alpha(i) = i + m
+ 1$, $\lambda(i) = i + m - 1$, $\mu(i) = mi$ and $\nu(i) = (n -
m)i$ (i.e., $\alpha = \alpha_{m^2 -1}^{m+1}$, $\lambda =
\alpha_{m^2 -1}^{m-1}$). Then $\alpha$ and $\lambda$ are
automorphisms of $T_{m^2-1, 1, k}$. Also, by the above claim,
$\mu$ and $\nu$ are distinct automorphisms of $T_{m^2-1, 1, k}$.
Clearly, the orders of $\alpha$, $\lambda$, $\mu$ and $\nu$ are
$m-1$, $m+1$, 2 and 2 respectively. Again, $\alpha\mu =
\mu\alpha$, $\alpha\lambda = \lambda\alpha$, $\lambda\nu =
\nu\lambda$, $\mu\nu = \nu\mu$, $\nu\alpha\nu = \alpha^{-1}$ and
$\mu\lambda\mu = \lambda^{-1}$. Thus, $<\alpha, \lambda, \mu, \nu>
= <\lambda, \mu>\times <\alpha, \nu> \cong D_{m+1}\times D_{m-1}$.
This proves $(d)$ and $(e)$.

Let $\sigma \colon V(T_{m^2 - m + 1, 1, m - 1}) \to V(T_{m^2 - m +
1, 1, m - 1})$ be given by $\sigma(i) = mi$. Then \linebreak
$\sigma(A_{m^2 - m + 1, m - 1, i}) = B_{m^2-m+1, m-1, mi}$ and
$\sigma(B_{m^2-m+1, m-1, i}) = A_{m^2-m+1, m-1, mi-1}$. Thus,
$\sigma\in {\rm Aut}(T_{m^2-m+1, 1, m-1})$. Since 6 is the
smallest positive integer $n$ for which $m^n -1$ is divisible by
$m^2-m+1$, the order of $\sigma$ is 6. Now, if $\rho = \alpha_{m^2
- m +1}$ (i.e., $\rho(i) = i+1$) then $\rho$ is an automorphism of
order $m^2-m+1$ and $\sigma^{-1}\rho\sigma(i) = m((1-m)i + 1) = i
+ m = \rho^m(i)$. Thus, $<\rho, \sigma> = \ZZ_{m^2-m+1} \colon
\ZZ_6$. This proves $(f)$. \hfill $\Box$

\begin{eg}$\!\!${\rm {\bf :}} \label{e2} {\rm A series of weakly
regular orientable combinatorial 2-manifolds of Euler
characteristic $0$. For each $n \geq 4$ and each $k = 1, \dots, n
- 3$,}
$$
T_{n,2,k} =  \{u_{i}u_{i+1}v_{i+1}, u_{i}v_{i}v_{i+1},
u_{i+k}u_{i+k+1}v_{i}, u_{i+k+1}v_{i}v_{i+1} ~ : ~ 1\leq i\leq n
\},
$$
{\rm where $V(T_{n,2,k}) = \{u_1, \dots, u_n\}\cup\{v_1, \dots,
v_n\}$. (Addition in the subscripts are modulo $n$.) Since ${\rm
lk}(u_i) = C_6(u_{i-1}, v_i, v_{i+1}, u_{i+1}, v_{n+i-k-1},
v_{n+i-k})$ and ${\rm lk}(v_j) = C_6(v_{j-1}, u_{j-1}, u_{j},
v_{j+1}$, $u_{j+k+1}, u_{j+k})$, $T_{n, 2, k}$ is a degree-regular
combinatorial $2$-manifold on $2n$ vertices. Clearly, $T_{n, 2,
k}$ triangulates the torus and hence it is orientable. If
$\alpha$, $\beta \colon V(T_{n, 2, k}) \to V(T_{n, 2, k})$ are the
maps given by $\alpha(u_i) = u_{i+1}$, $\alpha(v_i) = v_{i+1}$,
$\beta(u_i) = v_i$, $\beta(v_i) = u_{i+k}$ for $1\leq i\leq n$,
then $\alpha$ and $\beta$ are automorphisms of $T_{n, 2, k}$ and
hence $<\alpha, \beta>$ is a subgroup of ${\rm Aut}(T_{n, 2, k})$.
Clearly, $<\alpha, \beta>$ acts transitively on vertices. Thus
$T_{n, 2, k}$ is weakly regular. }
\end{eg}

\begin{lemma}$\!\!${\bf .} \label{l2.3} Let $T_{n, 1, j}$ and
$T_{m, 2, k}$ be as in Examples $\ref{e1}$, $\ref{e2}$. We have
the following\,:
\begin{enumerate}
     \item[$(a)$] If $n$ and $k$ are relatively prime or $n$
     and $k+2$ are relatively prime then $T_{n, 2, k}$ is
     isomorphic to $T_{2n, 1, j}$ for some $j$.
     \item[$(b)$] $T_{6, 2, 2} \not\cong T_{12, 1, i}$ for all
     $i$.
     \item[$(c)$] $T_{8, 2, 4} \cong T_{8, 2, 2} \not\cong
     T_{16, 1, i}$ for all $i$.

\end{enumerate}
\end{lemma}


\noindent {\bf Proof.} If $(n, k) =1$ and $k \leq n-3$, then there
exists $p\in \{1, \dots, n - 2\}$ such that $pk \equiv 1$ (mod
$n$). Since $k\leq n-3$ and $n\geq 4$, $2p \neq n$, $n-1$. ($2p =
n$ implies $2\equiv 2pk \equiv 0$ (mod $n$), a contradiction. $2p
= n-1$ implies $2 \equiv 2pk \equiv n-k$. This implies $k \equiv n
- 2$ (mod $n$) and hence $k=n-2$.) Let $\varphi \colon V(T_{n, 2,
k}) \to \{1, \dots, 2n\}$ be given by $\varphi(u_{i}) = 1 + 2p(i -
1)$ and $\varphi(v_{i}) = 2+ 2p(i - 1)$ (modulo $2n$). Since $(n,
p) = 1$, $\varphi$ is a bijection. Now, $\varphi(u_i u_{i + 1}
v_{i + 1}) = \{1+2(i-1)p, 1+2ip, 2+ 2ip\}$, $\varphi(u_i v_{i}
v_{i+1}) = \{1+2(i-1)p, 2+2(i-1)p, 2+ 2ip\}$, $\varphi(u_{i+k}
u_{i+k+1} v_{i}) = \{1+2(i-1)p+2, 1+2ip+2, 2+ 2(i-1)p\}
=\{2(i-1)p+2, 2(i-1)p+3, 2ip+3\}$, $\varphi(u_{i+k+1} v_{i}
v_{i+1}) = \{2ip+3, 2(i-1)p +2, 2+ 2ip\} =\{2(i-1)p+2, 2ip+2,
2ip+3\}\in T_{2n, 1, 2p}$. This shows that $\varphi \colon T_{n,
2, k} \to T_{2n, 1, 2p}$ is an isomorphism.

If $(n, k + 2) =1$ then assume that $(n, k) \neq 1$ (otherwise
there is nothing to prove). Let $p(k + 2) \equiv 1$ (mod $n$).
Observe that $2p \neq n, n - 1$. ($2p = n$ implies $2\equiv 2p(k +
2) \equiv 0$ (mod $n$), a contradiction. $2p = n - 1$ implies $2
\equiv 2p(k + 2) \equiv (n - 1)(k + 2) \equiv n - k - 2$. This
implies $k \equiv n - 4$ (mod $n$) and hence $k = n - 4$. Since
$(n, k) \neq 1$, $n$ and $k$ are even and hence $(n, k + 2) \neq
1$, a contradiction.) Let $\psi \colon V(T_{n, 2, k}) \to \{1,
\dots, 2n\}$ be given by $\psi(u_{i}) = 1 + 2p(i - 1)$ and
$\psi(v_{i}) = 2 + 2p(i - 2)$ (modulo $2n$). Since $(n, p) = 1$,
$\psi$ is a bijection. Similar argument as before shows that $\psi
\colon T_{n, 2, k} \to T_{2n, 1, 2p - 1}$ is an isomorphism. This
proves $(a)$.

Since $G_4({\rm EG}(T_{6, 2, 2})) = 3K_4$ and $G_4({\rm
EG}(T_{12, 1, i})) = C_{12}$ for $i= 2, 3$, $T_{6, 2, 2}\not\cong
T_{12, 1, i}$ for $i=2, 3$.

Now, $G_3({\rm EG}(T_{12, 1, 4})) = 4C_3$ (with edges $\{i, i +
4\}$, $1 \leq i \leq 12$). So, $G_3({\rm EG}(T_{12, 1, 4}))$ is a
subgraph of ${\rm EG}(T_{12, 1, 4})$. Whereas, $G_3({\rm
EG}(T_{6, 2, 2})) = 4C_3$ (with edges $\{u_i, u_j\}$, $\{v_i,
v_j\}$, where $i-j\equiv 0$ (mod 2)). So, edges of $G_3({\rm
EG}(T_{6, 2, 2}))$ are not in ${\rm EG}(T_{6, 2, 2})$. Thus,
$T_{6, 2, 2}\not\cong T_{12, 1, 4}$. This proves $(b)$.

The map $u_i\mapsto u_i$, $v_i\mapsto v_{i+3}$ defines an
isomorphism from $T_{8, 2, 2}$ to $T_{8, 2, 4}$. (As usual,
addition in the subscripts are modulo 8.) Observe that $G_1({\rm
EG}(T_{8, 2, 2}))$ is a null graph, whereas $G_1({\rm EG}(T_{16,
1, 2})) = 2C_8$, $G_1({\rm EG}(T_{16, 1, 3})) = 2C_8$, $G_1({\rm
EG}(T_{16, 1, 6})) = 4C_4$. Thus $T_{8, 2, 2} \not\cong T_{16, 1,
j}$ for $j = 2, 3$ and  6. Therefore, by Lemma \ref{l2.1} $(g)$,
$T_{8, 2, 2} \not\cong T_{16, 1, i}$ for all $i$. This proves
$(c)$. \hfill $\Box$

\bigskip

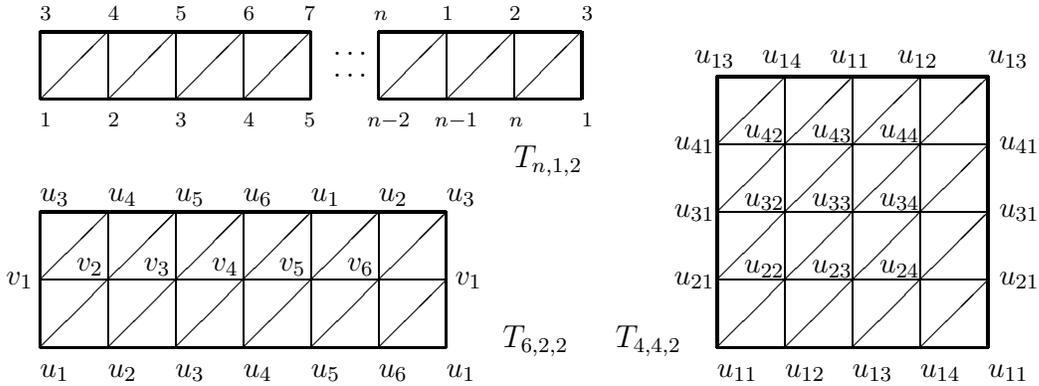
\hrule

\setlength{\unitlength}{3mm}

\begin{picture}(45,19)(-2,-2.5)


\thicklines

\put(0,11){\line(1,0){12}} \put(0,11){\line(0,1){3}}
\put(0,14){\line(1,0){12}} \put(12,11){\line(0,1){3}}

\put(13,12){$\dots$} \put(13,13){$\dots$}

 \put(15,11){\line(1,0){9}} \put(15,11){\line(0,1){3}}
\put(15,14){\line(1,0){9}} \put(24,11){\line(0,1){3}}

\thinlines

\put(3,11){\line(0,1){3}} \put(6,11){\line(0,1){3}}
\put(9,11){\line(0,1){3}} \put(18,11){\line(0,1){3}}
\put(21,11){\line(0,1){3}}

\put(0,11){\line(1,1){3}} \put(3,11){\line(1,1){3}}
\put(6,11){\line(1,1){3}} \put(9,11){\line(1,1){3}}
\put(15,11){\line(1,1){3}} \put(18,11){\line(1,1){3}}
\put(21,11){\line(1,1){3}}

\put(0,10){\mbox{${}_1$}} \put(3,10){\mbox{${}_2$}}
\put(6,10){\mbox{${}_3$}} \put(9,10){\mbox{${}_4$}}
\put(11.7,10){\mbox{${}_5$}} \put(14.5,10){\mbox{${}_{n-2}$}}
\put(17.5,10){\mbox{${}_{n-1}$}} \put(20.8,10){\mbox{${}_{n}$}}
\put(24,10){\mbox{${}_1$}}

\put(0,14.7){\mbox{${}_3$}} \put(3,14.7){\mbox{${}_4$}}
\put(6,14.7){\mbox{${}_5$}} \put(9,14.7){\mbox{${}_6$}}
\put(11.7,14.7){\mbox{${}_7$}} \put(14.8,14.7){\mbox{${}_{n}$}}
\put(17.8,14.7){\mbox{${}_{1}$}} \put(20.8,14.7){\mbox{${}_{2}$}}
\put(24,14.7){\mbox{${}_3$}}

{\large \put(21,8){\mbox{$T_{n,1,2}$}} }


\thicklines

\put(0,0){\line(1,0){18}} \put(0,0){\line(0,1){6}}
\put(0,6){\line(1,0){18}} \put(18,0){\line(0,1){6}}

\thinlines

\put(0,3){\line(1,0){18}} \put(3,0){\line(0,1){6}}
\put(6,0){\line(0,1){6}} \put(9,0){\line(0,1){6}}
\put(12,0){\line(0,1){6}} \put(15,0){\line(0,1){6}}

\put(0,0){\line(1,1){6}} \put(0,3){\line(1,1){3}}
\put(3,0){\line(1,1){6}} \put(6,0){\line(1,1){6}}
\put(9,0){\line(1,1){6}} \put(12,0){\line(1,1){6}}
\put(15,0){\line(1,1){3}}

\put(0,-1.3){\mbox{$u_{1}$}} \put(3,-1.3){\mbox{$u_{2}$}}
\put(6,-1.3){\mbox{$u_{3}$}} \put(9,-1.3){\mbox{$u_{4}$}}
\put(12,-1.3){\mbox{$u_{5}$}} \put(15,-1.3){\mbox{$u_{6}$}}
\put(18,-1.3){\mbox{$u_{1}$}}

\put(-1.5,2.8){\mbox{$v_{1}$}} \put(1.6,3.4){\mbox{$v_{2}$}}
\put(4.6,3.4){\mbox{$v_{3}$}} \put(7.6,3.4){\mbox{$v_{4}$}}
\put(10.6,3.4){\mbox{$v_{5}$}} \put(13.6,3.4){\mbox{$v_{6}$}}
\put(18.4,2.8){\mbox{$v_{1}$}}

\put(0,6.5){\mbox{$u_{3}$}} \put(3,6.5){\mbox{$u_{4}$}}
\put(6,6.5){\mbox{$u_{5}$}} \put(9,6.5){\mbox{$u_{6}$}}
\put(12,6.5){\mbox{$u_{1}$}} \put(15,6.5){\mbox{$u_{2}$}}
\put(18,6.5){\mbox{$u_{3}$}}

{\large \put(20.5,0){\mbox{$T_{6,2,2}$}} }


\thicklines

\put(30,0){\line(1,0){12}} \put(30,0){\line(0,1){12}}
\put(30,12){\line(1,0){12}} \put(42,0){\line(0,1){12}}

\thinlines

\put(30,3){\line(1,0){12}} \put(30,6){\line(1,0){12}}
 \put(30,9){\line(1,0){12}} \put(33,0){\line(0,1){12}}
 \put(36,0){\line(0,1){12}} \put(39,0){\line(0,1){12}}

\put(30,9){\line(1,1){3}} \put(30,0){\line(1,1){12}}
\put(30,3){\line(1,1){9}} \put(30,6){\line(1,1){6}}
\put(33,0){\line(1,1){9}} \put(36,0){\line(1,1){6}}
\put(39,0){\line(1,1){3}}


\put(30,-1.3){\mbox{$u_{11}$}} \put(33,-1.3){\mbox{$u_{12}$}}
\put(36,-1.3){\mbox{$u_{13}$}} \put(39,-1.3){\mbox{$u_{14}$}}
\put(42,-1.3){\mbox{$u_{11}$}}

\put(28.1,2.8){\mbox{$u_{21}$}} \put(31.2,3.3){\mbox{$u_{22}$}}
\put(34.2,3.3){\mbox{$u_{23}$}} \put(37.2,3.3){\mbox{$u_{24}$}}
\put(42.5,2.8){\mbox{$u_{21}$}}

\put(28.1,5.8){\mbox{$u_{31}$}} \put(31.2,6.3){\mbox{$u_{32}$}}
\put(34.2,6.3){\mbox{$u_{33}$}} \put(37.2,6.3){\mbox{$u_{34}$}}
\put(42.5,5.8){\mbox{$u_{31}$}}

\put(28.1,8.8){\mbox{$u_{41}$}} \put(31.2,9.3){\mbox{$u_{42}$}}
\put(34.2,9.3){\mbox{$u_{43}$}} \put(37.2,9.3){\mbox{$u_{44}$}}
\put(42.5,8.8){\mbox{$u_{41}$}}

\put(29,12.6){\mbox{$u_{13}$}} \put(32,12.6){\mbox{$u_{14}$}}
\put(35,12.6){\mbox{$u_{11}$}} \put(38,12.6){\mbox{$u_{12}$}}
\put(42,12.6){\mbox{$u_{13}$}}

{\large \put(25.5,0){\mbox{$T_{4,4,2}$}} }

\end{picture}

\hrule

\begin{eg}$\!\!${\rm {\bf :}} \label{e3} {\rm Some more weakly
regular orientable combinatorial 2-manifolds of Euler
characteristic $0$. For $n, m \geq 3$ and $k = 0, \dots, n - 1$,}
\begin{eqnarray*} T_{n, m, k} & = & \{u_{i, j} u_{i, j + 1}u_{i +
1, j + 1}, u_{i, j} u_{i + 1, j} u_{i + 1, j + 1} ~ : ~ 1 \leq i
\leq m - 1, ~ 1 \leq j \leq n \} \\
&& \cup \{u_{m, j} u_{m, j + 1} u_{1, j + k + 1}, u_{m, j}u_{1, j
+ k} u_{1, j + k + 1} ~ : ~ 1 \leq j \leq n \},
\end{eqnarray*}
{\rm where $V(T_{n, m, k}) = \{u_{i, j} ~ : 1 \leq i \leq m, 1
\leq j \leq n\}$. (Addition in the second subscripts are modulo
$n$.) Clearly, $T_{n, m, k}$ triangulates the torus and hence it
is an orientable combinatorial $2$-manifold on $mn$ vertices.
Since the degree of each vertex is $6$, $T_{n, m, k}$ is
degree-regular. ($T_{n, m, 0}$ was earlier defined in \cite{dn} as
$A_{m, n}$.) If $\sigma$, $\gamma \colon V(T_{n, m, k}) \to
V(T_{n, m, k})$ are the maps given by $\gamma(u_{i,j}) = u_{i + 1,
j}$ for $1 \leq i \leq m - 1$, $\gamma(u_{m, j}) = u_{1, j + k}$
and $\sigma(u_{i, j}) = u_{i, j + 1}$ then $\sigma$ and $\gamma$
are automorphisms of $T_{n, m, k}$ and hence $<\sigma, \gamma>$ is
a subgroup of ${\rm Aut}(T_{n, m, k})$. For any vertex $u_{i, j}$,
$\sigma^{j - 1} \gamma^{i - 1}(u_{1, 1}) = u_{i, j}$. This implies
that the action of $<\sigma, \gamma>$ is vertex-transitive. Thus
$T_{n, m, k}$ is weakly regular.

}
\end{eg}

\begin{lemma}$\!\!${\bf .} \label{l2.4} Let $T_{n, m, k}$ be as
in Examples $\ref{e1}$, $\ref{e2}$, $\ref{e3}$. Then
\begin{enumerate}
     \item[$(a)$] If $n$ and $k$ are relatively prime or $n$
     and $k+m$ are relatively prime then $T_{n, m, k}$ is
     isomorphic to $T_{nm, 1, j}$ for some $j$.
    \item[$(b)$] $T_{4, 4, 2}\cong T_{8, 2, 2}$.
    \item[$(c)$] $T_{16, 1, k} \not\cong T_{4, 4, 0} \not\cong
    T_{8, 2, j}$ for all $k$ and $j$.
\end{enumerate}
\end{lemma}

\noindent {\bf Proof.} Since $(n, k) =1$, there exists $p \in \{1,
\dots, n - 1\}$ such that $pk \equiv 1$ (mod $n$). Let $\varphi
\colon V(T_{n, m, k}) \to \{1, \dots, mn\}$ be given by
$\varphi(u_{i,j}) = i + mp(j - 1)$ (modulo $mn$). Since $(n, p) =
1$, $\varphi$ is a bijection. By the similar argument as in the
proof of Lemma \ref{l2.3}, $\varphi \colon T_{n, m, k} \to T_{mn,
1, mp}$ is an isomorphism.

Let $(n, k + m) = 1$ where $m \geq 3$. Let $p \in \{1, \dots, n -
1\}$ be such that $p(k+m) \equiv 1$ (mod $n$). Let $\psi \colon
V(T_{n, 2, k}) \to \{1, \dots, mn\}$ be given by $\psi(u_{i, j}) =
i + mp(j - i)$ (modulo $mn$). Since $(n, p) = 1$, $\psi$ is a
bijection. Now, $\psi(u_{i, j} u_{i, j + 1} u_{i + 1, j + 1}) =
\{i + mp(j - i), i + mp(j + 1 - i), i + 1 + mp(j - i)\}$,
$\psi(u_{i, j} u_{i + 1, j} u_{i + 1, j + 1}) = \{i + mp(j - i), i
+ 1 + mp(j - i - 1), i + 1 + mp(j - i)\}$, $\psi(u_{m, j} u_{m, j
+ 1} u_{1, j + k + 1}) = \{m + mp(j - m), m + mp(j + 1 - m), 1 +
mp(j + k)\} = \{m + mp(j - m), m+mp(j + 1 - m), 1 + m + mp(j -
m)\}$, $\psi(u_{m, j} u_{1, j+k} u_{1, j+k+1}) = \{m+mp(j-m), 1+
mp(j+k-1), 1+ mp(j+k)\} = \{m + mp(j - m), 1 + m + mp(j - m - 1),
1 + m + mp(j - m)\} \in T_{mn, 1, mp-1}$. (Clearly, if $mn$ is
even then $mn/2-1 \neq mp-1 \neq mn/2$ and if $mn$ is odd then
$mp-1 \neq (mn-1)/2$.) Thus $\psi \colon T_{n, m, k} \to T_{mn, 1,
mp - 1}$ is an isomorphism. This proves $(a)$.

The map $u_{i, 1} \mapsto u_i$, $u_{i, 3} \mapsto u_{i+4}$, $u_{i,
2} \mapsto v_i$ and $u_{i, 4} \mapsto v_{i+4}$ for $1 \leq i \leq
4$ defines an isomorphism from $T_{4, 4, 2}$ to $T_{8, 2, 2}$.
This proves $(b)$.

Since $G_1({\rm EG}(T_{16, 1, 2})) = 2C_8$, $G_1({\rm EG}(T_{16,
1, 3})) = 2C_8$, $G_1({\rm EG}(T_{16, 1, 6})) = 4C_4$ and
$G_1({\rm EG}(T_{4, 4, 0}))$ is a null graph, $T_{4, 4, 0}$ is not
isomorphic to any of $T_{16, 1, 2}$, $T_{16, 1, 3}$ and $T_{16,
1, 6})$. Again $G_4({\rm EG}(T_{4, 4, 0}))$ is a null graph
whereas $G_4({\rm EG}(T_{8, 2, 2})) = 8K_2$. Thus $T_{4, 4, 0}
\not\cong T_{8, 2, 2}$. $(c)$ now follows from parts $(a)$, $(c)$
of Lemma \ref{l2.3} and part $(g)$ of Lemma \ref{l2.1}. \hfill
$\Box$

\bigskip

Now, we will present three series of degree-regular triangulations
of the Klein bottle. Among these, $B_{m, n}$ were defined earlier
in \cite{dn}. The smallest among $Q_{m, n}$'s, namely $Q_{5,2}$,
also defined in \cite{dn} as $Q$.

\begin{eg}$\!\!${\rm {\bf :}} \label{e4} {\rm A series of degree-regular
non-orientable combinatorial 2-manifolds of Euler characteristic
$0$. For $m, n\geq 3$,}
\begin{eqnarray*} B_{m,n} &= &\{v_{i, j}v_{i+1, j}v_{i+1, j+1},
v_{i, j}v_{i, j+1}v_{i+1, j+1} ~ : ~ 1\leq i\leq n, ~ 1\leq j\leq m-1 \} \\
&& \cup \, \{v_{i, m}v_{n+2-i, 1}v_{n+1-i, 1}, v_{i,m}v_{i+1, m}
v_{n+1-i, 1} ~ : ~ 1\leq i\leq n \},
\end{eqnarray*}
{\rm where $V(B_{m,n}) = \{v_{i,j} ~ : 1 \leq i\leq n, 1\leq j\leq
m\}$. (Addition in the first subscripts are modulo $n$.) Clearly,
$B_{m, n}$ triangulates the Klein bottle and hence it is a
non-orientable combinatorial $2$-manifold. Since the degree of
each vertex is $6$, $B_{m, n}$ is degree-regular.}
\end{eg}

\begin{eg}$\!\!${\rm {\bf :}} \label{e5} {\rm A series of degree-regular
non-orientable combinatorial 2-manifolds of Euler characteristic
$0$. For $m\geq 3$ and $n\geq 2$,}
\begin{eqnarray*} K_{m, 2n} &= &\{v_{i, j}v_{i, j+1}v_{i+1, j},
v_{i, j+1}v_{i+1, j}v_{i+1, j+1} ~ : ~ 1\leq i\leq n, ~ 1\leq j\leq m-1 \} \\
&& \cup \, \{v_{i, m}v_{i+1, m}v_{2n+2-i, 1},
v_{i+1, m}v_{2n+2-i,1}v_{2n+1-i,1} ~ : ~ 1\leq i\leq n \} \\
& & \cup \, \{v_{i, j}v_{i+1, j}v_{i+1, j+1},
v_{i, j}v_{i, j+1}v_{i+1, j+1} ~ : ~ n+1\leq i\leq 2n, ~ 1\leq j\leq m-1 \} \\
&& \cup \, \{v_{i, m}v_{i+1, m}v_{2n+1-i, 1},
v_{i,m}v_{2n+2-i,1}v_{2n+1-i,1} ~ : ~ n+1\leq i\leq 2n \},
\end{eqnarray*}
{\rm where $V(K_{m, 2n}) = \{v_{i,j} ~ : 1 \leq i\leq 2n, 1\leq
j\leq m\}$. (Addition in the first subscripts are modulo $2n$.)
Clearly, $K_{m, 2n}$ triangulates the Klein bottle and hence it is
a non-orientable combinatorial $2$-manifold. Since the degree of
each vertex is $6$, $K_{m, 2n}$ is degree-regular.}
\end{eg}

\smallskip

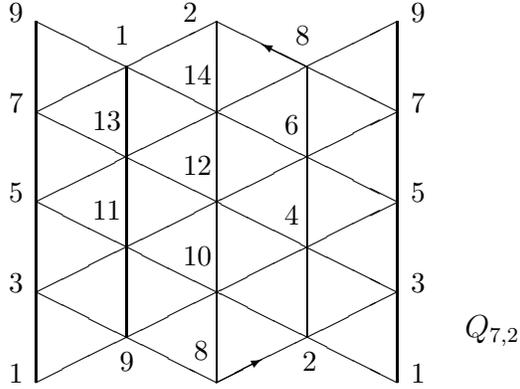
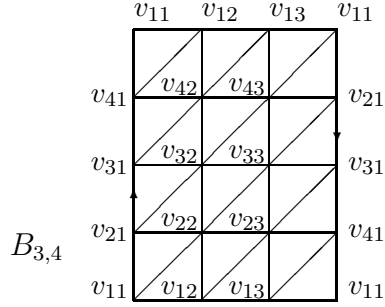
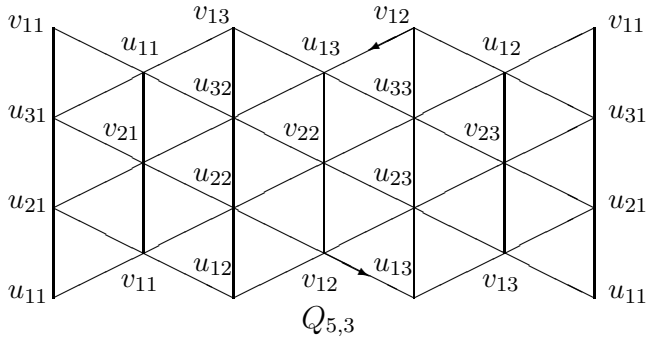
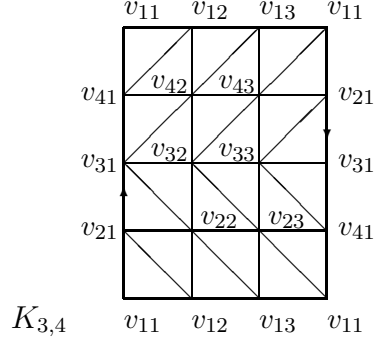
\hrule

\setlength{\unitlength}{3mm} 

\begin{picture}(45,20)(-2,-1.5)


\put(8,0){\vector(2,1){2}} \put(12,14){\vector(-2,1){2}}

\thicklines

\put(0,0){\line(0,1){16}} \put(16,0){\line(0,1){16}}

\thinlines

\put(4,2){\line(0,1){12}} \put(8,0){\line(0,1){16}}
\put(12,2){\line(0,1){12}}

\put(0,12){\line(2,1){8}} \put(0,8){\line(2,1){16}}
\put(0,4){\line(2,1){16}} \put(0,0){\line(2,1){16}}
\put(8,0){\line(2,1){8}}

\put(0,4){\line(2,-1){8}} \put(0,8){\line(2,-1){16}}
\put(0,12){\line(2,-1){16}} \put(0,16){\line(2,-1){16}}
\put(8,16){\line(2,-1){8}}

\put(-1.2,0){\mbox{$1$}} \put(-1.2,4){\mbox{$3$}}
\put(-1.2,8){\mbox{$5$}} \put(-1.2,12){\mbox{$7$}}
\put(-1.2,16){\mbox{$9$}}

\put(3.7,0.5){\mbox{$9$}} \put(2.5,7.2){\mbox{$11$}}
\put(2.5,11.2){\mbox{$13$}} \put(3.5,15){\mbox{$1$}}

\put(7,1){\mbox{$8$}} \put(6.5,5.2){\mbox{$10$}}
\put(6.5,9.2){\mbox{$12$}} \put(6.5,13.2){\mbox{$14$}}
\put(6.5,16){\mbox{$2$}}

\put(11.8,0.5){\mbox{$2$}} \put(11,7){\mbox{$4$}}
\put(11,11){\mbox{$6$}} \put(11.5,15){\mbox{$8$}}

\put(16.6,0){\mbox{$1$}} \put(16.6,4){\mbox{$3$}}
\put(16.6,8){\mbox{$5$}} \put(16.6,12){\mbox{$7$}}
\put(16.6,16){\mbox{$9$}}

{\large \put(19,2){\mbox{$Q_{7,2}$}} }

\put(33,3){\vector(0,1){2}} \put(42,9){\vector(0,-1){2}}

\thicklines

\put(33,0){\line(1,0){9}} \put(33,0){\line(0,1){12}}
\put(33,12){\line(1,0){9}} \put(42,0){\line(0,1){12}}

\thinlines

\put(33,3){\line(1,0){9}} \put(33,6){\line(1,0){9}}
 \put(33,9){\line(1,0){9}} \put(36,0){\line(0,1){12}}
\put(39,0){\line(0,1){12}}

\put(33,0){\line(1,1){9}} \put(33,3){\line(1,1){9}}
\put(33,6){\line(1,1){6}} \put(33,9){\line(1,1){3}}
\put(36,0){\line(1,1){6}} \put(39,0){\line(1,1){3}}


\put(31.1,0.3){\mbox{$v_{11}$}} \put(34.2,0.3){\mbox{$v_{12}$}}
\put(37.2,0.3){\mbox{$v_{13}$}} \put(42.5,0.3){\mbox{$v_{11}$}}

\put(31.1,2.8){\mbox{$v_{21}$}} \put(34.2,3.3){\mbox{$v_{22}$}}
\put(37.2,3.3){\mbox{$v_{23}$}} \put(42.5,2.8){\mbox{$v_{41}$}}

\put(31.1,5.8){\mbox{$v_{31}$}} \put(34.2,6.3){\mbox{$v_{32}$}}
\put(37.2,6.3){\mbox{$v_{33}$}} \put(42.5,5.8){\mbox{$v_{31}$}}

\put(31.1,8.8){\mbox{$v_{41}$}} \put(34.2,9.3){\mbox{$v_{42}$}}
\put(37.2,9.3){\mbox{$v_{43}$}} \put(42.5,8.8){\mbox{$v_{21}$}}

\put(33,12.6){\mbox{$v_{11}$}} \put(36,12.6){\mbox{$v_{12}$}}
\put(39,12.6){\mbox{$v_{13}$}} \put(42,12.6){\mbox{$v_{11}$}}

{\large \put(27.5,2){\mbox{$B_{3,4}$}} }

\end{picture}

\hrule

\bigskip

\hrule

\setlength{\unitlength}{3mm}

\begin{picture}(45,17.5)(-2,-2.5)


\put(12,2){\vector(2,-1){2}} \put(16,12){\vector(-2,-1){2}}

\thicklines

\put(0,0){\line(0,1){12}} \put(24,0){\line(0,1){12}}

\thinlines

\put(4,2){\line(0,1){8}} \put(8,0){\line(0,1){12}}
\put(12,2){\line(0,1){8}} \put(16,0){\line(0,1){12}}
\put(20,2){\line(0,1){8}}

\put(0,8){\line(2,1){8}} \put(0,4){\line(2,1){16}}
\put(0,0){\line(2,1){24}} \put(8,0){\line(2,1){16}}
\put(16,0){\line(2,1){8}}

\put(0,4){\line(2,-1){8}} \put(0,8){\line(2,-1){16}}
\put(0,12){\line(2,-1){24}} \put(8,12){\line(2,-1){16}}
\put(16,12){\line(2,-1){8}}

\put(-2,0){\mbox{$u_{11}$}} \put(-2,4){\mbox{$u_{21}$}}
\put(-2,8){\mbox{$u_{31}$}} \put(-2,12){\mbox{$v_{11}$}}

\put(3,0.5){\mbox{$v_{11}$}} \put(2.2,7.2){\mbox{$v_{21}$}}
\put(3,11){\mbox{$u_{11}$}}

\put(6.2,1.2){\mbox{$u_{12}$}} \put(6.2,5.2){\mbox{$u_{22}$}}
\put(6.2,9.2){\mbox{$u_{32}$}} \put(6.2,12.4){\mbox{$v_{13}$}}

\put(11,0.5){\mbox{$v_{12}$}} \put(10.2,7.2){\mbox{$v_{22}$}}
\put(11,11){\mbox{$u_{13}$}}

\put(14.2,1.2){\mbox{$u_{13}$}} \put(14.2,5.2){\mbox{$u_{23}$}}
\put(14.2,9.2){\mbox{$u_{33}$}} \put(14.2,12.4){\mbox{$v_{12}$}}

\put(19,0.5){\mbox{$v_{13}$}} \put(18.2,7.2){\mbox{$v_{23}$}}
\put(19,11){\mbox{$u_{12}$}}

\put(24.6,0){\mbox{$u_{11}$}} \put(24.6,4){\mbox{$u_{21}$}}
\put(24.6,8){\mbox{$u_{31}$}} \put(24.6,12){\mbox{$v_{11}$}}

{\large \put(11,-1.3){\mbox{$Q_{5,3}$}} }


\put(33,3){\vector(0,1){2}} \put(42,9){\vector(0,-1){2}}

\thicklines

\put(33,0){\line(1,0){9}} \put(33,0){\line(0,1){12}}
\put(33,12){\line(1,0){9}} \put(42,0){\line(0,1){12}}

\thinlines

\put(33,3){\line(1,0){9}} \put(33,6){\line(1,0){9}}
 \put(33,9){\line(1,0){9}} \put(36,0){\line(0,1){12}}
\put(39,0){\line(0,1){12}}

\put(36,0){\line(-1,1){3}} \put(39,0){\line(-1,1){6}}
\put(42,0){\line(-1,1){6}} \put(42,3){\line(-1,1){3}}
\put(33,6){\line(1,1){6}} \put(36,6){\line(1,1){6}}
\put(39,6){\line(1,1){3}} \put(33,9){\line(1,1){3}}


\put(33,-1.3){\mbox{$v_{11}$}} \put(36,-1.3){\mbox{$v_{12}$}}
\put(39,-1.3){\mbox{$v_{13}$}} \put(42,-1.3){\mbox{$v_{11}$}}

\put(31.1,2.8){\mbox{$v_{21}$}} \put(36.4,3.3){\mbox{$v_{22}$}}
\put(39.4,3.3){\mbox{$v_{23}$}} \put(42.5,2.8){\mbox{$v_{41}$}}

\put(31.1,5.8){\mbox{$v_{31}$}} \put(34.2,6.3){\mbox{$v_{32}$}}
\put(37.2,6.3){\mbox{$v_{33}$}} \put(42.5,5.8){\mbox{$v_{31}$}}

\put(31.1,8.8){\mbox{$v_{41}$}} \put(34.2,9.3){\mbox{$v_{42}$}}
\put(37.2,9.3){\mbox{$v_{43}$}} \put(42.5,8.8){\mbox{$v_{21}$}}

\put(33,12.6){\mbox{$v_{11}$}} \put(36,12.6){\mbox{$v_{12}$}}
\put(39,12.6){\mbox{$v_{13}$}} \put(42,12.6){\mbox{$v_{11}$}}

{\large \put(28,-1.3){\mbox{$K_{3,4}$}} }
\end{picture}

\hrule

\medskip

\begin{eg}$\!\!${\rm {\bf :}} \label{e6} {\rm A series of weakly
regular non-orientable combinatorial 2-manifolds of Euler
characteristic $0$. For each $m \geq 2$,}
\begin{eqnarray*} Q_{2m+1, 2} &= & \{\{i, i+1, i+2\}, \{i, i+2, i
+ 2m +2\}  ~ : ~ 1 \leq i \leq 4m +2\},
\end{eqnarray*}
{\rm where $V(Q_{2m+1, 2}) = \{1, \dots, 4m + 2\}$. (Addition
modulo $4m+2$.) Clearly, $Q_{2m+1, 2}$ triangulates the Klein
bottle and hence it is a non-orientable combinatorial
$2$-manifold. Since $\ZZ_{4m+2}$ acts transitively (by addition)
on vertices, $Q_{2m+1, 2}$ is weakly regular.}
\end{eg}

More generally, for each $n \geq 2$ we define the following.

\begin{eg}$\!\!${\rm {\bf :}} \label{e7} {\rm A series of degree-regular
non-orientable combinatorial 2-manifolds of Euler characteristic
$0$. For $m, n\geq 2$,}
\begin{eqnarray*} Q_{2m+1, n} &= & \{u_{i, j}u_{i+1, j}v_{i, j},
u_{i, j+1}u_{i+1, j+1}v_{i, j} ~ : ~ 1\leq i\leq m, 1\leq j\leq n\} \\
&& \cup \, \{v_{i, j}v_{i+1, j}u_{i+1, j},
v_{i, j}v_{i+1, j}u_{i+1, j+1} ~ : ~ 1\leq i\leq m-1, 1\leq j\leq n\} \\
&& \cup \, \{u_{m+1, j}u_{1, n+2-j}v_{1, n+2-j}, ~ u_{m+1,
j+1}u_{1, n+2-j}v_{1, n+1-j}, \\
&& ~~~~~~ u_{m+1,j}u_{1,n+2-j}v_{m,j}, ~ u_{m+1, j+1}u_{1,
n+2-j}v_{m, j} ~ : ~ 1\leq j\leq n\},
\end{eqnarray*}
{\rm where $V(Q_{2m+1, n}) = \{u_{i, j} ~ : ~ 1\leq i\leq m+1,
1\leq j\leq n\}\cup\{v_{i,j} ~ : 1 \leq i\leq m, 1\leq j\leq n\}$.
(Addition in the second subscripts are modulo $n$.) Clearly,
$Q_{2m+1, n}$ triangulates the Klein bottle and hence it is a
non-orientable combinatorial $2$-manifold. Since the degree of
each vertex is $6$, $Q_{2m+1, n}$ is degree-regular.

For $n \geq 3$ there are two induced 3-cycles (induced
subcomplexes which are 3-cycles) through $v_{2,1}$ in $Q_{5, n}$,
namely, $C_3(v_{2,1}, u_{1,1}, u_{2,1})$ and $C_3(v_{2,1},
u_{3,1}, v_{1,1})$. But there is no induced 3-cycle through
$v_{2,n}$. So, there does not exist any automorphism of $Q_{5, n}$
which sends $v_{2,1}$ to $v_{2,n}$. Thus, $Q_{5, n}$ is not weakly
regular for $n \geq 3$.

Observe that $G_3({\rm EG}(Q_{7, 2t-1})) = C_7(u_{1, 1}, u_{3,
1}, v_{1, 1}, v_{3, 1}, u_{2, 1}, u_{4, 1}, v_{2, 1}) \sqcup \,
C_7(u_{1, t + 1}$, \linebreak $u_{3, t + 1}, v_{1, t}, v_{3, t},
u_{2, t + 1}, u_{4, t + 1}, v_{2, t})$ and $G_3({\rm EG}(Q_{7, 2t
- 2})) = C_7(u_{1, 1}, u_{3, 1}, v_{1, 1}, v_{3, 1}, u_{2, 1},
u_{4, 1}$, \linebreak $v_{2, 1})  \sqcup C_7(u_{1, t}, u_{3, t},
v_{1, t}, v_{3, t}, u_{2, t}, u_{4, t}, v_{2, t})$ for $t \geq
2$. So, for $n \geq 2$, $G_3({\rm EG}(Q_{7, n}))$ consists of two
disjoint 7-cycles. This implies that $G_3({\rm EG}(Q_{7, n}))$ is
not regular for $n \geq 3$ and hence $Q_{7, n}$ is not weakly
regular for $n \geq 3$.

For each $m \geq 2$, there are exactly two $(2m + 1)$-cycles
(namely, $C_{2m + 1}(u_{1, 1}, \dots, u_{m + 1, 1}$, $v_{1, 1},
\dots, v_{m, 1})$, $C_{2m+1}(u_{1, t+1}, \dots, u_{m+1, t+1},
v_{1, t}, \dots, v_{m, t})$ in $Q_{2m+1, 2t-1}$ and
$C_{2m+1}(u_{1, 1}, \dots$, $u_{m+1, 1}, v_{1, 1}, \dots, v_{m,
1})$, $C_{2m+1}(u_{1, t+1}, \dots, u_{m+1, t+1}, v_{1, t+1},
\dots, v_{m, t+1})$ in $Q_{2m+1, 2t}$) each of which is the
boundary of a $(2m+1)$-vertex M\"{o}bius strip. In other words,
there are exactly two $(2m+1)$-cycle, say $C_1$ and $C_2$, such
that $|Q_{2m-1, n}|\setminus |C_i|$ is union of two disjoint open
M\"{o}bius strips for $i=1, 2$. Thus, if $n\geq 3$ and $u$ is a
vertex in $C_1$ and $v$ is a vertex outside $C_1\cup C_2$ then
there does not exist any automorphism of $Q_{2m+1, n}$ which
sends $u$ to $v$. So, $Q_{2m+1, n}$ is not weakly regular for all
$m\geq 2$ and $n\geq 3$.

}
\end{eg}

\begin{lemma}$\!\!${\bf .} \label{l2.5} Let $B_{m, n}$, $K_{m, 2n}$
and $Q_{2m+1, n}$ be as above. We have the following\,:
\begin{enumerate}
    \item[$(a)$] $B_{3, 4} \not\cong B_{4, 3} \not \cong K_{3, 4}
    \not\cong B_{3, 4}$.
    \item[$(b)$] $B_{3, 5} \not\cong B_{5, 3} \not \cong Q_{5, 3}
    \not\cong B_{3, 5}$.
    \item[$(c)$] None of $B_{3, 4}$, $B_{4, 3}$, $K_{3, 4}$,
    $B_{3, 5}$, $B_{5, 3}$, $Q_{5, 3}$ are weakly regular.
\end{enumerate}
\end{lemma}

\noindent {\bf Proof.} Observe that $G_4({\rm EG}(B_{3, 4})) =
C_3(v_{11}, v_{42}, v_{23}) \cup C_3(v_{41}, v_{12}, v_{13})\cup
C_3(v_{31}, v_{22}, v_{43}) \cup C_3(v_{21}, v_{32}, v_{33})$,
$G_4({\rm EG}(B_{4, 3})) = C_8(v_{11}, v_{32}, v_{23}, v_{14},
v_{21}, v_{12}, v_{33}, v_{24}) \cup K_4(\{v_{31}, v_{22},
v_{13}$, $v_{34}\})$ and $G_4({\rm EG}(K_{3, 4})) = K_4(\{v_{11},
v_{22}, v_{33}, v_{42}\})\cup K_4(\{v_{12}, v_{23}, v_{31},
v_{43}\}) \cup K_4(\{v_{13}, v_{21}, v_{32}$, $v_{41}\})$. These
proves $(a)$ (since $M\cong N$ implies $G_4({\rm EG}(M)) \cong
G_4({\rm EG}(N))$).

It is easy to see the following: (i) $G_4({\rm EG}(B_{5,3})) =
C_{5}(v_{21}, v_{12}, v_{33}, v_{24}, v_{15}) \cup C_{10}(v_{11},
v_{32}$, $v_{23}, v_{14}, v_{35}, v_{31}, v_{22}, v_{13}, v_{34},
v_{25})$, (ii) $G_4({\rm EG}(B_{3, 5}))$ is $C_3(v_{11}, v_{52},
v_{23}) \cup C_3(v_{31}, v_{32}, v_{43}) \cup C_3(v_{21},
v_{42}$, $v_{33}) \cup C_3(v_{51}, v_{12}, v_{13})$ together with
the three isolated vertices $v_{41}$, $v_{22}$, $v_{53}$ and (iii)
$G_4({\rm EG}(Q_{5, 3}))$ is $C_5(u_{11}, u_{21}, u_{31}, v_{11},
v_{21}) \cup C_5(u_{13}, u_{23}, u_{33}, v_{12}, v_{22})$
together with five isolated vertices. These prove $(b)$.

If $M$ is weakly regular then ${\rm Aut}(M)$ acts
vertex-transitively on $G_n({\rm EG}(M))$ for all $n\geq 0$.
Since $G_4({\rm EG}(B_{4, 3})) = C_8 \sqcup K_4$, no group can act
vertex-transitively on $G_4({\rm EG}(B_{4, 3}))$. So, $B_{4, 3}$
is not weakly regular. Similarly, $B_{5, 3}$ is not weakly
regular. Since $G_4({\rm EG}(B_{3, 5}))$ and $G_4({\rm
EG}(Q_{5,3}))$ are not regular graphs, $B_{3, 5}$ and $Q_{5, 3}$
are not weakly regular.

Observe that $G_3({\rm EG}(B_{3, 4})) \cap {\rm EG}(B_{3, 4}) =
C_6(v_{11}, v_{12}, v_{23}, v_{41}, v_{42}, v_{13}) \cup
C_6(v_{31}, v_{32}, v_{43}$, $v_{21}, v_{22}, v_{33})$. If
possible let there be $\sigma \in {\rm Aut}(B_{3, 4})$ such that
$\sigma(v_{11}) = v_{12}$. Since ${\rm Aut}(B_{3, 4})$ acts
vertex-transitively on the graph $G_3({\rm EG}(B_{3, 4})) \cap
{\rm EG}(B_{3, 4})$, $\sigma(v_{12}) = v_{11}$ or $v_{23}$. In the
first case, $\sigma(v_{41}) = v_{42}$ and hence $\sigma(v_{11}
v_{12}v_{41}) = v_{11}v_{12}v_{42}$. But $v_{11}v_{12}v_{41}$ is a
face, whereas $v_{11}v_{12}v_{42}$ is not a face, a
contradiction. In the second case, $(v_{11}, v_{12}, v_{23},
v_{41}$, $v_{42}, v_{13})$ is a cycle in (the permutation)
$\sigma$. Then $\sigma(v_{11}v_{12}v_{41}) = v_{12}v_{23}v_{42}$.
But $v_{11}v_{12}v_{41}$ is a face, whereas $v_{12}v_{23}v_{42}$
is not a face, a contradiction. So, there is no automorphism
which maps $v_{11}$ to $v_{12}$. Therefore, $B_{3, 4}$ is not
weakly regular.

If possible let there be $\tau \in {\rm Aut}(K_{3, 4})$ such that
$\tau(v_{11}) = v_{21}$. Since $G_3({\rm EG}(K_{3, 4})) =
C_3(v_{11}, v_{12}, v_{13}) \cup C_3(v_{21}, v_{23}, v_{42}) \cup
C_3(v_{31}, v_{32}, v_{33}) \cup C_3(v_{41}, v_{22}, v_{43})$,
$\tau(v_{12}) = v_{23}$ or $v_{42}$. In either case $\tau$ maps
the edge $v_{11}v_{12}$ of $K_{3, 4}$ to a non-edge of $K_{3,
4}$, a contradiction. So, there is no automorphism which maps
$v_{11}$ to $v_{21}$. Thus, $K_{3, 4}$ is not weakly regular.
\hfill $\Box$

\begin{eg}$\!\!${\rm {\bf :}} \label{e8} {\rm A triangulation
($E$) of the plane $\RR^2$. The vertex-set $V(E) = \{u_{m,
2n}=(m, n\sqrt{3}), u_{m, 2n-1} = (m + \frac{1}{2},
\frac{(2n-1)\sqrt{3}}{2}) : m, n \in \ZZ\}$ and the faces are
$\{u_{m, 2n}u_{m+1, 2n}u_{m, 2n+1}$, $u_{m+1, 2n}u_{m, 2n+1}
u_{m+1, 2n+1}, \, u_{m, 2n-1}u_{m+1, 2n-1}u_{m+1, 2n}, \, u_{m,
2n-1}u_{m, 2n} u_{m+1, 2n} : m, n \in \ZZ\}$. \linebreak The group
$H$ of translations generated by $\alpha_1 : u \mapsto u + u_{1,
0}$ and $\alpha_2 : u \mapsto u + u_{0,1}$ is a subgroup of ${\rm
Aut}(E)$. Clearly, $H$ acts transitively on $V(E)$. The
stabilizer of any vertex $u$ in ${\rm Aut}(E)$ is isomorphic to
the dihedral group $D_6$ (of order 12) which acts transitively on
the set of flags containing $u$. So, $E$ is combinatorially
regular. Let $G_0$ denote the stabilizer of $u_{0, 0}$. Since $H$
acts transitively on $V(E)$, ${\rm Aut}(E) = \langle H, G_0
\rangle$. This implies that if $\sigma \in {\rm Aut}(E)$ has no
fixed element in $E$ (vertex, edge or face) then either $\sigma
\in H \setminus \{{\rm Id}\}$ or is a glide reflection (i.e., an
automorphism of the form $t_a \circ r_l$, where $r_l \in {\rm
Aut}(E)$ is the reflection about a line $l$ through some vertex
and of slope a multiple of $\pi/6$ and $t_a \in H$ is the
translation by a nonzero vector $a$ parallel to $l$). }
\end{eg}

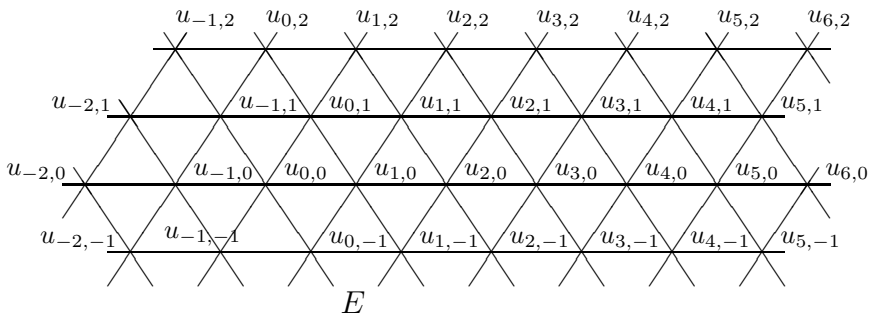
\hrule

\setlength{\unitlength}{3mm}

\begin{picture}(45,15)(-3,-0.2)


\thicklines

\put(0,6){\line(1,0){34}}

\thinlines

\put(2,3){\line(1,0){30}} \put(2,9){\line(1,0){30}}
\put(4,12){\line(1,0){30}}

\put(0,4.5){\line(2,3){5.5}} \put(2,1.5){\line(2,3){7.5}}
\put(6,1.5){\line(2,3){7.5}} \put(10,1.5){\line(2,3){7.5}}
\put(14,1.5){\line(2,3){7.5}} \put(18,1.5){\line(2,3){7.5}}
\put(22,1.5){\line(2,3){7.5}} \put(26,1.5){\line(2,3){7.5}}
\put(30,1.5){\line(2,3){3.5}}

\put(4,1.5){\line(-2,3){3.5}} \put(8,1.5){\line(-2,3){5.5}}
\put(12,1.5){\line(-2,3){7.5}} \put(16,1.5){\line(-2,3){7.5}}
\put(20,1.5){\line(-2,3){7.5}} \put(24,1.5){\line(-2,3){7.5}}
\put(28,1.5){\line(-2,3){7.5}} \put(32,1.5){\line(-2,3){7.5}}
\put(34,4.5){\line(-2,3){5.5}} \put(34,10.5){\line(-2,3){1.5}}



\put(-1,3.4){\mbox{\small $u_{-2,-1}$}}
 \put(4.5,3.6){\mbox{\small $u_{-1,-1}$}}
 \put(11.8,3.4){\mbox{\small $u_{0,-1}$}}
 \put(15.8,3.4){\mbox{\small $u_{1,-1}$}}
 \put(19.8,3.4){\mbox{\small $u_{2,-1}$}}
 \put(23.8,3.4){\mbox{\small $u_{3,-1}$}}
 \put(27.8,3.4){\mbox{\small $u_{4,-1}$}}
 \put(31.8,3.4){\mbox{\small $u_{5,-1}$}}

 \put(-2.5,6.4){\mbox{\small $u_{-2,0}$}}
 \put(5.8,6.4){\mbox{\small $u_{-1,0}$}}
 \put(9.8,6.4){\mbox{\small $u_{0,0}$}}
 \put(13.8,6.4){\mbox{\small $u_{1,0}$}}
 \put(17.8,6.4){\mbox{\small $u_{2,0}$}}
 \put(21.8,6.4){\mbox{\small $u_{3,0}$}}
 \put(25.8,6.4){\mbox{\small $u_{4,0}$}}
 \put(29.8,6.4){\mbox{\small $u_{5,0}$}}
 \put(33.8,6.4){\mbox{\small $u_{6,0}$}}

 \put(-0.5,9.4){\mbox{\small $u_{-2,1}$}}
 \put(7.8,9.4){\mbox{\small $u_{-1,1}$}}
 \put(11.8,9.4){\mbox{\small $u_{0,1}$}}
 \put(15.8,9.4){\mbox{\small $u_{1,1}$}}
 \put(19.8,9.4){\mbox{\small $u_{2,1}$}}
 \put(23.8,9.4){\mbox{\small $u_{3,1}$}}
 \put(27.8,9.4){\mbox{\small $u_{4,1}$}}
 \put(31.8,9.4){\mbox{\small $u_{5,1}$}}

\put(5,13.2){\mbox{\small $u_{-1,2}$}}
 \put(9,13.2){\mbox{\small $u_{0,2}$}}
 \put(13,13.2){\mbox{\small $u_{1,2}$}}
 \put(17,13.2){\mbox{\small $u_{2,2}$}}
 \put(21,13.2){\mbox{\small $u_{3,2}$}}
 \put(25,13.2){\mbox{\small $u_{4,2}$}}
 \put(29,13.2){\mbox{\small $u_{5,2}$}}
 \put(33,13.2){\mbox{\small $u_{6,2}$}}


{\large \put(12.3,0.3){\mbox{$E$}} }

\end{picture}

\hrule

\section{Proofs.}

\begin{lemma}$\!\!${\bf .} \label{disc}
There is no triangulation of the closed $2$-disk such that $(i)$
the degree of each vertex $($except one$)$ in the boundary is $4$
and $(ii)$  the degree of each interior vertex is $6$.
\end{lemma}

\noindent {\bf Proof.} If possible let there be a triangulation
$K$ of the closed 2-disk on $m+n+1$ vertices with $n$ interior
vertices such that the degree of each interior vertex is 6, the
degree of one vertex in the boundary is $k$ ($\geq 2$) and the
degree of each of the remaining $m$ vertices in the boundary is
4. Then $f_0(K) = n + m + 1$, $f_1(K) = \frac{6n + 4m + k}{2}$ and
$f_2(K) = \frac{6n + 3m + k - 1}{3}$. Therefore, $1 = \chi(K) =
f_0(K) - f_1(K) + f_2(K) = n + m +1 - (3n + 2m + k/2) + (2n + m +
(k - 1)/3)$. This implies that $k = - 2$, a contradiction. This
proves the lemma \hfill $\Box$

\begin{lemma}$\!\!${\bf .} \label{plane}
Let $E$ be as in Example $8$ and let $M$ be a triangulation of the
plane $\RR^2$. If the degree of each vertex of $M$ is $6$ then $M$
is isomorphic to $E$.
\end{lemma}

\noindent {\bf Proof.} Choose an edge, say $v_{0, 0}v_{1, 0}$.
Then there exists a unique vertex, say $v_{2, 0}$, in ${\rm
lk}(v_{1, 0})$ such that each side of the segment $v_{0, 0} v_{1,
0} v_{2, 0}$ (union of two line segments) contains three faces
from ${\rm st}(v_{1, 0})$ (i.e., ${\rm lk}(v_{1, 0})$ is of the
form $C_6(v_{0, 0}, x, y, v_{2, 0}, z, w)$). Now, given $v_{1,
0}$ and $v_{2, 0}$ there exists unique vertex $v_{3, 0}$ in ${\rm
lk}(v_{2, 0})$ such that each side of the segment $v_{1, 0} v_{2,
0} v_{3, 0}$ contains three faces from ${\rm st}(v_{2, 0})$.
Similarly, given $v_{1, 0}$ and $v_{0, 0}$ there exists unique
vertex  $v_{-1, 0}$ in ${\rm lk}(v_{0, 0})$ such that each side of
the segment $v_{1, 0} v_{0, 0} v_{-1, 0}$ contains three faces
from ${\rm st}(v_{0, 0})$. Continuing this way we get vertices
$v_{i, 0}$, $i\in \ZZ$, such that each side of the segment
$v_{i-1, 0} v_{i, 0} v_{i+1, 0}$ contains three faces from ${\rm
st}(v_{i, 0})$. Because of Lemma \ref{disc}, all these vertices
are distinct. So, we get a triangulation of a line (see the
figure).

\bigskip

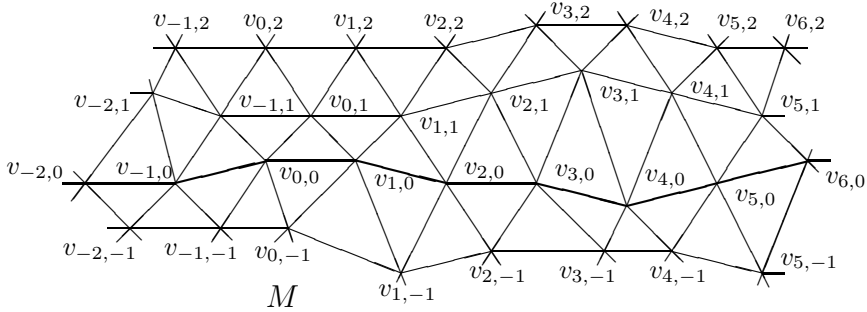
\hrule

\setlength{\unitlength}{3mm}

\begin{picture}(45,15)(-3,-0.2)


\thicklines

\put(0,6){\line(1,0){5}} \put(5,6){\line(4,1){4}}
\put(9,7){\line(1,0){4}} \put(13,7){\line(4,-1){4}}
\put(17,6){\line(1,0){4}} \put(21,6){\line(4,-1){4}}
\put(25,5){\line(4,1){8}} \put(33,7){\line(1,0){1}}

\thinlines

\put(2,4){\line(1,0){8}} \put(10,4){\line(5,-2){5}}
\put(15,2){\line(4,1){4}} \put(19,3){\line(1,0){8}}
\put(27,3){\line(4,-1){4}} \put(31,2){\line(1,0){1}}

\put(3,10){\line(1,0){1}} \put(4,10){\line(3,-1){3}}
\put(7,9){\line(1,0){8}} \put(15,9){\line(4,1){8}}
\put(23,11){\line(4,-1){8}} \put(31,9){\line(1,0){1}}

\put(4,12){\line(1,0){13}} \put(17,12){\line(4,1){4}}
\put(21,13){\line(1,0){4}} \put(25,13){\line(4,-1){4}}
\put(29,12){\line(1,0){4}}

\put(4,10){\line(-3,-4){3.6}} \put(4,10){\line(1,2){1.3}}
\put(5,6){\line(2,3){4.4}} \put(5,6){\line(-1,-1){2.6}}
\put(11,9){\line(2,3){2.4}} \put(9,7){\line(1,1){2}}
\put(9,7){\line(-2,-3){2.4}} \put(15,9){\line(-1,-1){5.6}}
\put(15,9){\line(2,3){2.4}} \put(19,10){\line(-1,-2){4.2}}
\put(19,10){\line(2,3){2.4}} \put(21,6){\line(-2,-3){2.4}}
\put(21,6){\line(2,5){2}} \put(23,11){\line(1,1){2.5}}
\put(25,5){\line(2,5){2}} \put(25,5){\line(-1,-2){1.3}}
\put(27,10){\line(1,1){2.5}} \put(31,9){\line(-2,-3){4.3}}
\put(31,9){\line(1,3){1.2}} \put(31,2){\line(2,5){2.2}}
\put(33,7){\line(-2,-5){2.2}}

\put(3.5,3.5){\line(-1,1){3}} \put(5,6){\line(1,-1){2.5}}
\put(5,6){\line(-1,4){1.15}} \put(9,7){\line(-1,1){2}}
\put(9,7){\line(1,-3){1.2}} \put(7,9){\line(-2,3){2.4}}
\put(13,7){\line(-1,1){2}} \put(13,7){\line(2,-5){2.2}}
\put(11,9){\line(-2,3){2.4}} \put(17,6){\line(-2,3){4.4}}
\put(17,6){\line(2,-3){2.4}} \put(21,6){\line(1,-1){3.5}}
\put(21,6){\line(-1,2){2}} \put(19,10){\line(-1,1){2.5}}
\put(25,5){\line(1,-1){2.5}} \put(25,5){\line(-1,3){2}}
\put(23,11){\line(-1,1){2.5}} \put(27,10){\line(1,-2){4.3}}
\put(27,10){\line(-2,3){2.4}} \put(31,9){\line(1,-1){2.5}}
\put(31,9){\line(-2,3){2.4}} \put(31.5,12.5){\line(1,-1){1.2}}


\put(0,2.8){\mbox{\small $v_{-2,-1}$}}
 \put(4.5,2.8){\mbox{\small $v_{-1,-1}$}}
 \put(8.5,2.7){\mbox{\small $v_{0,-1}$}}
 \put(14,1){\mbox{\small $v_{1,-1}$}}
 \put(18,1.8){\mbox{\small $v_{2,-1}$}}
 \put(22,1.8){\mbox{\small $v_{3,-1}$}}
 \put(26,1.8){\mbox{\small $v_{4,-1}$}}
 \put(31.8,2.4){\mbox{\small $v_{5,-1}$}}

 \put(-2.5,6.4){\mbox{\small $v_{-2,0}$}}
 \put(2.35,6.4){\mbox{\small $v_{-1,0}$}}
 \put(9.6,6.2){\mbox{\small $v_{0,0}$}}
 \put(13.8,5.9){\mbox{\small $v_{1,0}$}}
 \put(17.8,6.4){\mbox{\small $v_{2,0}$}}
 \put(21.8,6.3){\mbox{\small $v_{3,0}$}}
 \put(25.8,6.1){\mbox{\small $v_{4,0}$}}
 \put(29.8,5.2){\mbox{\small $v_{5,0}$}}
 \put(33.8,6){\mbox{\small $v_{6,0}$}}

 \put(0.5,9.4){\mbox{\small $v_{-2,1}$}}
 \put(7.8,9.4){\mbox{\small $v_{-1,1}$}}
 \put(11.8,9.4){\mbox{\small $v_{0,1}$}}
 \put(15.8,8.4){\mbox{\small $v_{1,1}$}}
 \put(19.8,9.4){\mbox{\small $v_{2,1}$}}
 \put(23.8,9.8){\mbox{\small $v_{3,1}$}}
 \put(27.8,10){\mbox{\small $v_{4,1}$}}
 \put(31.8,9.4){\mbox{\small $v_{5,1}$}}

\put(4,13){\mbox{\small $v_{-1,2}$}}
 \put(8,13){\mbox{\small $v_{0,2}$}}
 \put(12,13){\mbox{\small $v_{1,2}$}}
 \put(16,13){\mbox{\small $v_{2,2}$}}
 \put(21.6,13.5){\mbox{\small $v_{3,2}$}}
 \put(26,13.2){\mbox{\small $v_{4,2}$}}
 \put(29,13){\mbox{\small $v_{5,2}$}}
 \put(32,13){\mbox{\small $v_{6,2}$}}


{\large \put(9,0.5){\mbox{$M$}} }

\end{picture}

\hrule

\bigskip

Let ${\rm lk}(v_{1, 0}) = C_6(v_{0, 0}, v_{0, 1}, v_{1, 1}, v_{2,
0}, v_{1, -1}, v_{0, -1})$. By the same argument as above, there
exists a unique vertex, say $v_{2, 1}$, in ${\rm lk}(v_{1, 1})$
such that each side of the segment $v_{0, 1} v_{1, 1} v_{2, 1}$
contains three faces from ${\rm st}(v_{1, 1})$. This implies that
${\rm lk}(v_{2, 0})$ is of the form $C_6(v_{1, 0}, v_{1, 1}, v_{2,
1}, v_{3, 0}$, $x, v_{1, -1})$. If we continue this way we get
vertices $v_{i, 1}$, $i \in \ZZ$, such that each side of the
segment $v_{i - 1, 1} v_{i, 1} v_{i + 1, 1}$ contains three faces
from ${\rm st}(v_{i, 1})$ and $v_{i, 0}v_{i + 1, 0}v_{i, 1}$,
$v_{i + 1,0} v_{i, 1}v_{i + 1,1}$ are faces for all $i \in \ZZ$.

Similarly, we get: (i) vertices $v_{i, 2}$, $i \in \ZZ$, such
that each side of the segment $v_{i - 1, 2} v_{i, 2} v_{i + 1,
2}$ contains three faces from ${\rm st}(v_{i, 2})$ and $v_{i, 1}
v_{i + 1, 1}v_{i+1, 2}$, $v_{i, 1} v_{i, 2} v_{i + 1, 2}$ are
faces for all $i \in \ZZ$. (ii) vertices $v_{i, -1}$, $i \in
\ZZ$, such that each side of the segment $v_{i - 1, -1} v_{i, -1}
v_{i + 1, -1}$ contains three faces from ${\rm st}(v_{i, -1})$
and $v_{i, -1} v_{i + 1, -1}v_{i+1, 0}$, $v_{i, -1} v_{i, 0} v_{i
+ 1, 0}$ are faces for all $i \in \ZZ$.

Continuing this way we get vertices $v_{i, j}$, $i, j \in \ZZ$, of
$M$ such that each side of the segment $v_{i - 1, j} v_{i, j} v_{i
+ 1, j}$ contains three faces from ${\rm st}(v_{i, j})$ and $v_{i,
2k} v_{i + 1, 2k} v_{i, 2k + 1}$, $v_{i + 1, 2k} v_{i, 2k + 1}
v_{i + 1, 2k + 1}$, $v_{i, 2k + 1} v_{i + 1, 2k + 1}v_{i + 1, 2k
+ 2}$, $v_{i, 2k + 1} v_{i, 2k + 2} v_{i + 1, 2k + 2}$ are faces
for all $i, j, k \in \ZZ$. Since $M$ is connected $\{v_{i, j} : i,
j \in \ZZ\}$ is the vertex-set of $M$. Then $\varphi \colon V(M)
\to V(E)$, given by $\varphi(v_{i, j}) = u_{i, j}$, is an
isomorphism. This proves the lemma. \hfill $\Box$

\bigskip

\noindent {\bf Proof of Theorem \ref{t1}.} Let $K$ be a
degree-regular triangulation of the torus. Since $\RR^2$ is the
universal cover of the torus, there exists a triangulation $M$ of
$\RR^2$ and a simplicial covering map $\eta \colon M \to K$ (cf.
\cite[Page 144]{sp}). Since the degree of each vertex in $K$ is 6,
the degree of each vertex in $M$ is 6. Because of Lemma
\ref{plane}, we may assume that $M = E$.

Let $\Gamma$ be the group of covering transformations. Then $|K|
= |E|/\Gamma$. For $\sigma\in \Gamma$, $\eta \circ \sigma = \eta$.
So, $\sigma$ maps the geometric carrier of a simplex to the
geometric carrier of a simplex. This implies that $\sigma$
induces an automorphism $\sigma$ of $E$. Thus, we can identify
$\Gamma$ with a subgroup of ${\rm Aut}(E)$. So, $K$ is a quotient
of $E$ by a subgroup $\Gamma$ of ${\rm Aut}(E)$, where $\Gamma$
has no fixed element (vertex, edge or face). Hence $\Gamma$
consists of translations and glide reflections. Since $K =
E/\Gamma$ is orientable, $\Gamma$ does not contain any glide
reflection. Thus $\Gamma \leq H$ (the group of translations). Now
$H$ is commutative. So, $\Gamma$ is a normal subgroup of $H$.
Since $H$ acts transitively on $V(E)$, $H/\Gamma$ acts
transitively on the vertices of $E/\Gamma$. Thus, $K$ is weakly
regular. \hfill $\Box$

\begin{lemma}$\!\!${\bf .} \label{prime}
For a prime $n \geq 7$, if $M$ is an $n$-vertex weakly regular
combinatorial $2$-manifold of Euler characteristic 0 then $M$ is
isomorphic to $T_{n, 1, k}$ for some $k$.
\end{lemma}

\noindent {\bf Proof.}  Since $M$ is weakly regular, it is
degree-regular. Let $d$ be the degree of each vertex. Then $nd =
2f_1(M)$ and $n- f_1(M) + f_2(M) = \chi(M) = 0$. Since each edge
is in two triangles, $2f_1(M) = 3f_2(M)$. These imply that $d=6$.

Let $G = {\rm Aut}(M)$. Then $G$ is isomorphic to a subgroup of
the permutation group $S_n$ and hence $G$ is a finite group. Fix a
vertex $u$ of $M$. Let $H$ be the stabilizer of $u$ in $G$. Since
$M$ is weakly regular, the orbit of $u$ under the action of $G$
contains all the $n$ vertices and hence the index of $H$ in $G$ is
$n$. Thus, $n$ divides the order of $G$. Since $n$ is prime, $G$
has an element, say $\tau$, of order $n$.

Let $v$ be a vertex in $M$ such that $\tau(v)\neq v$. Then $V(M) =
\{v, \tau(v), \dots, \tau^{n-1}(v)\}$. Choose an edge $e$
containing $v$. Let $e=v\tau^k(v)$. Then $\sigma = \tau^k$ is
again an automorphism of order $n$ and $V(M) = \{v_0=v,
v_1=\sigma(v), \dots, v_{n-1}= \sigma^{n-1}(v)\}$. For each $i=0,
\dots, n-1$, $\sigma^i$ is an automorphism. Thus, $v_kv_l$ is an
edge implies $v_{k+i}v_{l+i}$ is an edge and $v_kv_lv_j$ is a face
implies $v_{k+i}v_{l+i}v_{j+i}$ is a face for each $i$. Since
$v_0v_1$ is an edge, $v_iv_{i+1}$ is an edge for each $i$.
(Addition in the subscripts are modulo $n$.)

\smallskip

\noindent {\bf Claim. }  $v_0v_1v_2$, $v_0v_1v_{n-1}$ and
$v_0v_1v_{\frac{n +1}{2}}$ are not faces.

If $v_0v_1v_2$ is a face then $v_{n-2}v_{n-1}v_0$, $v_{n-1}v_0v_1$
are faces. Let $v_0v_2v_i$ ($\neq v_0v_1v_2$) be the second face
containing $v_0v_2$. Then $v_{n-2}v_0v_{i-2}$ is a face and hence
${\rm lk}(v_0)= C_6(v_i, v_2, v_1, v_{n-1}$, $v_{n-2}, v_{i-2})$.
Then $v_{i-2}v_iv_0$ is a face and hence $v_0v_2v_{n-i+2}$ is a
face. This implies that $n-i+2 =1$ or $i$. Since $i\neq 1$, $n-i+2
= i$. Then $n = 2i - 2$. This is not possible since $n$ is a
prime. So, $v_0v_1v_2$ is not a face and hence $v_0v_1v_{n-1}$ is
not a face.

Let $c= \frac{n+1}{2}$. If $v_0v_1v_c$ is a face then $v_{n-
1}v_{0}v_{c-1}$, $v_{c-1}v_cv_0$ are faces. Let $v_0v_1v_i$ ($\neq
v_0v_1v_c$) be the second face containing $v_0v_1$. Then
$v_{n-1}v_0v_{i-1}$ is a face and hence ${\rm lk}(v_0)= C_6(v_i,
v_1, v_c, v_{c-1}$, $v_{n-1}, v_{i-1})$. Then $v_{i-1}v_iv_0$ is a
face and hence $v_0v_1v_{n-i+1}$ is a face. This implies that
$n-i+1 =c$ or $i$. In either case, we get $i = c$. This is not
possible since $v_0v_1v_i \neq v_0v_1v_c$. So, $v_0v_1v_c$ is not
a face, where $c= \frac{n+1}{2}$. This proves the claim.

Let $v_0v_1v_k$ be a face containing $v_0v_1$. Then, by the claim,
$k\in \{3, \dots, \frac{n-1}{2}, \frac{n+3}{2}, \dots, n-2\}$.
Now, $v_0v_1v_k\in M$ implies $v_{n-1}v_0v_{k-1}$,
$v_{n-k}v_{n-k+1}v_0\in M$. Then $V({\rm lk}(v_0)) = \{v_{k- 1}$,
$v_{n-1}, v_{n-k}, v_{n - k + 1}, v_1, v_{k}\}$ and hence $V({\rm
lk}(v_1)) = \{v_{k}, v_{0}, v_{n-k+1}, v_{n - k+2}, v_2,
v_{k+1}\}$. Thus, $V({\rm lk}(v_0)) \cap V({\rm lk}(v_1)) = \{v_k,
v_{n-k+1}\}$ and hence $v_0v_1v_{n-k+1}\in M$. This gives $v_{n-
1}v_0v_{n- k}\in M$. Then ${\rm lk}(v_0) = C_6(v_{k- 1}, v_{n-1},
v_{n-k}, v_{n - k + 1}, v_1, v_{k})$ and hence ${\rm lk}(v_i) =
C_6(v_{i+k-1}, v_{n+i-1}$, $v_{n+i-k}, v_{n+i-k+1}, v_{i+1}, v_{i
+ k})$ for all $i$. Now, $M \cong T_{n, 1, k - 1}$ by the map
$\varphi\colon V(M)\to \{1, \dots, n\}$ given by $\varphi(v_j) =
j$ for $1\leq i\leq n-1$ and $\varphi(v_0) = n$. \hfill $\Box$

\begin{lemma} $\!\!${\bf .} \label{23x24} Let $G$ be a group of
order $24\times 23$. Then $G$ has a unique $($and hence normal$)$
subgroup of order $23$.
\end{lemma}

\noindent {\bf Proof.} Clearly, the number of Sylow 23-subgroups
of $G$ is 24 or 1. If possible let there be 24 Sylow 23-subgroups.
Let $H$ be a Sylow 23-subgroup. Let $N(H)$ be the normalizer of
$H$ in $G$. Since all the Sylow 23-subgroups are conjugates of $H$
and $|G| = |N(H)| \times |\{\mbox{conjugates of $H$}\}|$, $|N(H)|
= 23$ and hence $N(H) = H$.

Let $A$ be the set of Sylow 3-subgroups. Then $H$ acts on $A$ by
conjugation. Since there is a Sylow 2-subgroup, the number of
elements of order 3 is at most 16 and hence $\#(A) < 23$. This
implies that the action of $H$ on $A$ is trivial. Let $K \in A$.
Then $xKx^{-1} = K$ for all $x\in H$. So, $H$ acts on $K$ by
conjugation. Since $|K| = 3$, this action of $H$ on $K$ is
trivial. So, $xy = yx$ for all $x \in H$ and $y \in K$. This
implies that $K \subseteq N(H)$. This is a contradiction since
$N(H) = H$. This proves the lemma. \hfill $\Box$

\begin{lemma} $\!\!${\bf .} \label{flag}
Let $M$ be a connected combinatorial $2$-manifold. Then the
number of flags in $M$ is divisible by the order of ${\rm
Aut}(M)$.
\end{lemma}

\noindent {\bf Proof.} Let $G = {\rm Aut}(M)$ and let ${\cal F}$
denote the set of flags of $M$. Then $G$ acts on ${\cal F}$. Let
$\sigma \in G$. If there exists a flag $F = (u, uv, uvw)$ such
that $\sigma(F) = F$ then $\sigma(v) = v$ and $\sigma(w) = w$.
This implies that $\sigma|_{{\rm lk}(u)} \equiv {\rm Id}$. Since
$M$ is connected, this implies that $\sigma \equiv {\rm Id}$.
Thus, no element of ${\cal F}$ is fixed by a non-identity element
of $G$. Therefore, the length of each orbit in ${\cal F}$ is same
as the order of $G$. This proves the lemma. \hfill $\Box$

\begin{lemma} $\!\!${\bf .} \label{2p} Let $K$ be a degree-regular
triangulation of the torus on $2p$ vertices. If $p$ is prime and
$\geq 13$ then ${\rm Aut}(K)$ has a normal subgroup of order $p$.
\end{lemma}

\noindent {\bf Proof.} Let $G= {\rm Aut}(K)$. By Lemma
\ref{flag}, $|G|$ is a factor of $2p \times 6 \times 2 =24p$. By
Theorem 1, $G$ acts transitively on $V(K)$. So, the index of the
stabilizer of a vertex is $2p$. Thus, $2p$ (and hence $p$)
divides the order of $G$. Since $p$ is prime, $G$ has an element,
say $\sigma$, of order $p$.

Since $|G|$ is a factor of $24p$, by Sylow's theorem, $G$ has a
unique Sylow $p$-subgroup for $p = 13, 17, 19$ or $p > 23$. If $p
= 23$ and $|G|< 24p$ then, by Sylow's theorem, $G$ has a unique
Sylow $p$-subgroup. Finally, if $p = 23$ and $|G| = 24 \times 23$
then, by Lemma \ref{23x24}, $G$ has a unique Sylow $23$-subgroup.
Therefore, $H = \langle \sigma \rangle$ is the unique (and hence
normal) subgroup of order $p$ in $G$. \hfill $\Box$

\bigskip

\noindent {\bf Proof of Theorem \ref{t2}.} Let $n \geq 9$ be a
composite number. Then either $n = mk$ for some $m, k \geq 3$ or
$n = 2p$ for some prime $p \geq 5$. For $m, k \geq 3$, $B_{m, k}$
(defined in Example \ref{e4}) is an $(mk)$-vertex degree-regular
triangulation of the Klein bottle. If $p \geq 5$ is a prime then
$Q_{p, 2}$ (defined in Example \ref{e6}) is an $(2p)$-vertex
degree-regular triangulation of the Klein bottle.

Let $p \geq 13$ be a prime. If possible let there be a $p$-vertex
degree-regular triangulation $X$ of the Klein bottle. Since the
torus is an orientable double cover of the Klein bottle, there
exists a $(2p)$-vertex degree-regular triangulation $K$ of the
torus and a simplicial covering map $\eta : K \to X$. Then $X$ is
a quotient of $K$ by a subgroup $\langle \tau \rangle$ of ${\rm
Aut}(K)$, where $\tau$ is an automorphism of order 2 without a
fixed element (vertex, edge or face). Then $u \tau(u)$ is a
non-edge for each $u \in V(K)$. If there exist $u, v\in V(K)$
such that $uv$ and $u \tau(v)$ are edges in $K$ then
$\deg_X(\eta(v)) <6$, a contradiction. So, $u$ and $\tau(u)$ are
not adjacent to a common vertex for all $u\in V(K)$.

By Lemma \ref{2p}, there exists a normal subgroup $H \leq {\rm
Aut}(K)$ of order $p$. Let $H = \langle \sigma \rangle$. Then
$\langle \sigma, \tau \rangle = \langle \sigma \rangle \langle
\tau \rangle$ is subgroup of order $2p$. If $\sigma \circ \tau =
\tau \circ \sigma$ then $\langle \sigma \rangle$ acts on
$K/\langle \tau \rangle$ non-trivially. This implies that $X$ is
weakly regular. But, this is not possible by Lemma \ref{prime}.
So, $\sigma \circ \tau \neq \tau \circ \sigma$ and hence $\langle
\sigma, \tau \rangle \cong D_p$.

\smallskip

\noindent {\bf Claim 1.} No vertex is fixed by $\sigma$.

If possible let $\sigma$ has a fixed vertex. Since $K$ is
connected, there is an edge of the form $uv$ such that $\sigma(u)
= u$ and $\sigma(v) \neq v$. This implies that $u\sigma^i(v)$ is
an edge for all $i$. Then $\deg(u) \geq p > 6$, a contradiction.
This proves the claim.

\smallskip

\noindent {\bf Claim 2.}  There exists $w \in V(K)$ and $i \neq 0$
such that $w \sigma^i(w)$ is an edge.

By Claim 1, $\sigma$ can be written as $\sigma = (u, \sigma(u),
\dots, \sigma^{p-1}(u))(v, \sigma(v), \dots, \sigma^{p-1}(v))$ (a
permutations on $V(K)$). If $u\sigma^i(u)$ is a non-edge for all
$i$ then the link of $u$ is of the form $C_6(\sigma^{i_1}(v),
\dots, \sigma^{i_6}(v))$. Then $\sigma^{i_1}(v)\sigma^{i_2}(v)$ is
an edge. This implies $v\sigma^{i_2-i_1}(v)$ is an edge. This
proves the claim.

By Claim 2, there exists $i\neq 0$ and $w_0\in V(K)$ such that
$w_0\sigma^i(w_0)$ is an edge. Let $\alpha = \sigma^i$. Then $w_0
\alpha(w_0)$ is an edge and hence $\alpha^{j-1}(w_0)
\alpha^j(w_0)$ is an edge for all $j$. Since $p$ is prime,
$\langle \alpha, \tau \rangle = \langle \sigma, \tau \rangle
\cong D_p$. Then $\alpha^j \circ \tau = \tau \circ \alpha^{p -
j}$ for all $j$.

Since $p$ is odd there exists $i_0$ such that
$\tau(\alpha^{i_0}(w_0)) \neq \alpha^j(w_0)$ for any $j$. Let $u_0
= \alpha^{i_0}(w_0)$, $v_0 =\tau(u_0)$, $u_i =\alpha^i(u_0)$ and
$v_i = \alpha^i(v_0)$ for $1\leq i\leq p-1$. Therefore, $\alpha =
(u_0, u_1, \dots, u_{p-1})(v_0, v_1, \dots, v_{p-1})$. Then
$\tau(u_i) = \tau(\alpha^i(u_0)) = \alpha^{p-i}(\tau(u_0)) =
\alpha^{p-i}(v_0) = v_{p-i}$ and $\tau(v_i) = \tau(\alpha^i(v_0))
= \alpha^{p - i}(\tau(v_0)) = \alpha^{p-i}(u_0) = u_{p-i}$.

Since $K$ is connected, there exists an edge of the form $u_iv_j$
and hence there exists an edge of the form $u_0v_k$ for some $k
\in \{0, \dots, p-1\}$. If $k$ is odd, let $l = \frac{p-k}{2}$.
Then $\alpha^{\,l}(u_0v_k)$ is an edge. But, $\alpha^{\,l}(u_0v_k)
= u_lv_{k+l} = u_lv_{p-l} = u_l \tau(u_l)$ is a non-edge, a
contradiction. If $k$ is even, let $m = \frac{p-k-1}{2}$. Then
$\alpha^{\,m}(u_0v_k)$ is an edge. But, $\alpha^{\,m}(u_0v_k) =
u_mv_{k+m} = u_mv_{p-m-1} = u_m \tau(u_{m+1})$. This is not
possible since $u_mu_{m+1}$ is an edge. This proves that there is
no $p$-vertex degree regular triangulation of the Klein bottle
for a prime $p \geq 13$.

If $n = 7, 8$ or 11 then, by Proposition \ref{datta2}, there does
not exist any $n$-vertex degree regular triangulation of the Klein
bottle. This completes the proof. \hfill $\Box$

\bigskip

\noindent {\bf Proof of Theorem \ref{t3}.} Since $T_{n, 1, k}$ is
a weakly regular orientable combinatorial 2-manifold of Euler
characteristic 0, Part $(a)$ follows from Part $(b)$ of Lemma
\ref{l2.1} and Part $(b)$ follows from Parts $(b)$, $(c)$, $(i)$
and $(j)$ of Lemma \ref{l2.1}. Part $(c)$ follows from Example 6.
\hfill $\Box$

\bigskip

\noindent {\bf Proof of Theorem \ref{t4}.} Follows from Lemma
\ref{prime}. \hfill $\Box$

\bigskip

\noindent {\bf Proof of Corollary \ref{t5}.} Let $M$ be an
$n$-vertex degree-regular combinatorial 2-manifold of Euler
characteristic 0. If $n$ is prime then by Theorem 2, $M$
triangulates the torus and hence, by Theorem 1, $M$ is weakly
regular. Then, by Theorem 4, $M$ is isomorphic to $T_{n, 1, k}$
for some $k$. Now, Part $(a)$ follows from Parts $(a)$, $(b)$,
$(e)$ and $(h)$ of Lemma \ref{l2.1} and Part $(b)$ follows from
Parts $(a)$, $(b)$ and $(j)$ of Lemma \ref{l2.1}. \hfill $\Box$

\begin{lemma} $\!\!${\bf .} \label{l3.7} Let $M$ be a
combinatorial $2$-manifold and $a_1, \dots, a_5$ be five vertices.
If the degree of each vertex is $6$ then the number of faces in
${\rm st}(a_1) \cup {\rm st}(a_2) \cup {\rm st}(a_3)$ is $\geq
12$ and the number of faces in ${\rm st}(a_1) \cup \cdots \cup
{\rm st}(a_5)$ is $ > 12$.
\end{lemma}

\noindent {\bf Proof.} Let $n$ be the number of faces in ${\rm
st}(a_1) \cup {\rm st}(a_2)\cup {\rm st}(a_3)$. If $a_1a_2a_3$ is
a face in $M$ then clearly $n = 13$. If $a_1a_2$, $a_2a_3$,
$a_1a_3$ are edges in $M$ but $a_1a_2a_3$ is not a face then $n =
12$. In the other cases, $n \geq 14$.

Let $m$ be the number of faces in ${\rm st}(a_1) \cup \cdots \cup
{\rm st}(a_5)$. By above, $m \geq 12$ and $m = 12$ if and only if
the induced subcomplex on a set of any three vertices is a $K_3$.
So, if $m=12$ then the induced subcomplex of $M$ on $\{a_1, \dots,
a_5\}$ is a $K_5$ and hence the induced subcomplex of ${\rm
lk}(a_1)$ on $\{a_2, a_3, a_4, a_5\}$ is a null graph on 4
vertices. This is not possible since $\deg(a_1)=6$. This proves
the lemma. \hfill $\Box$

\begin{lemma} $\!\!${\bf .} \label{l3.8} Let $M$ be an $n$-vertex
connected combinatorial $2$-manifold. If the degree of each vertex
is $6$ and $n>7$ then for any vertex $u$ there exist faces of the
form $uab$, $vab$ where $uv$ is a non-edge.
\end{lemma}

\noindent {\bf Proof.} Let ${\rm lk}(u) = C_6(1, \dots, 6)$. Since
the degree of each vertex is 6, $123, \dots, 456, 561$ are not
faces. We want to show that there exists $v\not\in \{u, 1, \dots,
6\}$ such that $12v, \dots, 56v$ or $16v$ is a face. If not then
$124$ or $125$ is a face. Assume, without loss of generality, that
$124$ is a face. Then (since $146\in M$ $\Rightarrow$ $\deg(1)
=4$) the second face containing $16$ is $136$. Inductively, $256$,
$145$, $346$ and $235$ are faces. These imply that $xy$ is an edge
for $x\neq y\in U := \{u, 1, \dots, 6\}$. Since the degree of each
vertex is 6, for $x\in U$ and $z \not\in U$, $xz$ is a non-edge
and hence (since $M$ has more than 7 vertices) $M$ is not
connected. This completes the proof. \hfill $\Box$

\begin{lemma} $\!\!${\bf .} \label{l3.9} Let $M$ be a connected
combinatorial $2$-manifold and let $U$ be a set of $m$ vertices of
$M$. If the degree of each vertex in $M$ is $6$ and $m < f_0(M)$
then the number of edges in $M[U]$ is at most $3m - 3$.
\end{lemma}

\noindent {\bf Proof.} Let $V$ be the vertex-set of $M$. Let $n$
be the number of edges in $M[U]$ and let $k$ be the number of
edges of the form $ab$, where $a\in U$, $b\in V\setminus U$ (i.e.,
$k$ is the number of connecting edge between $U$ and $V\setminus
U$). Since $M$ is connected and $U\neq V$, $k\neq 0$.

Now, for any connecting edge $ab$, there exists two faces (with
vertices both in $U$ and $V\setminus U$) containing $ab$. On the
other hand, for each such face there exists exactly two connecting
edges. This implies $k\geq 3$. If $k\leq 5$ then clearly all faces
containing the connecting edges have to be of the form $ab_1b_2,
\dots, ab_{k-1}b_k, ab_kb_1$, where $a\in U$ and $b_1, \dots,
b_k\in V\setminus U$ or $a\in V\setminus U$ and $b_1, \dots,
b_k\in U$. Then $C_5(b_1, \dots, b_k)\subseteq {\rm lk}(a)$. This
is not possible. So, $k\geq 6$.

Counting two ways the number of pairs of the form $(u, e)$, where
$u\in U$ and $e$ is an edge containing $u$, we get $6m = n\times 2
+ k$ or $2n= 6m -k\leq 6m -6$. This proves the lemma. \hfill
$\Box$

\bigskip

\begin{lemma} $\!\!${\bf .} \label{tkb12v} If $M$ is a $12$-vertex
degree-regular combinatorial $2$-manifolds of Euler
characteristic $0$ then $M$ is isomorphic to $T_{12, 1, 2}$,
$T_{12, 1, 3}$, $T_{12, 1, 4}$, $T_{6, 2, 2}$, $B_{3, 4}$, $B_{4,
3}$ or $K_{3, 4}$.
\end{lemma}

\noindent {\bf Proof.} Let $M$ be a 12-vertex degree regular
combinatorial 2-manifold of Euler characteristic 0. Let the
vertex set $V$ of $M$ be $\{0, \dots, 9, u, v\}$. Let $\varphi
\colon V \to \{1, \dots, 12\}$ be given by $\varphi(i)=i$ for
$1\leq i\leq 9$, $\varphi(0) = 10$, $\varphi(u)=11$ and
$\varphi(v)=12$.

Since $\chi(M) = 0$, the degree of each vertex is 6. Assume,
without loss of generality, that ${\rm lk}(0) = C_6(1, \dots, 6)$.
Since the degree of each vertex is 6, $123, \dots, 456, 561\not\in
M$. Since each component contains at least 7 vertices, $M$ is
connected. So, by Lemma \ref{l3.8}, we may assume that $127$ is a
face. Then ${\rm lk}(1)$ has the form $C_6(7, 2, 0, 6, x, y)$, for
some $x, y \in V$. It is easy to see that $(x, y)= (3, 4)$, $(3,
5)$, $(3, 8)$, $(4, 3)$, $(4, 5)$, $(4, 8)$, $(8, 4)$, $(8, 9)$,
$(8, 3)$, $(8, 5)$. The cases  $(x, y) = (8, 3)$ and $(8, 5)$ are
isomorphic to the case $(x, y) = (4, 8)$ by the map $(0, 1)(2,
6)(3, 4, 8)(5, 7)$ and $(0, 1)(3, 7)(4, 8, 5)$ respectively. So,
we need not consider the last two cases.

\smallskip

\noindent {\bf Claim.} $(x, y)= (3, 4)$, $(3, 8)$, $(4, 3)$, $(4,
8)$, $(8, 4)$ or $(8, 9)$.

If $(x, y) = (3, 5)$, then $045$, $056$, $135$, $157$ are faces
and hence ${\rm lk}(5) = C_6(4, 0, 6, 3, 1, 7)$. This implies that
$C_3(1, 5, 6) \subseteq {\rm lk}(3)$. This is not possible.

If $(x, y) = (4, 5)$ then ${\rm lk}(4) = C_6(6, 1, 5, 0, 3, z)$,
where $z = 7, 8, 9, u$ or $v$. If $z = 7$ then, as in the previous
case, we get a contradiction. So, we may assume that ${\rm lk}(4)
= C_6(6, 1, 5, 0, 3, 8)$ and hence ${\rm lk}(6) = C_6(8, 4, 1, 0,
5, w)$, ${\rm lk}(5) = C_6(7, 1, 4, 0, 6, w)$ for some $w \in V$.
It is easy to see that $w = 9, u$ or $v$. In any case, we get 15
faces not containing any from $\{9, u, v\}\setminus \{w\}$. This
is not possible since $M$ has 24 faces. This proves the claim.

\smallskip

\noindent {\bf Case 1.} $(x, y) = (3, 4)$, i.e., ${\rm lk}(1) =
C_6(7, 2, 0, 6, 3, 4)$. Now, ${\rm lk}(3) = C_6(6, 1, 4, 0, 2,
z)$ for some $z\in V$. If $z = 5$ then $C_4(5, 3, 1, 0) \subseteq
{\rm lk}(6)$. If $z = 7$ then $C_4(7, 3, 0, 1) \subseteq {\rm
lk}(2)$. This implies that $z\in \{8, 9, u, v\}$. Assume, without
loss of generality, that ${\rm lk}(3) = C_6(6, 1, 4, 0, 2, 8)$.
Now, ${\rm lk}(2) = C_6(8, 3, 0, 1, 7, w)$, for some $w \in V$.
If $w \not\in\{9, u, v\}$ then we get 14 faces not containing any
of $9$, $u$, $v$. This is not possible by Lemma \ref{l3.7}. So,
assume without loss of generality, that $z=9$, i.e., ${\rm lk}(2)
= C_6(8, 3, 0, 1, 7, 9)$.

Completing successively, we get ${\rm lk}(6) = C_6(8, 3, 1, 0, 5,
u)$, ${\rm lk}(8) = C_6(u, 6, 3, 2, 9, v)$, ${\rm lk}(4) = C_6(7,
1, 3, 0, 5, v)$, ${\rm lk}(5) = C_6(v, 4, 0, 6, u, 9)$, ${\rm
lk}(9) = C_6(7, u, 5, v, 8, 2)$, ${\rm lk}(u) = C_6(7, 9, 5, 6, 8,
v)$. Here $M \cong B_{3, 4}$ by the map $\varphi_{34}\circ (0, 9,
3, 4, 7, u, 2, v)(1, 8)(5, 6)$, where $\varphi_{34}\colon V\to
V(B_{3, 4})$ is given by $\varphi_{34}(i)=v_{1i}$, $\varphi_{34}(3
+i) = v_{2i}$, $\varphi_{34}(6+i)=v_{3i}$, for $1\leq i\leq 3$,
$\varphi_{34}(0) = v_{41}$, $\varphi_{34}(u) = v_{42}$,
$\varphi_{34}(v) = v_{43}$.

\smallskip

\noindent {\bf Case 2.} $(x, y) = (3, 8)$, i.e., ${\rm lk}(1) =
C_6(7, 2, 0, 6, 3, 8)$. Now, $023$, $034$, $136$ and $138$ are
faces in $M$. So, ${\rm lk}(3) = C_6(2, 0, 4, 8, 1, 6)$ or $C_6(2,
0, 4, 6, 1, 8)$. In the first case, ${\rm lk}(2) = C_6(6, 3, 0, 1,
7, z)$ for some $z \in V$. As in Case 1, $z$ is a new vertex, say,
$9$. Then ${\rm lk}(2) = C_6(6, 3, 0, 1, 7, 9)$ and ${\rm lk}(6) =
C_6(9, 2, 3, 1, 0, 5)$. This gives 15 faces not containing $u$ or
$v$. This is not possible. Thus, ${\rm lk}(3) = C_6(2, 0, 4, 6, 1,
8)$. Now, ${\rm lk}(6) = C_6(4, 3, 1, 0, 5, w)$ for some $w \in
V$. If $w =2, 7$ or $8$, then we get 14 faces not containing any
of $9$, $u$, $v$. This is not possible by Lemma \ref{l3.7}. So,
assume without loss, that ${\rm lk}(6) = C_6(4, 3, 1, 0, 5, 9)$.

Completing successively, we get ${\rm lk}(4) = C_6(9, 6, 3, 0, 5,
u)$, ${\rm lk}(5) = C_6(4, 0, 6, 9, v, u)$, ${\rm lk}(9) = C_6(4,
6, 5, v, 7, u)$, ${\rm lk}(7) = C_6(9, u, 8, 1, 2, v)$, ${\rm
lk}(2) = C_6(7, 1, 0, 3, 8, v)$, ${\rm lk}(8) = C_6(2, 3, 1, 7, u,
v)$. Here $M \cong T_{12, 1, 2}$ by the map $\varphi\circ (0, 4,
6, 5, 7, u, 9, 8, v)(1, 2)$.

\smallskip

\noindent {\bf Case 3.} $(x, y) = (4, 3)$, i.e., ${\rm lk}(1) =
C_6(7, 2, 0, 6, 4, 3)$. Now, ${\rm lk}(4) = C_6(6, 1, 3, 0, 5,
z)$ for some $z \in V$. If $z = 2$, then ${\rm lk}(2)$ has 7
vertices. If $z = 7$, then ${\rm lk}(4) = C_6(6, 1, 3, 0, 5, 7)$
and hence $127$, $137$, $457$, $467$ are faces in $M$. This
implies that ${\rm lk}(7)= C_6(2, 1, 3, 5, 4, 6)$ or $C_6(2, 1,
3, 6, 4, 5)$. In either cases we get 14 faces not containing any
of $8$, $9$, $u$ or $v$. This is not possible by Lemma
\ref{l3.7}. This implies that $z = 8, 9, u$ or $v$. Assume,
without loss, that $z = 8$. This case is now isomorphic to Case 2
by the map $(0, 3, 6, 2, 4, 1)(5, 8, 7)(9, v)$.

\smallskip

\noindent {\bf Case 4.} $(x, y) = (4, 8)$, i.e., ${\rm lk}(1) =
C_6(6, 0, 2, 7, 8, 4)$. This gives ${\rm lk}(4) = C_6(8, 1, 6, 3,
0, 5)$.

Now, ${\rm lk}(6) = C_6(3, 4, 1, 0, 5, z)$ for some $z \in V$. By
using Lemma \ref{l3.7}, $z=9$, $u$ or $v$. So, assume that ${\rm
lk}(6) = C_6(3, 4, 1, 0, 5, 9)$. This implies that ${\rm lk}(5) =
C_6(9, 6, 0, 4, 8, u)$. This case is now isomorphic to Case 1 by
the map $(1, 6)(2, 5)(3, 4)(9, 7, u)$.

\smallskip

\noindent {\bf Case 5.} $(x, y) = (8, 4)$, i.e., ${\rm lk}(1) =
C_6(6, 0, 2, 7, 4, 8)$. Now, $034$, $045$, $147$ and $148$ are
faces in $M$. So, ${\rm lk}(4) = C_6(3, 0, 5, 7, 1, 8)$ or
$C_6(3, 0, 5, 8, 1, 7)$.

\smallskip

\noindent {\bf Subcase 5.1.} ${\rm lk}(4) = C_6(3, 0, 5, 7, 1,
8)$. Then (by using Lemma \ref{l3.7}) ${\rm lk}(7) = C_6(5, 4, 1,
2, z$, $w)$, where $z, w \in \{9, u, v\}$. So, assume without
loss, that ${\rm lk}(7) = C_6(5, 4, 1, 2, 9, u)$. Completing
successively, we get ${\rm lk}(2) = C_6(9, 7, 1, 0, 3, v)$, ${\rm
lk}(3) = C_6(2, 0, 4, 8, u, v)$, ${\rm lk}(u) = C_6(8, 9, 7, 5, v,
3)$, ${\rm lk}(8) = C_6(9, u, 3, 4, 1, 6)$, ${\rm lk}(5) = C_6(7,
4, 0, 6, v, u)$, ${\rm lk}(6) = C_6(5, 0, 1, 8, 9$, $v)$. Here $M
\cong K_{3, 4}$ by the map $\psi_{34}\circ (0, 8)(1, 9)(2, v, 3,
u)(4, 7)$, where $\psi_{34}$ is same as $\varphi_{34}$ (of Case 1)
on the vertex-set.

\smallskip

\noindent {\bf Subcase 5.2.} ${\rm lk}(4) = C_6(3, 0, 5, 8, 1,
7)$. Then ${\rm lk}(7) = C_6(3, 4, 1, 2, z, w)$, for some $z, w
\in V$. As in the previous case, $z, w\in \{9, u, v\}$. So, assume
without loss, that ${\rm lk}(7) = C_6(3, 4, 1, 2, 9, u)$. Then
${\rm lk}(2) = C_6(9, 7, 1, 0, 3, a)$ for some $a \in V$. It is
easy to see that $a=8$ or $v$. If $a =8$ then, considering ${\rm
lk}(8)$, we get 19 faces not containing $v$. This is not possible
since $f_2(M)=24$. So, ${\rm lk}(2) = C_6(9, 7, 1, 0, 3, v)$ and
hence ${\rm lk}(3) = C_6(2, 0, 4, 7, u, v)$. Then, ${\rm lk}(8) =
C_6(5, 4, 1, 6, b, c)$, where $b, c \in \{9, u, v\}$. Since the
set of known faces is invariant under $(1,4)(2,3)(5, 6)(9, u)$,
we may assume that $(b, c) = (9, u)$, $(u, 9)$, $(v, u)$ or $(v,
9)$.

\smallskip

\noindent {\bf Subcase 5.2.1.} ${\rm lk}(8) = C_6(5, 4, 1, 6, 9,
u)$. Completing successively, we get ${\rm lk}(9) = C_6(2, 7, u,
8, 6, v)$, ${\rm lk}(6) = C_6(9, 8, 1, 0, 5, v)$, ${\rm lk}(5) =
C_6(6, 0, 4, 8, u, v)$. Here $M \cong T_{12, 1, 4}$ by the map
$\varphi\circ (0, 2, 4, u, 9, 6)(1, v, 8)(3, 5, 7)$.

\smallskip

\noindent {\bf Subcase 5.2.2.} ${\rm lk}(8) = C_6(5, 4, 1, 6, u,
9)$. Completing successively, we get ${\rm lk}(u) = C_6(6, 8, 9,
7, 3, v)$, ${\rm lk}(9) = C_6(5, 8, u, 7, 2, v)$, ${\rm lk}(v) =
C_6(6, u, 3, 2, 9, 5)$. Here $M \cong B_{4, 3}$ by the map
$\varphi_{43} \circ (0, 9, u, 7, 6, 8, v, 3, 2)(1, 5, 4)$, where
$\varphi_{43} \colon V\to V(B_{4, 3})$ is given by
$\varphi_{43}(i)=v_{1i}$, $\varphi_{43}(4+i)= v_{2i}$, for $1\leq
i\leq 4$, $\varphi_{43}(9)=v_{31}$, $\varphi_{43}(0) = v_{32}$,
$\varphi_{43}(u) = v_{33}$, $\varphi_{43}(v) = v_{34}$.

\smallskip

\noindent {\bf Subcase 5.2.3.} ${\rm lk}(8) = C_6(5, 4, 1, 6, v,
9)$. Completing successively, we get ${\rm lk}(v) = $ \linebreak
$C_6(3, 2, 9, 8, 6, u)$, ${\rm lk}(9) = C_6(8, v, 2, 7, u, 5)$,
${\rm lk}(u) = C_6(6, 5, 9, 7, 3, v)$. Here $M \cong T_{12, 1,
3}$ by the map $\varphi\circ (1, 0, 8, 6)(2, 5, u, 4, 3, v, 9)$.

\smallskip

\noindent {\bf Subcase 5.2.4.} ${\rm lk}(8) = C_6(5, 4, 1, 6, v,
u)$. Completing successively, we get ${\rm lk}(v) = C_6(3, 2, 9,
6, 8, u)$, ${\rm lk}(u) = C_6(3, 7, 9, 5, 8, v)$, ${\rm lk}(9) =
C_6(6, 5, u, 7, 2, v)$. Here $M \cong B_{4, 3}$ by the map
$\varphi_{43} \circ (1, 2, 9, 4, 6, 3, 5, u, v, 8, 7)$, where
$\varphi_{43}$ is as in Subcase 5.2.2.

\smallskip

\noindent {\bf Case 6.} $(x, y)=(8,9)$, i.e., ${\rm lk}(1) =
C_6(6, 0, 2, 7, 9, 8)$. Now, ${\rm lk}(2) = C_6(7, 1, 0, 3, z,
w)$ for some $z, w \in V$. It is easy to see that $(z, w)= (5,
6)$, $(6, 8)$, $(5, 8)$, $(5, u)$, $(u, 8)$, $(6, 5)$, $(8, 6)$,
$(8, u)$, $(u, 5)$, $(9, u)$, $(u, 4)$, $(u, v)$. If $(z, w)= (5,
8)$, i.e., ${\rm lk}(2) = C_6(7, 1, 0, 3, 5, 8)$ then $045$,
$056$, $235$, $258$ are faces in $M$. This implies that ${\rm
lk}(5)= C_6(6, 0, 4, 8, 2, 3)$. Then $\deg(8)\geq 7$. Since the
set of known faces are invariant under $(0, 1)(3, 7)(4, 9)(5,
8)$, we may assume that $(z, w)= (5, 6)$, $(5, u)$, $(6, 5)$,
$(u, 5)$, $(9, u)$ or $(u, v)$.

\smallskip

\noindent {\bf Subcase 6.1.} ${\rm lk}(2) = C_6(7, 1, 0, 3, 5,
6)$. Now, it is easy to see that ${\rm lk}(6) = C_6(7, 2, 5, 0, 1,
8)$. Now, completing successively, we get ${\rm lk}(7) = C_6(8,
6, 2, 1, 9, u)$, ${\rm lk}(8) = C_6(u, 7, 6, 1, 9, v)$, ${\rm
lk}(9) = C_6(8, 1, 7, u, 4, v)$, ${\rm lk}(4) = C_6(3, 0, 5, v,
9, u)$, ${\rm lk}(5) = C_6(4, 0, 6, 2, 3, v)$, ${\rm lk}(3) =
C_6(5, 2, 0, 4, u, v)$. Here $M \cong T_{12, 1, 2}$ by the map
$\varphi\circ (0, 9, 3, u, 2, 8, 4, v, 1, 6, 7, 5)$.

\smallskip

\noindent {\bf Subcase 6.2.} ${\rm lk}(2) = C_6(7, 1, 0, 3, 5,
u)$. Now, it is easy to see that ${\rm lk}(5) = C_6(4, 0, 6, 3, 2,
u)$. Now, completing successively, we get ${\rm lk}(3) = C_6(2,
0, 4, v, 6, 5)$, ${\rm lk}(4) = C_6(3, 0, 5, u, 9, v)$, ${\rm
lk}(9) = C_6(8, 1, 7, v, 4, u)$, ${\rm lk}(7) = C_6(9, 1, 2, u,
8, v)$, ${\rm lk}(8) = C_6(1, 6, v, 7, u, 9)$, ${\rm lk}(u) =
C_6(2, 5, 4, 9, 8, 7)$. Here $M \cong B_{3, 4}$ by the map
$\varphi_{34} \circ (0, 9, 3, 4, v, 1, 6, 5, 8, 2, 7)$, where
$\varphi_{34}$ is as in Case 1.

\smallskip

\noindent {\bf Subcase 6.3.} ${\rm lk}(2) = C_6(7, 1, 0, 3, 6,
5)$. Completing successively, we get ${\rm lk}(6) = C_6(5, 0, 1$,
$8, 3, 2)$, ${\rm lk}(3) = C_6(2, 0, 4, u, 8, 6)$, ${\rm lk}(5) =
C_6(6, 0, 4, v, 7, 2)$, ${\rm lk}(8) = C_6(3, 6, 1, 9, v, u)$,
${\rm lk}(7) = C_6(2, 1, 9, u, v, 5)$, ${\rm lk}(4) = C_6(5, 0,
3, u, 9, v)$, ${\rm lk}(u) = C_6(4, 3, 8, v, 7, 9)$. Here $M\cong
B_{3, 4}$ by the map $\varphi_{34} \circ (0, 9, 3, 5, 7, u, 2, 8,
1, v)(4, 6)$, where $\varphi_{34}$ is as in Case 1.

\smallskip

\noindent {\bf Subcase 6.4.} ${\rm lk}(2) = C_6(7, 1, 0, 3, 9,
u)$. Completing successively, we get ${\rm lk}(9) = C_6(8, 1, 7$,
$3, 2, u)$, ${\rm lk}(3) = C_6(2, 0, 4, v, 7, 9)$, ${\rm lk}(7) =
C_6(3, 9, 1, 2, u, v)$, ${\rm lk}(u) = C_6(7, 2, 9, 8, 5, v)$,
${\rm lk}(8) = C_6(5, u, 9, 1, 6, 4)$, ${\rm lk}(5) = C_6(8, 4,
0, 6, v, u)$, ${\rm lk}(6) = C_6(5, 0, 1, 8, 4, v)$. Here $M\cong
B_{3,4}$ by the map $\varphi_{34}\circ (0, 7, 2)(1, u)(3, 6, 8,
v, 5, 4, 9)$, where $\varphi_{34}$ is as in Case 1.

\smallskip

\noindent {\bf Subcase 6.5.} ${\rm lk}(2) = C_6(7, 1, 0, 3, u,
5)$. Now, it is easy to see that ${\rm lk}(5)= C_6(4, 0, 6, u, 2,
7)$ or $C_6(4, 0, 6, 7, 2, u)$. The first case is isomorphic to
Subcase 5.1 by the map $(0, 1, 7, 5)(2, 4, 6)$ $(3, 8, u)$. The
second case is isomorphic to Subcase 5.2 by the map $(0, 4, 5)(1,
7, 2)(3, 8, u, 6)$.

\smallskip

\noindent {\bf Subcase 6.6.} ${\rm lk}(2) = C_6(7, 1, 0, 3, u,
v)$. Now, it is easy to see that ${\rm lk}(7) = C_7(v, 2, 1, 9, a,
b)$, where $(a, b)= (3,4)$, $(4,3)$, $(4,5)$, $(4,8)$, $(5,4)$,
$(5,6)$, $(5,8)$, $(6,5)$, $(6,8)$, $(u,3)$, $(u,4)$, $(u,5)$,
$(u,8)$. Since the set of known faces is invariant under the map
$(1, 2)(3,6)(4,5)(8,u)(9,v)$, we may assume that $(a, b) =(3,4)$,
$(4,3)$, $(4,5)$, $(4,8)$, $(5,4)$, $(5,8)$, $(6,8)$, $(u,8)$.

\smallskip

\noindent {\bf Claim.} $(a, b) = (3,4)$ or $(5,4)$.

If $(a, b) = (4, 3)$ then, considering ${\rm lk}(3)$, we get
$C_4(u, 3, 7, 2) \subseteq {\rm lk}(v)$. If $(a,b) = (4,5)$ then,
${\rm lk}(5)$ can not be a 6-cycle. If $(a, b) = (4, 8)$ then,
considering ${\rm lk}(9)$, we see that 0, 1, 4, 7, 8, $9 \not \in
{\rm lk}(u)$. This is not possible. If $(a, b) = (5, 8)$ then,
considering ${\rm lk}(8)$, we get $C_4(5, 8, 1, 0) \subseteq {\rm
lk}(6)$. If $(a,b) = (6,8)$ then, considering the links of $6, 9$
and $8$ successively, we get $C_4(u, 8, 7, 2) \subseteq {\rm
lk}(v)$. If $(a, b) = (u,8)$ then, considering ${\rm lk}(u)$, we
get $7$ vertices in ${\rm lk}(8)$. These prove the claim.

\smallskip

\noindent {\bf Subcase 6.6.1.} ${\rm lk}(7) = C_6(v, 2, 1, 9, 3,
4)$. Completing successively, we get ${\rm lk}(3) = C_6(0, 2$, $u,
9, 7, 4)$, ${\rm lk}(9) = C_6(1, 7, 3, u, 5, 8)$, ${\rm lk}(4) =
C_6(0, 3, 7, v, 8, 5)$, ${\rm lk}(5) = C_6(0, 4, 8, 9, u, 6)$,
${\rm lk}(u) = C_6(2, 3, 9, 5, 6, v)$, ${\rm lk}(6) = C_6(0, 1,
8, v, u, 5)$. Now, $M\cong K_{3, 4}$ by the map $\psi_{34} \circ
(0, 4, 5, 7, 3, 2, 1, v, 6, 9, u)$, where $\psi_{34}$ is as in
Subcase 5.1.

\smallskip

\noindent {\bf Subcase 6.6.2.} ${\rm lk}(7) = C_6(v, 2, 1, 9, 5,
4)$. Now, it is easy to see that ${\rm lk}(4) = C_6(3, 0, 5, 7$,
$v, 8)$. Now, $0, 1, 4, 7, 8 \not \in {\rm lk}(u)$. So, $5 \in
{\rm lk}(u)$ and hence ${\rm lk}(5) = C_6(6, 0, 4, 7, 9, u)$. Then
${\rm lk}(3) = C_6(8, 4, 0, 2, u, 6)$ or $C_6(8, 4, 0, 2, u, 9)$.

\smallskip

\noindent {\bf Subcase 6.6.2.1.} ${\rm lk}(3) = C_6(8, 4, 0, 2, u,
6)$. Completing successively, we get ${\rm lk}(6) = C_6(0, 1, 8,
3, u, 5)$, ${\rm lk}(8) = C_6(1, 6, 3, 4, v, 9)$, ${\rm lk}(9) =
C_6(1, 7, 5, u, v, 8)$. Now, $M\cong K_{3, 4}$ by the map
$\psi_{34} \circ (0, 8, 4, v, 1, 5, u)(2, 6, 7, 3, 9)$, where
$\psi_{34}$ is as in Subcase 5.1.

\smallskip

\noindent {\bf Subcase 6.6.2.2.} ${\rm lk}(3) = C_6(8, 4, 0, 2, u,
9)$. Completing successively, we get ${\rm lk}(9) = C_6(1, 7, 5,
u, 3, 8)$, ${\rm lk}(8) = C_6(1, 6, v, 4, 3, 9)$, ${\rm lk}(6) =
C_6(0, 1, 8, v, u, 5)$. Now, $M \cong T_{6, 2, 2}$ by the map
$\psi \circ (0, 1, u, 7, 4, 3, 8)(2, 6, v, 5)$, where $\psi\colon
V\to V(T_{6, 2, 2})$ is given by $\psi(i)= u_i$, for $1\leq i\leq
6$, $\psi(5+i) = v_{i}$, for $2\leq i\leq 4$, $\psi(0)=v_5$,
$\psi(u) = v_6$ and $\psi(v)=v_7$. \hfill $\Box$

\begin{lemma} $\!\!${\bf .}  \label{tkb14v} If $M$ is a $14$-vertex
degree-regular combinatorial $2$-manifolds of Euler
characteristic $0$ then $M$ is isomorphic to  $T_{14, 1, 2}$,
$T_{14, 1, 3}$ or $Q_{7, 2}$.
\end{lemma}

\noindent {\bf Proof.} Let $M$ be a 14-vertex degree regular
combinatorial 2-manifold of Euler characteristic 0. Let the
vertex set $V$ be $\{0, 1, \ldots, 9, u, v, w, z\}$. Let $\varphi
\colon V \to \{1, \dots, 14\}$ be given by $\varphi(i) = i$, for
$1\leq i\leq 9$, $\varphi(0) = 10$, $\varphi(u) = 11$,
$\varphi(v) = 12$, $\varphi(w) = 13$ and $\varphi(z) = 14$.

Since $\chi(M) = 0$, the degree of each vertex is 6. Assume
without loss that ${\rm lk}(0) = C_6(1, 2, 3, 4, 5, 6)$. By Lemma
\ref{l3.8}, ${\rm lk}(1)= C_6(6, 0, 2, 7, x, y)$, for some $x, y
\in V$. It is easy to see that $(x, y) = (3, 4)$, $(3, 8)$, $(4,
3)$, $(4, 8)$, $(5, 3)$, $(5, 4)$, $(5, 8)$, $(8, 3)$, $(8, 4)$,
$(8, 9)$. The case $(x, y) = (3, 8)$ is isomorphic to the case
$(x, y) = (5,8)$ by the map $(2, 6)(3, 5)(7, 8)$ and to the case
$(x,y) = (8, 4)$ by the map $(0, 1)(2, 6)(3, 4, 8)(5,7)$. Hence we
may assume that $(x, y) = (3, 4)$, $(3, 8)$, $(4, 3)$, $(4, 8)$,
$(5, 3)$, $(5, 4)$, $(8, 3)$, $(8, 9)$.

\smallskip

\noindent {\bf Claim.} $(x, y) = (3,4)$, $(8,3)$ or $(8,9)$

If $(x, y) = (4, 3)$ then, considering the links of $3, 6, 2, 4,
7, 8$, $u$ successively we get $C_3(v, u, 7) \subseteq {\rm
lk}(w)$. So, $(x, y) \neq (4, 3)$. If $(x, y) = (5, 3)$ then,
${\rm lk}(5) = C_6(1, 3, 6, 0, 4, 7)$. But then $C_4(0, 1, 3, 5)$
$ \subseteq {\rm lk}(6)$. So, $(x, y) \neq (5, 3)$. Similarly,
$(x, y) \neq (3, 8)$, $(4, 8)$ or $(5, 4)$. This proves the claim.

\smallskip

\noindent {\bf Case 1.} $(x, y) = (3, 4)$, i.e., ${\rm lk}(1)=
C_6(6, 0, 2, 7, 3, 4)$. Now, it easy to see that ${\rm lk}(3) =
C_6(7, 1, 4, 0, 2, 5)$ or $C_6(7, 1, 4, 0, 2, 8)$. In the first
case, ${\rm lk}(5) = C_6(4, 0, 6, 2, 3, 7)$ or $C_6(4, 0, 6$, $7,
3, 2)$. In both these cases we have $23$ edges in $M[\{0, \dots
7\}]$. This is not possible by Lemma \ref{l3.9}. Thus, ${\rm
lk}(3) = C_6(7, 1, 4, 0, 2, 8)$. Now, ${\rm lk}(2) = C_6(7, 1, 0,
3, 8, b)$ for some $b\in V$. It is easy to see that $b = 5$ or 9.
If $b = 5$ then, considering the links of $4$ and $5$, we get
$\geq 17$ faces not containing any of $9$, $u$, $v$, $w$, $z$.
This is not possible by Lemma \ref{l3.7}. Thus ${\rm lk}(2) =
C_6(7, 1, 0, 3, 8, 9)$. Again, by using Lemma \ref{l3.7}, we
successively get ${\rm lk}(7) = C_6(8, 3, 1, 2, 9, u)$, ${\rm
lk}(8) = C_6(9, 2, 3, 7$, $u, v)$, ${\rm lk}(9) = C_6(u, 7, 2, 8,
v, w)$, ${\rm lk}(u) = C_6(v, 8, 7, 9, w, z)$, ${\rm lk}(v) =
C_6(w, 9, 8$, $u, z, 5)$. Then $045$, $056$, $5vw$ and $5vz$ are
faces. So, ${\rm lk}(5) = C_6(6, 0, 4, w, v, z)$ or $C_6(6, 0, 4,
z, v, w)$. In the first case, considering the links of $4$ and
$6$, we get $C_5(w, 6, 5, v, u) \subseteq {\rm lk}(z)$. Thus
${\rm lk}(5) = C_6(6, 0, 4, z, v, w)$.

Now, completing successively, we get ${\rm lk}(w) = C_6(9, v, 5,
6$, $z, u)$ and ${\rm lk}(6) = C_6(1, 0, 5, w$, $z, 4)$. Here $M$
is isomorphic to $T_{14, 1, 2}$ by the map $\varphi \circ (0, 7,
3, 5)(1, 6, 9)(2, 4, 8)(u, z)(v, w)$.

\smallskip

\noindent {\bf Case 2.} $(x, y) = (8, 3)$. Since $023$, $034$,
$136$ and $138$ are faces, ${\rm lk}(3) = C_6(2, 0, 4, 8, 1, 6)$
or $C_6(2, 0, 4, 6, 1, 8)$. In the first case, considering the
links of 6, 2, 9, 7, 8 and 4 successively, we get $7$ vertices in
${\rm lk}(u)$. Thus, ${\rm lk}(3) = C_6(2, 0, 4, 6, 1, 8)$.

Now, completing successively, we get ${\rm lk}(6) = C_6(5, 0, 1,
3, 4, 9)$, ${\rm lk}(4) = C_6(5, 0, 3, 6, 9, u)$, ${\rm lk}(5) =
C_6(u, 4, 0, 6, 9, v)$, ${\rm lk}(9) = C_6(u, 4, 6, 5, v, w)$,
${\rm lk}(2) = C_6(8, 3, 0, 1, 7, z)$, ${\rm lk}(8) = C_6(7, 1, 3,
2, z, w)$, ${\rm lk}(u) = C_6(4, 5, v, z, w, 9)$, ${\rm lk}(z) =
C_6(2, 7, v, u, w, 8)$ and ${\rm lk}(7) = C_6(1, 2$, $z, v, w,
8)$. Here $M \cong T_{14, 1, 2}$ by the map $\varphi \circ (0, 6,
7, 1, 4, 8, 2, 3, 5, 9)$.

\smallskip

\noindent {\bf Case 3.} $(x, y) = (8, 9)$, i.e., ${\rm lk}(1) =
C_6(6, 0, 2, 7, 8, 9)$. Now, ${\rm lk}(6) = C_6(9, 1, 0, 5, a,
b)$, for some $a, b \in V$. It is easy to see that $(a, b) = (2,
3)$, $(2, 7)$, $(3, 2)$, $(3, 4)$, $(3, 7)$, $(3, u)$, $(7, 2)$,
$(7, 3)$, $(7, 4)$, $(7, u)$, $(8, 3)$, $(8, 4)$, $(8, 7)$, $(8,
u)$, $(u, 3)$, $(u, 4)$, $(u, 7)$, $(u, v)$. Since the set of
known faces is invariant under the map $(0, 1)(3, 7)(4, 8)(5,
9)$, we may assume that $(a, b) =$ $(2, 3)$, $(2, 7)$, $(3, 4)$,
$(3, 7)$, $(3, u)$, $(7, 3)$, $(7, 4)$, $(7, u)$, $(8, 4)$, $(8,
u)$ or $(u, v)$.

By the similar arguments as in the previous claim one gets $(a,
b) = (2,7)$ or $(u,v)$.

\smallskip

\noindent {\bf Subcase 3.1.} ${\rm lk}(6) = C_6(9, 1, 0, 5, 2,
7)$. Completing successively, we get ${\rm lk}(2) = C_6(1, 0, 3$,
$5, 6, 7)$, ${\rm lk}(5) = C_6(4, 0, 6, 2, 3, u)$, ${\rm lk}(3) =
C_6(2, 0, 4, v, u, 5)$, ${\rm lk}(4) = C_6(3, 0, 5, u, w, v)$,
${\rm lk}(u) = C_6(3, 5, 4, w, z, v)$, ${\rm lk}(7) = C_6(1, 2,
6, 9, z, 8)$, ${\rm lk}(9) = C_6(1, 6, 7, z, w, 8)$, ${\rm lk}(z)
= C_6(7, 8, v, u, w, 9)$ and ${\rm lk}(v) = C_6(3, 4, w, 8, z,
u)$. Here $M$ is isomorphic to $T_{14, 1, 2}$ by the map $\varphi
\circ (0, 4, 1, 7, 8)(2, 5$, $3)(u, z)(v, w)$.

\smallskip

\noindent {\bf Subcase 3.2.} ${\rm lk}(6) = C_6(9, 1, 0, 5, u,
v)$. Now, ${\rm lk}(5) = C_6(u, 6, 0, 4, c, d)$, for some $c, d
\in V$. It is easy to see that $(c, d) = (2, 3)$, $(2, 7)$, $(7,
2)$, $(7, 3)$, $(7, 8)$, $(7, w)$, $(8, 3)$, $(8, 7)$, $(8, 9)$,
$(8, w)$, $(9, 8)$, $(v, 3)$, $(v, 7)$, $(v, 8)$, $(v, 9)$, $(v,
w)$, $(w, 3)$, $(w, 7)$, $(w, 8)$, $(w, z)$. Since the set of
known faces is invariant under the map $(0, 6)(2, 9)(3, v)(4,
u)(7, 8)$, we may assume that $(c, d) = (2, 3)$, $(2, 7)$, $(7,
2)$, $(7, 3)$, $(7, 8)$, $(7, w)$, $(8, 3)$, $(8, 7)$, $(8, w)$,
$(v, 3)$, $(v, w)$ or $(w, z)$.

\smallskip

\noindent {\bf Claim.} $(c, d) = (7, 8)$, $(8, w)$ or $(w, z)$.

If $(c, d) = (2, 3)$ then, ${\rm lk}(5) = C_6(u, 6, 0, 4, 2, 3)$.
Considering the links of $2$, $4$, $3$, $u$, $w$ successively, we
get $C_5(8, w, 4, 2, 1) \subseteq {\rm lk}(7)$. If $(c, d) = (2,
7)$ then, considering ${\rm lk}(2)$ we get $C_4(3, 2, 5, 0)
\subseteq {\rm lk}(4)$. If $(c, d) = (7, 3)$ then, considering
${\rm lk}(3)$, we get $C_4(7, 3, 0, 1) \subseteq {\rm lk}(2)$. If
$(c, d) = (v, 3)$ then, considering ${\rm lk}(3)$, we get $7$
vertices in ${\rm lk}(v)$. So, $(c, d) \neq (2, 3), (2, 7)$, $(7,
3)$ or $(v, 3)$. Similarly, $(c, d) \neq (7, 2)$, $(7, w)$, $(8,
3)$, $(8, 7)$ or $(v, w)$. This proves the claim.

\smallskip

\noindent {\bf Subcase 3.2.1.} ${\rm lk}(5) = C_6(u, 6, 0, 4, 7,
8)$. Now, ${\rm lk}(8) = C_6(u, 5, 7, 1$, $9, x)$, for some $x
\in V$. It is easy to check that $x = 3, w$. By using Lemma
\ref{l3.7}, we get $x\neq 3$. So, ${\rm lk}(8) = C_6(u, 5, 7, 1,
9, w)$. This implies that ${\rm lk}(9) = C_6(w, 8, 1, 6, v, y)$,
for some $y \in V$. It is easy to see that $y = 3$ or $z$.

\smallskip

\noindent {\bf Subcase 3.2.1.1.} ${\rm lk}(9) = C_6(w, 8, 1, 6,
v, 3)$. This implies that ${\rm lk}(u) = C_6(w, 8, 5, 6, v, z)$.
Again, by using Lemma \ref{l3.7}, we get ${\rm lk}(7) = C_6(4, 5,
8, 1, 2, z)$. Then ${\rm lk}(4) = C_6(3, 0, 5, 7, z, a)$, for
some $a \in V$. Considering ${\rm lk}(3)$, we get $a = v$ or $w$.

\smallskip

\noindent {\bf A.} ${\rm lk}(4) = C_6(3, 0, 5, 7, z, v)$.
Completing successively we get ${\rm lk}(z) = C_6(2, 7, 4, v, u,
w)$, ${\rm lk}(2) = C_6(1, 0, 3, w, z, 7)$ and ${\rm lk}(3) =
C_6(2, 0, 4, v, 9, w)$. Here $M \cong T_{14,1,3}$ by the map
$\varphi \circ (0, 1, 4, v, w, 6, z, 9, 3, 2, 5, u)(7, 8)$.

\smallskip

\noindent {\bf B.} ${\rm lk}(4) = C_6(3, 0, 5, 7, z, w)$.
Completing successively, we get ${\rm lk}(3) = C_6(2, 0, 4, w, 9,
v)$, ${\rm lk}(2) = C_6(1, 0, 3, v, z, 7)$ and ${\rm lk}(v) =
C_6(2, 3, 9, 6, u, z)$. Now, $M \cong Q_{7, 2}$ by the map
$\varphi \circ (0, 1, 9, u, v, 5, 2, 7, 8)(3, w, 4, z, 6)$.

\smallskip

\noindent {\bf Subcase 3.2.1.2.} ${\rm lk}(9) = C_6(w, 8, 1, 6,
v, z)$. Then it follows that ${\rm lk}(u) = C_6(v, 6, 5, 8$, $w,
3)$. This case is now isomorphic to the Subcase 3.2.1.1 by the
map $(1, 5)(2,4)(u,9)$.

\smallskip

\noindent {\bf Subcase 3.2.2.} ${\rm lk}(5) = C_6(u, 6, 0, 4, 8,
w)$. Since, $178$, $189$, $458$ and $58w$ are faces, ${\rm lk}(8)
= C_6(4, 5, w, 9, 1, 7)$ or $C_6(4, 5, w, 7, 1, 9)$.

\smallskip

\noindent {\bf Subcase 3.2.2.1.} ${\rm lk}(8) = C_6(4, 5, w, 9,
1, 7)$. Then ${\rm lk}(4) = C_6(3, 0, 5, 8, 7, x)$, for some $x
\in V$. It is easy to see that $x = v, z$. If $x = v$ then,
considering the links of $4$, $7$, $v$ successively, we obtain
$29$ faces in $M$, which is not possible. Thus ${\rm lk}(4) =
C_6(3, 0, 5, 8, 7, z)$.

Completing successively, we get ${\rm lk}(7) = C_6(2, 1, 8, 4, z,
v)$, ${\rm lk}(2) = C_6(3, 0, 1, 7, v, u)$, ${\rm lk}(v) = C_6(2,
u, 6, 9, z, 7)$, ${\rm lk}(9) = C_6(1, 6, v, z, w, 8)$, ${\rm
lk}(z) = C_6(3, 4, 7, v, 9, w)$ and ${\rm lk}(w) = C_6(3, u, 5, 8,
9, z)$. Here $M\cong T_{14, 1, 3}$ by the map $\varphi \circ (0,
6, 2, 9, 1, 5, 3)(4, 7, 8)(u, w, z)$.

\smallskip

\noindent {\bf Subcase 3.2.2.2.} ${\rm lk}(8) = C_6(4, 5, w, 7,
1, 9)$. Then ${\rm lk}(4) = C_6(3, 0, 5, 8, 9, z)$ and ${\rm
lk}(9) = C_6(1, 6, v, z, 4, 8)$. Now, ${\rm lk}(w) = C_6(u, 5, 8,
7, a, b)$, for some $a, b \in V$. It is easy to check that $(a, b)
= (3, 2)$, $(3, z)$, $(v, z)$ or $(z, 3)$. The set of known faces
is invariant under the map $(0, 9)(1, 6)(2, v)(3, z)(5, 8)(7,
u)$. So, we may assume that $(a, b) = (3, 2)$, $(3, z)$ or $(z,
3)$. If $(a, b) = (3, 2)$ then, considering the links of $2$ and
$3$, we get $7$ vertices in ${\rm lk}(7)$. If $(a, b) = (3, z)$
then, considering ${\rm lk}(3)$ we get $C_4(7, 3, 0, 1) \subseteq
{\rm lk}(2)$. So, ${\rm lk}(w) = C_6(u, 5, 8, 7, z, 3)$.

Completing successively, we get ${\rm lk}(3) = C_6(2, 0, 4, z, w,
u)$, ${\rm lk}(u) = C_6(2, 3, w, 5, 6, v)$, ${\rm lk}(2) = C_6(1,
0, 3, u, v, 7)$ and ${\rm lk}(z) = C_6(3, 4, 9, v, 7, w)$. Here $M
\cong Q_{7,2}$ by the map $\varphi \circ (0, 1, w)(2, z, 4, 3)(5,
9)(6, 7, v)(8, u)$.

\smallskip

\noindent {\bf Subcase 3.2.3.} ${\rm lk}(5) = C_6(u, 6, 0, 4, w,
z)$. This implies that ${\rm lk}(4) = C_6(w, 5, 0, 3, x, y)$, for
some $x, y \in V$. It is easy to check that $(x, y) = (7, 2)$,
$(7, 8)$, $(7, v)$, $(8, 7)$, $(8, 9)$, $(9, v)$, $(u, v)$, $(v,
7)$, $(v, 8)$, $(v, 9)$, $(v, u)$, $(z, 7)$, $(z, 8)$, $(z, u)$
or $(z, v)$. By similar arguments as in the previous claims one
gets $(x, y) = (8, 7)$, $(8, 9)$ or $(9, v)$.

\smallskip

\noindent {\bf Subcase 3.2.3.1.} ${\rm lk}(4) = C_6(w, 5, 0, 3, 8,
9)$. Completing successively, we get ${\rm lk}(9) = C_6(w, 4, 8,
1, 6, v)$, ${\rm lk}(w) = C_6(v, 9, 4, 5, z, 7)$, ${\rm lk}(8) =
C_6(3, 4, 9, 1, 7, z)$, ${\rm lk}(7) = C_6(2, 1, 8, z$, $w, v)$,
${\rm lk}(z) = C_6(3, 8, 7, w$, $5, u)$, ${\rm lk}(v) = C_6(2, 7,
w, 9, 6, u)$ and ${\rm lk}(u) = C_6(2, 3, z, 5, 6, v)$. Here $M
\cong T_{14, 1, 3}$ by the map $\varphi \circ (0, 5, 1, 6, 2, 9,
3, 8, 7)(u, v, w, z)$.

\smallskip

\noindent {\bf Subcase 3.2.3.2.} ${\rm lk}(4) = C_6(w, 5, 0, 3, 9,
v)$. Completing successively, we get ${\rm lk}(9) = C_6(3, 4, v,
6, 1, 8)$, ${\rm lk}(3) = C_6(2, 0, 4, 9, 8, z)$, ${\rm lk}(v) =
C_6(u, 6, 9, 4, w, 7)$, ${\rm lk}(2) = C_6(1, 0, 3, z$, $u, 7)$,
${\rm lk}(8) = C_7(1, 9, 3, z$, $w, 7)$, ${\rm lk}(z) = C_6(2, 3,
8, w, 5, u)$, ${\rm lk}(w) = C_6(4, 5, z, 8, 7, v)$ and ${\rm
lk}(7) = C_6(1, 2, u, v, w, 8)$. Here $M \cong T_{14, 1, 3}$ by
the map $\varphi \circ (0, 1, 5, z, w)(3, v, 7, 6, 4$, $u)(8, 9)$.

\smallskip

\noindent {\bf Subcase 3.2.3.3.} ${\rm lk}(4) = C_6(w, 5, 0, 3, 8,
7)$. Completing successively we get ${\rm lk}(7) = C_6(w, 4, 8, 1,
2, v)$, ${\rm lk}(8) = C_6(3, 4, 7, 1, 9, z)$, ${\rm lk}(9) =
C_6(z, 8, 1, 6, v, w)$, ${\rm lk}(v) = C_6(2, 7, w, 9$, $6, u)$,
${\rm lk}(z) = C_6(5, w, 9, 8$, $3, u)$, ${\rm lk}(u) = C_6(2, v,
6, 5, z, 3)$ and ${\rm lk}(2) = C_6(1, 0, 3, u, v, 7)$. Here $M
\cong Q_{7, 2}$ by the map given by $\varphi \circ (0, 1, 3, 7,
u, 8, 5, z, 6, 2, 9, 4, w, v)$. \hfill $\Box$

\begin{lemma} $\!\!${\bf .}  \label{tkb15v} If $M$ is a $15$-vertex
degree-regular combinatorial $2$-manifolds of Euler
characteristic $0$ then $M$ is isomorphic to $T_{15, 1, 2},
\dots, T_{15, 1, 5}$, $B_{3, 5}$, $B_{5, 3}$ or $Q_{5, 3}$.
\end{lemma}

\noindent {\bf Proof.} Let $M$ be a 15-vertex degree regular
combinatorial 2-manifold of Euler characteristic 0. Let the
vertex set $V$ be $\{0, 1, \dots, 9, u, v, w, z, s\}$. Let
$\varphi \colon V \to \{1, \dots, 15\}$ be given by $\varphi(i) =
i$, for $1 \leq i \leq 9$, $\varphi(0) = 10$, $\varphi(u) = 11$,
$\varphi(v) = 12$, $\varphi(w) = 13$, $\varphi(z) = 14$ and
$\varphi(s) = 15$.

Since $\chi(M) = 0$, the degree of each vertex is 6. As earlier,
we may assume that ${\rm lk}(0) = C_6(1, 2, 3, 4, 5, 6)$. By Lemma
\ref{l3.8}, ${\rm lk}(1) = C_6(7, 2, 0, 6, x, y)$, for some $x, y
\in V$. It is easy to see that $(x, y) = (3, 4)$, $(3, 5)$, $(3,
8)$, $(4, 3)$, $(4, 5)$, $(4, 8)$, $(8, 3)$, $(8, 4)$, $(8, 5)$,
$(8, 9)$.

If $(x, y) = (3, 5)$ then, considering ${\rm lk}(3)$ we get $7$
vertices in ${\rm lk}(5)$. The case $(x, y) = (8, 3)$ is
isomorphic to the case $(x, y) = (4, 8)$ by the map $(0, 1) (2,
6)(3, 4, 8)(5, 7)$ and the case $(x, y) = (8, 5)$ is isomorphic to
the case $(x, y) = (8, 3)$ by the map $(2, 6)(3, 5)(7, 8)$. So, we
may assume that $(x, y) = (3, 4)$, $(3, 8)$, $(4, 3)$, $(4, 5)$,
$(4, 8)$, $(8, 4)$ or $(8, 9)$.

\smallskip

\noindent {\bf Case 1.} $(x, y) = (3, 4)$, i.e., ${\rm lk}(1) =
C_6(7, 2, 0, 6, 3, 4)$. Then ${\rm lk}(3) = C_6(2, 0, 4, 1, 6,
8)$, ${\rm lk}(6) = C_6(5, 0, 1, 3, 8, 9)$, ${\rm lk}(4) = C_6(5,
0, 3, 1, 7, u)$ and ${\rm lk}(2) = C_6(8, 3, 0, 1, 7, v)$. Now, it
is easy to see that ${\rm lk}(8) = C_6(9, 6, 3, 2, v, u)$ or
$C_6(9, 6, 3, 2, v, w)$. In the first case, we get $34$ edges in
$M[\{0, \dots, 9, u, v\}]$, a contradiction to Lemma \ref{l3.9}.
So, ${\rm lk}(8) = C_6(9, 6, 3, 2, v, w)$. Now, completing
successively, we get ${\rm lk}(7) = C_6(u, 4, 1, 2, v, z)$, ${\rm
lk}(v) = C_6(w, 8, 2, 7, z, s)$, ${\rm lk}(5) = C_6(4, 0, 6, 9,
s, u)$, ${\rm lk}(s) = C_6(5, 9, z, v, w, u)$, ${\rm lk}(u) =
C_6(7, 4, 5, s, w, z)$ and ${\rm lk}(9) = C_6(5, 6, 8, w, z, s)$.
Here $M \cong B_{3, 5}$ by the map $\psi_{35} \circ (2, 9, z)(0,
v, 6, u, 1, 8)(3, 7, 5, s)$, where $\psi_{35} \colon V \to
V(B_{3, 5})$ given by $\psi_{35}(i) = v_{1i}$, $\psi_{35}(3+i) =
v_{2i}$, $\psi_{35}(6+i) = v_{3i}$, $1 \leq i \leq 3$,
$\psi_{35}(0) = v_{41}$, $\psi_{35}(u) = v_{42}$, $\psi_{35}(v) =
v_{43}$, $\psi_{35}(w) = v_{51}$, $\psi_{35}(z) = v_{52}$ and
$\psi_{35}(s) = v_{53}$.

\smallskip

\noindent {\bf Case 2.} $(x, y) = (3, 8)$. Then ${\rm lk}(3) =
C_6(2, 0, 4, 8, 1, 6)$ or $C_6(2, 0, 4, 6, 1, 8)$.

\smallskip

\noindent {\bf Subcase 2.1.} ${\rm lk}(3) = C_6(2, 0, 4, 8, 1,
6)$. Completing successively, we get ${\rm lk}(2) = C_6(7, 1, 0$,
$3, 6, 9)$, ${\rm lk}(6) = C_6(5, 0, 1, 3, 2, 9)$, ${\rm lk}(9) =
C_6(5, 6, 2, 7, v, u)$, ${\rm lk}(5) = C_6(4, 0, 6, 9, u, w)$,
${\rm lk}(7) = C_6(8, 1, 2, 9, v, z)$, ${\rm lk}(8) = C_6(4, 3,
1, 7, z, s)$, ${\rm lk}(4) = C_6(5, 0, 3, 8, s, w)$, ${\rm lk}(u)
= C_6(5, 9, v, s, z, w)$, ${\rm lk}(s) = C_6(4, 8, z, u, v, w)$
and ${\rm lk}(z) = C_6(7, 8, s, u, w, v)$. Here $M$ is isomorphic
to $Q_{5, 3}$ by the map  $\psi \circ (5, 8, s, v, 9, 7, 6)$ $(0,
1, 3, 2, 4, z, w, u)$, where $\psi \colon V \to V(Q_{5, 3})$ is
given by $\psi(i) = u_{i1}$, $\psi(5+i) = u_{i2}$, $1 \leq i \leq
3$, $\psi(3+j) = v_{j1}$, $1 \leq j \leq 2$, $\psi(9) = v_{12}$,
$\psi(0) = v_{22}$, $\psi(u) = u_{13}$, $\psi(v) = u_{23}$,
$\psi(w) = u_{33}$, $\psi(z) = v_{13}$ and $\psi(s) = v_{23}$.

\smallskip

\noindent {\bf Subcase 2.2.} ${\rm lk}(3) = C_6(2, 0, 4, 6, 1,
8)$. Completing successively, we get ${\rm lk}(2)$ $ = C_6(7, 1,
0$, $3, 8, 9)$, ${\rm lk}(8) = C_6(7, 1, 3, 2, 9, u)$, ${\rm
lk}(6) = C_6(5, 0, 1, 3, 4, v)$, ${\rm lk}(4) = C_6(5, 0, 3, 6,
v, w)$, ${\rm lk}(7) = C_6(9, 2, 1, 8, u, z)$, ${\rm lk}(5) =
C_6(v, 6, 0, 4, w, s)$, ${\rm lk}(9) = C_6(8, 2, 7, z, s, u)$,
${\rm lk}(v) = C_6(6, 4, w, z, s, 5)$, ${\rm lk}(s) = C_6(u, w,
5, v, z, 9)$ and ${\rm lk}(u) = C_6(7, 8, 9, s, w, z)$. Here $M
\cong T_{15,1,2}$ by the map $\varphi \circ (0, 7, 2, 4, 9, 1$,
$5)(3, 6, 8)(u, s, w, v)$.

\smallskip

\noindent {\bf Case 3.} $(x, y) = (4, 3)$, i.e., ${\rm lk}(1) =
C_6(7, 2, 0, 6, 4, 3)$. Then ${\rm lk}(4) = C_6(5, 0, 3, 1, 6,
8)$, ${\rm lk}(3) = C_6(2, 0, 4, 1, 7, 9)$, ${\rm lk}(6) = C_6(5,
0, 1, 4, 8, u)$, ${\rm lk}(2) = C_6(9, 3, 0, 1, 7, v)$, ${\rm
lk}(7) = C_6(9, 3, 1, 2, v, w)$, ${\rm lk}(5) = C_6(8, 4, 0, 6,
u, z)$, ${\rm lk}(8) = C_6(u, 6, 4, 5, z, s)$, ${\rm lk} (9) =
C_6(v, 2, 3, 7$, $w, s)$. Thus, ${\rm lk}(s) = C_6(v, 9, w, z, 8,
u)$ or $C_6(v, 9, w, u, 8, z)$. In the first case, considering
the links of $s$, $u$, $v$ successively, we get $7$ vertices in
${\rm lk}(v)$. So, ${\rm lk}(s) = C_6(v, 9, w, u, 8, z)$. Now,
completing successively, we get ${\rm lk}(u) = C_6(5, 6, 8, s, w,
z)$ and ${\rm lk}(z) = C_6(5, 8, s, v, w, u)$. Here $M\cong
T_{15,1,2}$ by the map $\varphi \circ (0, 5, 2, 8, 1, 6, 3, 7,
9)(u, s, w, v)$.

\smallskip

\noindent {\bf Case 4.} $(x, y) = (4, 5)$, i.e., ${\rm lk}(1) =
C_6(7, 2, 0, 6, 4, 5)$. Now, completing successively, we get
${\rm lk}(4) = C_6(3, 0, 5, 1, 6, 8)$, ${\rm lk}(6) = C_6(5, 0,
1, 4, 8, 9)$, ${\rm lk}(5) = C_6(7, 1, 4, 0, 6, 9)$, ${\rm lk}(7)
= C_6(2, 1, 5, 9, v, u)$, ${\rm lk}(2) = C_6(3, 0, 1, 7, u, w)$,
${\rm lk}(3) = C_6(8, 4, 0, 2, w, z)$, ${\rm lk}(8) = C_6(9, 6,
4, 3$, $z, s)$, ${\rm lk}(9) = C_6(7, 5, 6, 8, s, v)$, ${\rm
lk}(v) = C_6(u, 7, 9, s$, $w, z)$, ${\rm lk}(z) = C_6(w, 3, 8, s,
u, v)$ and ${\rm lk} (u) = C_6(2, 7, v, z, s, w)$. Now, $M \cong
Q_{5, 3}$ by the map $\psi \circ (0, 1, 2, z)(3, 8, 7, s, 9, 6,
4, 5)(u, v$, $w)$, where $\psi$ is as in Subcase 2.1.

\smallskip

\noindent {\bf Case 5.} $(x, y) = (4, 8)$, i.e., ${\rm lk}(1) =
C_6(7, 2, 0, 6, 4, 8)$. Now, completing successively, we get
${\rm lk}(4) = C_6(5, 0, 3, 6, 1, 8)$, ${\rm lk}(6) = C_6(5, 0,
1, 4, 3, 9)$, ${\rm lk}(3) = C_6(2, 0, 4, 6, 9, u)$, ${\rm lk}(2)
= C_6(7, 1, 0, 3, u, v)$, ${\rm lk}(5) = C_6(8, 4, 0, 6, 9, w)$,
${\rm lk}(8) = C_6(7, 1, 4, 5, w, z)$, ${\rm lk}(9) = C_6(u, 3,
6, 5$, $w, s)$, ${\rm lk}(7) = C_6(v, 2, 1, 8, z, s)$, ${\rm
lk}(s) = C_6(w, 9, u, z, 7, v)$, ${\rm lk}(w) = C_6(8, 5, 9, s,
v, z)$ and ${\rm lk}(z) = C_6(7, 8, w, v, u, s)$. Now, $M$ is
isomorphic to $B_{3, 5}$ by the map $\psi_{35}\circ(1, u)(2,
s)(5, 9)(0$, $v, 3, 4, 7, z, w, 6, 8)$, where $\psi_{35}$ is as
in Case 1.

\smallskip

\noindent {\bf Case 6.} $(x, y) = (8, 4)$, i.e., ${\rm lk}(1) =
C_6(7, 2, 0, 6, 8, 4)$. Now, ${\rm lk}(4) = C_6(5, 0, 3, 8, 1,
7)$ or $C_6(5, 0, 3, 7, 1, 8)$. In the first case, we get ${\rm
lk}(3) = C_6(8, 4, 0, 2, a, b)$, for some $a, b \in V$. Using
Lemma \ref{l3.9}, we may assume that $(a, b) = (9, u)$.
Considering the links of $2$, $7$, $5$, $6$, $8$, $9$, $s$, $u$
and $v$ successively, we get $7$ vertices in ${\rm lk}(v)$. Thus,
${\rm lk}(4) = C_6(5, 0, 3, 7, 1, 8)$.

Again, by using Lemma \ref{l3.9}, we get ${\rm lk}(3) = C_6(7, 4,
0, 2, 9, u)$, ${\rm lk}(2) = C_6(7, 1, 0, 3, 9, v)$, ${\rm lk}(7)
= C_6(2, 1, 4, 3, u, v)$, ${\rm lk}(9) = C_6(2, 3, u, z, w, v)$,
${\rm lk}(v) = C_6(u, 7, 2, 9, w, s)$, ${\rm lk}(u) = C_6(9, 3, 7,
v, s, z)$. Then ${\rm lk}(z) = C_6(s, u, 9, w, a, b)$, for some
$a, b \in V$. It is easy to see that $(a, b) = (5, 6)$, $(5, 8)$,
$(6, 5)$, $(6, 8)$, $(8, 5)$ or $(8, 6)$. Since the set of known
faces remain invariant under the map $(0, 4)(2, 7)(6, 8)(9, u)(w,
s)$, we can assume that $(a, b) = (5, 6)$, $(5, 8)$, $(6, 8)$ or
$(8, 6)$.

\smallskip

\noindent {\bf Subcase 6.1.} $(a, b) = (5, 6)$. Completing
successively, we get ${\rm lk}(5) = C_6(8, 4, 0, 6, z, w)$, ${\rm
lk}(6) = C_6(8, 1, 0, 5, z, s) $ and ${\rm lk}(8) = C_6(5, 4, 1,
6, s, w)$. Now, $M$ is isomorphic to $B_{5, 3}$ by the map
$\psi_{53} \circ (0, s, 3, 1)(2, 6, 9, 7, u)(4, 5, z, 8)$, where
$\psi_{53}\colon V\to V(B_{5, 3})$ is given by $\psi_{53}(i) =
v_{1i}$, for $1 \leq i \leq 5$, $\psi_{53}(5+i) = v_{2i}$, for $1
\leq i \leq 4$, $\psi_{53}(0) = v_{25}$, $\psi_{53}(u) = v_{31}$,
$\psi_{53}(v) = v_{32}$, $\psi_{53}(w) = v_{33}$, $\psi_{53}(z) =
v_{34}$, $\psi_{53}(s) = v_{35}$.

\smallskip

\noindent {\bf Subcase 6.2.} $(a, b) = (5, 8)$. Completing
successively, we get ${\rm lk}(5) = C_6(w, z, 8, 4, 0, 6)$, ${\rm
lk}(8) = C_6(6, 1, 4, 5, z, s)$ and ${\rm lk}(6) = C_6(5, 0, 1,
8, s, w)$. Here $M$ is isomorphic to $T_{15, 1, 5}$ by the map
$\varphi \circ (0, 2, 1, 7, 6, 8, w, 9)(3, u, 5)(4, v, s, z)$.

\smallskip

\noindent {\bf Subcase 6.3.} $(a, b) = (6, 8)$. Completing
successively, we get ${\rm lk}(6) = C_6(5, 0, 1, 8, z$, $w)$,
${\rm lk}(8) = C_6(1, 4, 5, s, z, 6)$ and ${\rm lk}(5) = 4, 0, 6,
w, s, 8)$. Here $M \cong B_{5, 3}$ by the map $\psi_{53} \circ (1,
9)(0, z, 2, s, v, 6, 8, 3, 5, w, 7)$, where $\psi_{53}$ is as in
Subcase 6.1.

\smallskip

\noindent {\bf Subcase 6.4.} $(a, b) = (8, 6)$. Completing
successively, we get ${\rm lk}(6) = C_6(5, 0, 1, 8, z, s)$, ${\rm
lk}(8) = C_6(5, 4, 1, 6, z, w)$ and ${\rm lk}(5) = C_6(8, 4, 0,
6, s, w)$. Here $M$ is isomorphic to $T_{15,1,4}$ by the map
given by $\varphi \circ (0, 7, 1, 2, 6, 3, u, s, 4, v, 5, 8, w,
9)$.

\smallskip

\noindent {\bf Case 7.} $(x, y) = (8, 9)$, i.e., ${\rm lk}(1) =
C_6(7, 2, 0, 6, 8, 9)$. Now, ${\rm lk}(6) = C_6(8, 1, 0, 5, a,
b)$. It is easy to check that $(a, b) = (2, 3)$, $(2, 7)$, $(3,
2)$, $(3, 7)$, $(3, 4)$, $(3, u)$, $(7, 2)$, $(7, 3)$, $(7, 4)$,
$(7, u)$, $(9, 3)$, $(9, 4)$, $(9, 7)$, $(9, u)$, $(u, 3)$, $(u,
4)$, $(u, 7)$, $(u, v)$.

If $(a, b) = (3, 7)$ then, considering links of $3$ and $7$, we
get $7$ vertices in ${\rm lk}(7)$. If $(a, b) = (7, 3)$, $(7, 4)$
or $(9, 4)$ then, considering ${\rm lk}(b)$, we get $7$ vertices
in ${\rm lk}(a)$. Since the set of known faces remain invariant
under that map $(0, 1)(3, 7)(4, 9)(5, 8)$, we may assume that
$(a, b) = (2, 3)$, $(2, 7)$, $(3, 4)$, $(3, u)$, $(7, u)$, $(9,
u)$, $(u, v)$.

\smallskip

\noindent {\bf Subcase 7.1.} ${\rm lk}(6) = C_6(8, 1, 0, 5, 2,
3)$. Completing successively, we get ${\rm lk}(2) = C_6(1, 0, 3$,
$6, 5, 7)$, ${\rm lk}(5) = C_6(4, 0, 6, 2, 7, u)$, ${\rm lk}(3) =
C_6(4, 0, 2, 6, 8, v)$, ${\rm lk}(8) = C_6(9, 1, 6, 3, v, w)$,
${\rm lk}(4) = C_6(u, 5, 0, 3, v, z)$, ${\rm lk}(v) = C_6(w, 8,
3, 4, z, s)$, ${\rm lk}(7) = C_6(9, 1, 2, 5, u, s)$, ${\rm lk}(9)
= C_6(w, 8, 1, 7$, $s, z)$, ${\rm lk}(z) = C_6(4, v, s, 9, w, u)$
and ${\rm lk}(s) = C_6(w, u, 7, 9, z, v)$. Here $M \cong B_{3,
5}$ by the map $\psi_{35} \circ (0, 7, 5, 9, 1, 4)(2, 8, s)(3, u,
6, v, z, w)$, where $\psi_{35}$ is as in Case 1.

\smallskip

\noindent {\bf Subcase 7.2.} ${\rm lk}(6) = C_6(8, 1, 0, 5, 2,
7)$. Completing successively, we get ${\rm lk}(2) = C_6(5, 6, 7$,
$1, 0, 3)$, ${\rm lk}(5) = C_6(4, 0, 6, 2, 3, u)$, ${\rm lk}(3) =
C_6(4, 0, 2, 5, u, v)$, ${\rm lk}(4) = C_6(u, 5, 0, 3, v, w)$,
${\rm lk}(u) = C_6(v, 3, 5, 4, w, z)$, ${\rm lk}(7) = C_6(9, 1,
2, 6, 8, s)$, ${\rm lk}(v) = C_6(w, 4, 3, u, z, s)$, ${\rm lk}(8)
= C_6(9, 1, 6, 7$, $s, z)$, ${\rm lk}(9) = C_6(8, 1, 7, s, w, z)$
and ${\rm lk}(s) = C_6(8, 7, 9, w$, $v, z)$. Here $M \cong
T_{15,1,2}$ by the map $\varphi \circ (0, 7, 3, 9, 1, 4)(2, 6, 5,
8)$.

\smallskip

\noindent {\bf Subcase 7.3.} ${\rm lk}(6) = C_6(8, 1, 0, 5, 3,
4)$. Completing successively, we get ${\rm lk}(3) = C_6(2, 0, 4$,
$6, 5, u)$, ${\rm lk}(4) = C_6(5, 0, 3, 6, 8, v)$, ${\rm lk}(5) =
C_6(3, 6, 0, 4, v, u)$, ${\rm lk}(2) = C_6(7, 1, 0, 3, u, w)$,
${\rm lk}(u) = C_6(v, 5, 3, 2, w, z)$, ${\rm lk}(v) = C_6(8, 4,
5, u, z, s)$, ${\rm lk}(8) = C_6(9, 1, 6, 4, v, s)$, ${\rm lk}(9)
= C_6(7, 1, 8, s$, $w, z)$, ${\rm lk}(w) = C_6(2, 7, s, 9, z, u)$
and ${\rm lk}(7) = C_6(9, 1, 2, w$, $s, z)$. Here $M\cong Q_{5,
3}$ by the map $\psi \circ (0, 1, 8, 7, u, s, 9)(2, z, w, v, 6,
5, 3)$, where $\psi$ is as in Subcase 2.1.

\smallskip

\noindent {\bf Subcase 7.4.} ${\rm lk}(6) = C_6(8, 1, 0, 5, 3,
u)$. Completing successively, we get ${\rm lk}(3) = C_6(4, 0, 2$,
$5, 6, u)$, ${\rm lk}(2) = C_6(7, 1, 0, 3, 5, v)$, ${\rm lk}(5) =
C_6(4, 0, 6, 3, 2, v)$, ${\rm lk}(4) = C_6(u, 3, 0, 5, v, w)$,
${\rm lk}(u) = C_6(8, 6, 3, 4, w$, $z)$, ${\rm lk}(8) = C_6(9, 1,
6, u, z, s)$, ${\rm lk}(v) = C_6(7, 2, 5, 4, w, s)$, ${\rm lk}(7)
= C_6(9, 1, 2, v$, $s, z)$, ${\rm lk}(s) = C_6(9, w, v, 7, z, 8)$
and ${\rm lk}(z) = C_6(w, 9, 7, s$, $8, u)$. Here $M \cong B_{3,
5}$ by the map $\psi_{35} \circ (0, 7, z, 2, u, 5, v, s, 3, 8, 6,
9, w, 1)$, where $\psi_{35}$ is as defined in Case 1.

\smallskip

\noindent {\bf Subcase 7.5.} ${\rm lk}(6) = C_6(8, 1, 0, 5, 7,
u)$. Now, it is easy to see that ${\rm lk}(7) = C_6(9, 1, 2, 5, 6,
u)$ or $C_6(9, 1, 2, u, 6, 5)$. The first case is isomorphic to
Subcase 6.2 by that map $(0, 3, 9, 6, 4, u, 8$, $5, 7, 1)$. The
second case is isomorphic to Subcase 6.1 by the map $(0, 3, 9, 6,
4, u, 7, 1)(5, 8)$.

\smallskip

\noindent {\bf Subcase 7.6.} ${\rm lk}(6) = C_6(8, 1, 0, 5, 9,
u)$. Completing successively, we get ${\rm lk}(9) = C_6(8, 1, 7$,
$u, 6, 5)$, ${\rm lk}(5) = C_6(4, 0, 6, 9, 8, v)$, ${\rm lk}(8) =
C_6(6, 1, 9, 5, v, u)$, ${\rm lk}(u) = C_6(7, 9, 6, 8, v, w)$,
${\rm lk}(7) = C_6(2, 1, 9, u, w, z)$, ${\rm lk}(v) = C_6(4, 5,
8, u, w, s)$, ${\rm lk}(2) = C_6(3, 0, 1, 7, z, s)$, ${\rm lk}(s)
= C_6(3, 2, z, 4$, $v, w)$, ${\rm lk}(4) = C_6(3, 0, 5, v, s, z)$
and ${\rm lk}(3) = C_6(4, 0, 2, s, w, z)$. Here $M\cong B_{3, 5}$
by the map $\psi_{35} \circ (1, z, v, 5, 6, w, 4, 9, 3, 7, s, 8,
2, u)$, where $\psi_{35}$ is as in Case 1.

\smallskip

\noindent {\bf Subcase 7.7.} ${\rm lk}(6) = C_6(8, 1, 0, 5, u,
v)$. Then, ${\rm lk}(5) = C_6(u, 6, 0, 4, c, d)$, for some $c, d
\in V$. It is easy to see that $(c, d) = (2, 3)$, $(2, 7)$, $(7,
2)$, $(7, 3)$, $(7, 9)$, $(7, w)$, $(8, 9)$, $(9, 3)$, $(9, 7)$,
$(9, 8)$, $(9, w)$, $(v, 3)$, $(v, 7)$, $(v, 8)$, $(v, 9)$, $(v,
w)$, $(w, 3)$, $(w, 7)$, $(w, 9)$, $(w, z)$. Since the set of
known faces remain invariant under the map $(0, 6)(2, 8)(3, v)(4,
u)(7, 9)$, we may assume that $(c, d) = (2, 3)$, $(2, 7)$, $(7,
2)$, $(7, 3)$, $(7, 9)$, $(7, w)$, $(9, 3)$, $(9, 7)$, $(9, w)$,
$(v, 3)$, $(v, w)$, $(w, z)$.

\smallskip

\noindent {\bf Claim.} $(c, d) = (2, 3)$, $(9, w)$, $(v, w)$ or
$(w, z)$.

If $(c, d) = (2, 7)$ then, considering ${\rm lk}(2)$ we get
$C_4(3, 2, 5, 0) \subseteq {\rm lk}(4)$. If $(c, d) = (7, 3)$
then, considering ${\rm lk}(3)$, we get $C_4(7, 3, 0, 1) \subseteq
{\rm lk}(2)$.

If $(c, d) = (7, 2)$ then, considering the links of $2$, $u$, $3$,
$4$, $7$, $w$, $9$, $v$ successively, we get $7$ vertices in ${\rm
lk}(8)$. If $(c, d) = (9, 3)$ or $(v, 3)$ then, considering ${\rm
lk}(3)$, we get $7$ vertices in ${\rm lk}(c)$ . If $(c, d) = (9,
7)$ then, considering the links of $9$, $4$, $7$, $2$, $v$, $8$,
$s$ successively, we get $7$ vertices in ${\rm lk}(u)$. If $(c, d)
= (7, 9)$ then, considering the links of $5$, $7$, $2$, $4$, $3$,
$v$ successively, we get ${\rm lk}(u) = C_6(9, 5, 6, v, z, x)$,
where $x = z$ or $s$. In either case ${\rm lk}(x)$ has $\geq 7$
vertices.

If $(c, d) = (7, w)$ then, ${\rm lk}(7) = C_6(2, 1, 9, w, 5, 4)$
or $C_6(2, 1, 9, 4, 5, w)$. In the first case, considering the
links of $4$ and $2$ we get $C_4(z, 2, 0, 4) \subseteq {\rm
lk}(3)$. In the second case, considering links of $7$, $4$, $2$,
$w$, $z$, $u$ successively, we get $7$ vertices in ${\rm lk}(3)$.
This proves the claim.

\smallskip

\noindent {\bf Subcase 7.7.1.} ${\rm lk}(5) = C_6(u, 6, 0, 4, 2,
3)$. Now, completing successively we get ${\rm lk}(2) = C_6(1, 0,
3, 5, 4, 7)$, ${\rm lk}(4) = C_6(3, 0, 5, 2, 7, w)$, ${\rm lk}(3)
= C_6(5, 2, 0, 4, w, u)$, ${\rm lk}(u) = C_6(6, 5, 3, w$, $s, v)$,
${\rm lk}(7) = C_6(9, 1, 2, 4, w, z)$, ${\rm lk}(w) = C_6(3, 4,
7, z, s, u)$, ${\rm lk}(8) = C_6(9, 1, 6, v, z, s)$, ${\rm lk}(9)
= C_6(8, 1, 7, z, v, s)$ and ${\rm lk}(v) = C_6(8, 6, u, s, 9,
z)$. Here $M$ is isomorphic to $Q_{5, 3}$ by the map $\psi \circ
(0, 1, z, w, 6, 8, u, 7, s, 9, v)(3, 4)$, where $\psi$ is as in
Subcase 2.1.

\smallskip

\noindent {\bf Subcase 7.7.2.} ${\rm lk}(5) = C_6(u, 6, 0, 4, 9,
w)$. This implies that ${\rm lk}(9) = C_6(7, 1, 8, 4, 5, w)$ or
$C_6(7, 1, 8, w, 5, 4)$. In the first case, considering links of
$9$, $4$, $8$, $v$ successively we see that ${\rm lk}(v)$ can not
be a 6-cycle. Thus ${\rm lk}(9) = C_6(7, 1, 8, w, 5, 4)$.

Now, completing successively, we get ${\rm lk}(4) = C_6(3, 0, 5,
9, 7, z)$, ${\rm lk}(7) = C_6(2, 1, 9, 4, z, s)$, ${\rm lk}(2) =
C_6(3, 0, 1, 7, s, v)$, ${\rm lk}(3) = C_6(z, 4, 0, 2, v, u)$,
${\rm lk}(u) = C_6(5, 6, v, 3, z, w)$, ${\rm lk}(v) = C_6(6, u, 3,
2, s, 8)$, ${\rm lk}(8) = C_6(9, 1, 6, v, s, w)$ and ${\rm lk}(z)
= C_6(7, 4, 3, u, w, s)$. Here $M \cong T_{15, 1, 3}$ by the map
$\varphi \circ (0, 5, 8)(1, 6, 9$, $7, 3)(u, v, w)(z, s)$.

\smallskip

\noindent {\bf Subcase 7.7.3.} ${\rm lk}(5) = C_6(u, 6, 0, 4, v,
w)$. Now, completing successively we get ${\rm lk}(v) = C_6(8, 6,
u, 4, 5, w)$, ${\rm lk}(4) = C_6(3, 0, 5, v, u, z)$, ${\rm lk}(u)
= C_6(5, 6, v, 4, z, w)$, ${\rm lk}(w) = C_6(8, v, 5, u$, $z, s)$,
${\rm lk}(8) = C_6(9, 1, 6, v, w, s)$, ${\rm lk}(z) = C_6(3, 4,
u, w, s, 7)$, ${\rm lk}(7) = C_6(2, 1, 9, 3, z, s)$, ${\rm lk}(3)
= C_6(2, 0, 4, z, 7, 9)$ and ${\rm lk}(2) = C_6(7, 1, 0, 3, 9,
s)$. Here $M$ is isomorphic to $B_{3, 5}$ by the map $\psi_{35}
\circ (1, 9, 8, 5, w)(2, 7, v)(3, u)(4, z, s)$, where $\psi_{35}$
is as in Case 1.

\smallskip

\noindent {\bf Subcase 7.7.4.} ${\rm lk}(5) = C_6(u, 6, 0, 4, w,
z)$. Then  ${\rm lk}(4) = C_6(w, 5, 0, 3, x, y)$, for some $x, y
\in V$. It is easy to see that $(x, y) = (7, 2)$, $(7, 9)$, $(7,
v)$, $(7, s)$, $(8, 9)$, $(8, v)$, $(9, 8)$, $(9, v)$, $(9, s)$,
$(u, v)$, $(v, 8)$, $(v, 9)$, $(v, u)$, $(v, s)$, $(z, 9)$, $(z,
u)$, $(z, v)$, $(z, s)$, $(s, 9)$, $(s, v)$. By the similar
arguments as before one gets $(x, y) = (7, 2)$, $(9, 8)$ or $(z,
u)$.

\smallskip

\noindent {\bf Subcase 7.7.4.1.} ${\rm lk}(4) = C_6(w, 5, 0, 3, 9,
8)$. Completing successively, we get ${\rm lk}(8)$ $ = C_6(6, 1,
9, 4, w, v)$, ${\rm lk}(9) = C_6(7, 1, 8, 4, 3, s)$, ${\rm lk}(3)
= C_6(2, 0, 4, 9, s, z)$, ${\rm lk}(2) = C_6(7, 1, 0, 3, z$,
$u)$, ${\rm lk}(u) = C_6(5, 6, v, 7, 2, z)$, ${\rm lk}(7) =
C_6(9, 1, 2, u, v, s)$, ${\rm lk}(v) = C_6(8, 6, u, 7, s, w)$ and
${\rm lk}(w) = C_6(5, 4, 8, v, s, z)$. Here $M \cong T_{15,1,3}$
by the map $\varphi \circ (0, 6, 9, 4, 7, 1, 5)(u, w)$.

\smallskip

\noindent {\bf Subcase 7.7.4.2.} ${\rm lk}(4) = C_6(w, 5, 0, 3, z,
u)$. Completing successively, we get ${\rm lk}(u) = C_6(6, 5, z,
4, w, v)$, ${\rm lk}(z) = C_6(w, 5, u, 4, 3, s)$, ${\rm lk}(w) =
C_6(u, 4, 5, z, s, v)$, ${\rm lk}(3) = C_6(2, 0, 4, z$, $s, 9)$,
${\rm lk}(9) = C_6(3, 2, 8, 1, 7, s)$, ${\rm lk}(s) = C_6(w, z,
3, 9, 7, v)$, ${\rm lk}(v) = C_6(6, u, w, s$, $7, 8)$ and ${\rm
lk}(7) = C_6(2, 1, 9, s, v, 8)$. Here $M \cong Q_{5, 3}$ by the
map $\psi \circ (0, z, w, 9, 3, s, 6, 8$, $5, u)(4, v, 7)$, where
$\psi$ is as in Subcase 2.1.

\smallskip

\noindent {\bf Subcase 7.7.4.3.} ${\rm lk}(4) = C_6(w, 5, 0, 3, 7,
2)$. Then, completing successively, we get ${\rm lk}(2) = C_6(3,
0, 1, 7, 4, w)$, ${\rm lk}(3) = C_6(7, 4, 0, 2, w, s)$, ${\rm
lk}(7) = C_6(1, 2, 4, 3, s, 9)$, ${\rm lk}(w) = C_6(5, 4, 2, 3$,
$s, z)$, ${\rm lk}(s) = C_6(9, 7, 3, w, z, v)$, ${\rm lk}(v) =
C_6(u, 6, 8, z, s, 9)$, ${\rm lk}(z) = C_6(5, w, s, v$, $8, u)$
and ${\rm lk}(8) = C_6(9, 1, 6, v, z, u)$. Here $M\cong B_{3, 5}$
by the map $\psi_{35} \circ (0, 9, 1, 5)(2$, $8)(3, v)(4$, $7)(u,
w)$, where $\psi_{35}$ is as in Case 1.  \hfill $\Box$

\bigskip

\noindent {\bf Proof of Theorem \ref{t6}.} Let $M$ be an
$n$-vertex degree-regular combinatorial 2-manifold of Euler
characteristic 0. Let $d$ be the degree of each vertex. Then $nd =
2f_1(M) = 3f_2(M)$ and $n - f_1(M) + f_2(M) = 0$. These imply that
$d =6$. Now, if $n\in \{12, 14, 15\}$ then, by Lemmas
\ref{tkb12v}, \dots, \ref{tkb15v}, $M$ is isomorphic to $T_{12,
1, 2}, \dots, T_{12, 1, 4}$, $T_{6, 2, 2}$, $T_{14, 1, 2}$,
$T_{14, 1, 3}$, $Q_{7, 2}$, $T_{15, 1, 2}, \dots, T_{15, 1, 5}$,
$Q_{5, 3}$, $B_{3, 4}$, $B_{4, 3}$, $B_{3, 5}$, $B_{5, 3}$ or
$K_{3, 4}$.

Since $B_{3, 4}$, $B_{4, 3}$, $B_{3, 5}$, $B_{5, 3}$, $Q_{7, 2}$,
$Q_{5, 3}$ and $K_{3, 4}$ are non-orientable and remaining 10 are
orientable, the second and third statements follow from Lemma
\ref{l2.1} (b), (d), (f), Lemma \ref{l2.3} (b) and Lemma
\ref{l2.5} (a), (b).

The last statement follows from the fact that $B_{m, n}$, $K_{m,
2k}$, $Q_{2k+1, n}$ are not weakly regular and $Q_{2k+1, 2}$ is
weakly regular for all $m, n\geq 3$ and $k\geq 2$. \hfill $\Box$

\bigskip

\noindent {\bf  Acknowledgement\,:} The authors are thankful to
B. Bagchi, S. P. Inamdar and N. S. N. Sastry for useful
conversations. The authors thank the anonymous referee for many
useful comments which led to substantial improvements in the
presentation of this paper. Theorem 1 was conjectured by the
referee. The second author thanks CSIR, New Delhi, India for its
research fellowship (Award No.: 9/79(797)/2001\,-\,EMR\,-\,I).

\end{document}